\documentclass[10pt]{article}

\usepackage[letterpaper, top=0.75in, bottom=0.75in, left=1in, right=1in]{geometry}

\usepackage{graphicx}
\usepackage{mathtools}
\usepackage{caption}
\usepackage{amsfonts} % For mathbb.
\usepackage{empheq} % For curly braces in equations (empheq command).
\usepackage{upgreek} % For non-italic and "boldable" greek letters.
\usepackage{physics} % For derivatives.
\usepackage{colortbl} % For cell colors in tables.
\usepackage[caption=false,font=scriptsize]{subfig} % For \subfloat command.
\usepackage[numbers]{natbib}
\usepackage{esint} % For oiint integrals.
\usepackage[shortlabels]{enumitem} % Change label of enumerate.
\usepackage{booktabs}  % For toprule.
\usepackage[colorlinks=true,linkcolor=red,anchorcolor=red,citecolor=blue,filecolor=red,urlcolor=blue,pdfauthor=author,hyperfootnotes=true]{hyperref} % Hyperlinks in footnotes doesn't work with footmisc package so I disabled them in hyperref settings.
\usepackage[bottom]{footmisc} % Put footnotes at bottom of page.

\makeatletter
\renewcommand*\env@matrix[1][\arraystretch]{%
	\edef\arraystretch{#1}%
	\hskip -\arraycolsep
	\let\@ifnextchar\new@ifnextchar
	\array{*\c@MaxMatrixCols c}}
\makeatother

\makeatletter
\def\env@dmatrix{\hskip -\arraycolsep
	\let\@ifnextchar\new@ifnextchar
	\extrarowheight=0ex
	\array{*\c@MaxMatrixCols{>{\displaystyle}c}}}

\makeatother

\newcommand{\cbox}[1]{{%
		\setlength{\fboxsep}{-\fboxrule}% Remove impact of space and rule inserted by box
		\fcolorbox{white}{yellow}{\strut#1}%
}}

\title{Topology Optimization of Fluid-Structure Interaction Problems with Total Stress Equilibrium}
\author{Mohamed Abdelhamid \\ mahamid@yorku.ca \\ Department of Mechanical Engineering, York University, \\ 4700 Keele St., Toronto, M3J 1P3, Ontario, Canada \and Aleksander Czekanski \\ alex.czekanski@lassonde.yorku.ca \\ Department of Mechanical Engineering, York University, \\ 4700 Keele St., Toronto, M3J 1P3, Ontario, Canada}
\date{May 2023}

\begin{document}
\maketitle

\section*{Abstract}
This work extends force coupling in the topology optimization of fluid-structure interaction problems from hydrostatic to total stresses through the inclusion of viscous stress components. The superconvergent patch recovery technique is implemented to remove the discontinuities in velocity derivatives over the finite elements boundaries. The sensitivity analysis is derived analytically for the superconvergent patch recovery approach and further verified through the use of the complex-step derivative approximation method. Numerical examples demonstrate a differentiation in the optimized designs using pressure vs. total stress coupling depending on the flow characteristics of the design problem.

\section*{Keywords} fluid-structure interactions, topology optimization, density-based methods, traction equilibrium, viscous stresses

\maketitle

\section{Introduction}
\label{sec:intro}

\textbf{F}luid-\textbf{s}tructure \textbf{i}nteractions (FSI) refer to problems where a deformable or a movable structure interacts with a flowing or stationary fluid. Examples of FSI in everyday life are abundant; the flow of blood through flexible arteries, the deflection of an elastic airplane wing under aerodynamic loads, and the vibration of bridges and tall buildings under wind loads \cite{Richter2017}. A large body of research work has been dedicated to the numerical analysis of these interactions \cite{Hou2012, Dowell2001}, and most multiphysics simulation packages - commercial and open source - include some degree of capability for solving these problems. However, the application of numerical design tools such as \textbf{t}opology \textbf{o}ptimization (TO) to these problems is still lagging behind. This is mainly due to the multidisciplinary nature of the combined \textbf{t}opology \textbf{o}ptimization of \textbf{f}luid-\textbf{s}tructure \textbf{i}nteraction (TOFSI) problem and its strong nonlinear behavior resulting in numerous stability and convergence issues.

Density-based TO methods first appeared in the seminal work by Bendsoe\cite{Bendsoe1989} which marked the start of modern TO in its current form. In density-based methods, each discretization unit (e.g. finite element) is interpolated between material and void throughout the optimization process. Penalization of a property of interest is used to force the final solution towards a discrete 0/1 state \cite{Abdelhamid2021}. Later, TO was extended to Stokes fluid flow in \ \cite{Borrvall2003} using a parameter to control permeability for design parametrization. In \ \cite{Gersborg-Hansen2005}, the authors extended the work to the Navier-Stokes equations and used the analogy of a 2D channel flow with varying thickness for design parametrization. They also recognized, along with  Evgrafov\cite{Evgrafov2005}, the similarity between this design parametrization and Brinkman equations \cite{Brinkman1947} of fluid flow in porous media \cite[p.~15]{Nield2017}. This approach - later termed \textit{Brinkman penalization} - became the de facto method for TO of fluid flow without much regard to its physical interpretation.

The first application of TO to FSI problems appeared in \ \cite{Maute2002}. The authors used a \textit{segregated} domain formulation to define the three-field FSI problem (i.e. fluid, solid, and mesh). In the segregated formulation, there is no overlap between the fluid and the solid computational domains thus only allowing for \textit{dry} TO (i.e. the fluid-structure interface is not optimized). This formulation was used in \ \cite{Maute2004} and extended to mechanism design in adaptive wings in \ \cite{Maute2006}. In \ \cite{Stanford2009, Stanford2009a}, the authors utilized the three-field, segregated domain formulation in the design of wing skeletons of micro air vehicles using a two-material formulation; wing skin or carbon fiber reinforcement. Although - technically speaking - this was a form of \textit{wet} TO, due to using a two-material formulation there were no fluid porous elements, and hence from a force coupling perspective it's similar to the dry TO described in \ \cite{Maute2002}. In \ \cite{Kreissl2010}, the authors performed TO on a flexible micro-fluidic device, however, they only considered the effect of structural deformations on the fluid flow and ignored the fluid forces on the structure. In \ \cite{Vu-Huu2018}, a polytree-based adaptive polygonal finite element method was used for TO of fluid-submerged breakwater barriers. However, they used a low fidelity model such that the initial fluid loading was used to calculate the forces without further updates based on the changes in the design domain. Hence the fluid forces were calculated on the fully solid structure and no porosity was incorporated in the fluid domain.

The first appearance of the \textit{unified} domain formulation was in \ \cite{Yoon2010} where the fluid and solid computational domains overlapped such that the whole design domain - including the fluid-structure interface - could be optimized, hence the term \textit{wet} TO. In \ \cite{Yoon2010}, the author extended the fluid-solid force coupling to the entire fluid domain including porous elements but only considered the hydrostatic component of the fluid stress tensor. The divergence theorem was utilized to transfer the integral from surface to volume. The same unified formulation and force coupling was later utilized in \ \cite{Yoon2014, Yoon2014a, Yoon2017} by the same author. More recently, in \ \cite{Lundgaard2018}, the authors revisited the density-based unified formulation with special emphasis on the multidisciplinary interaction between the fluid and the solid and its effects on the optimized designs. They also investigated different fluid flow conditions and demonstrated the significance of careful selection of interpolation functions and parameters to ensure convergence. Nonetheless, all these works were limited to hydrostatic traction equilibrium.

In the present work, we aim to extend the traction equilibrium condition to include the \textit{total} stress tensor such that viscous stresses are considered in addition to the previously-considered pressure stresses. The remainder of the paper is organized as follows; in \textbf{section \ref{sec:gov_eqns}}, we state the governing equations for the fluid-structure interaction problem. In \textbf{section \ref{sec:fin_elem_form}}, we detail the finite element formulations of the governing equations and the traction equilibrium condition as well as the relevant sensitivity analysis. In \textbf{section \ref{sec:tst_prblm}}, we present the test problem to be used in the following sections. In \textbf{section \ref{sec:total_traction_section}}, we discuss the implications of the discontinuities in velocity derivatives at the elemental boundaries and how to overcome this issue. In \textbf{section \ref{sec:numerical_exp}}, we demonstrate with numerical examples the effect of using pressure vs. total stress force coupling. Finally, in \textbf{section \ref{sec:cnclsns}}, we list the conclusions.

\section{Governing Equations of the Fluid-Structure Interaction Problem}
\label{sec:gov_eqns}

As further detailed in Sec. \ref{sec:fin_elem_form}, we limit our discussion to infinitesimal strains in this work. For large deformations, a mesh deformation technique needs to be implemented and the fluid governing equations need to be redefined in the deformed configuration, cf. \ \cite[p.~594]{Yoon2010} for such a derivation in a TOFSI context.

In our implementation of the unified domain formulation, the fluid computational domain $\mathit{\Omega}_f$ spans the whole computational domain $\mathit{\Omega}$. The solid computational domain $\mathit{\Omega}_s$ is fully contained within the fluid computational domain. The solid non-design domain $\mathit{\Omega}_n$, if any, along with the design domain constitute the solid domain. Such that:
\begin{gather}
	\mathit{\Omega}_s = \mathit{\Omega}_d \cup \mathit{\Omega}_n, \\
	\mathit{\Omega}_f := \mathit{\Omega}, \\
	\mathit{\Omega}_d \subseteq \mathit{\Omega}_s \subset \mathit{\Omega}_f.
\end{gather}

In $\mathit{\Omega}_d$, a design variable $\rho$ per each finite element determines whether it's pure fluid, pure solid, or in between. More specifically, $\rho$ represents the solid density of the design element; i.e. $\rho = 1$ for pure solid, $\rho = 0$ for pure fluid, and $0 < \rho < 1$ for porous solid. The solid/fluid interpolation required in density-based TO is achieved by:
\begin{enumerate}[(i)]
	\item appending the fluid momentum equation with a Brinkman penalization term $-\alpha (\rho) \mathbf{v}$ as a volume force where the inverse permeability\footnote{Yoon \cite[p.~502]{Yoon2010} called this term ``inverse permeability'' while Lundgaard et al. \cite[p.~971]{Lundgaard2018} called this term ``permeability''. By analogy to the original Brinkman equation \cite[p.~28]{Brinkman1947}, it should be called ``inverse permeability''.} $\alpha(\rho)$ is dependent on the design variables $\rho$, and
	\item interpolating the structural stiffness $E(\rho)$, i.e. Young's Modulus, as a function of the design variables $\rho$.	
\end{enumerate}

The Brinkman penalization term enforces zero fluid velocity for $\rho = 1$, porous flow conditions for $0 < \rho < 1$, and pure fluid flow for $\rho =0$. This way solid and fluid are allowed to co-exist within the same finite element, which is necessary to transform the discrete optimization problem into a continuous one. If a proper physical interpretation of the Brinkman penalization concept is sought, we assume that the porous media consists of a volumetric mixture of fluid and solid phases, the solid phase is homogeneous and isotropic on the macroscopic scale, and all the pores are interconnected \cite[p.~4]{Nield2017}.

For the fluid flow, the Navier-Stokes equations are used in their viscous, incompressible, steady-state form. The strong form of the partial differential equations is as follows \cite[p.~10]{Reddy2010}:
\begin{empheq}[right = \quad \empheqrbrace \, \text{in } \mathit{\Omega}_f]{gather}
	\boldsymbol{\nabla} \cdot \mathbf{v} = 0, \label{eq:contin} \\
	\rho_f \, (\mathbf{v} \cdot \boldsymbol{\nabla}) \mathbf{v} = \boldsymbol{\nabla} \cdot \boldsymbol{\upsigma}^f + \mathbf{f}^f - \alpha(\rho) \mathbf{v} , \label{eq:ns} \\ 
	\boldsymbol{\upsigma}^f = - p \mathbf{I} + \mu \left[ \boldsymbol{\nabla} \mathbf{v} + (\boldsymbol{\nabla} \mathbf{v})^T \right], \\
	\begin{split}
		\alpha(\rho) = \alpha_\text{max} + \left( 1 - \rho \right) \left( \alpha_\text{min} - \alpha_\text{max} \right) \frac{1 + p_\alpha}{1 - \rho + p_\alpha}. \label{eq:invrs_perm}
	\end{split}
\end{empheq}

\noindent where $\mathbf{v}$ is the fluid velocity, $\rho_f$ is the fluid density, $\boldsymbol{\upsigma}^f$ is the Cauchy fluid stress tensor, $\mathbf{f}^f$ is the external fluid force, $p$ is the hydrostatic pressure, and $\mu$ is the fluid dynamic viscosity. In Eq. \ref{eq:invrs_perm}, the inverse permeability $ \alpha (\rho) $ is dependent on the Brinkman penalization upper and lower limits, $\alpha_\text{max}$ and $\alpha_\text{min}$ respectively, and the Brinkman penalization interpolation parameter $p_\alpha$ that is dependent on the problem physics; e.g. Reynolds number\cite[p.~974]{Lundgaard2018}.

As for the structure, the Navier-Cauchy equations are used assuming linear elasticity with infinitesimal deformations under steady-state conditions. The strong form of the partial differential equations in the stress-divergence form is as follows \cite[p.~168]{Lai2010}:
\begin{empheq}[right = \quad \empheqrbrace \, \text{in } \mathit{\Omega}_s]{gather}
	\boldsymbol{\nabla} \cdot \boldsymbol{\upsigma}^s + \mathbf{f}^s = 0, \\
	\boldsymbol{\upsigma}^s = \mathbf{C}^s \, \boldsymbol{\upepsilon}, \\
	\boldsymbol{\upepsilon} = \frac{1}{2} \left[ \boldsymbol{\nabla} \mathbf{u} + (\boldsymbol{\nabla} \mathbf{u})^T \right], \\
	E(\rho) = E_\text{min} + (E_\text{max} - E_\text{min}) \rho^{\: \! p_E}. \label{eq:elastic_mod}
\end{empheq}

\noindent where $\boldsymbol{\upsigma}^s$ is the Cauchy solid stress tensor, $\mathbf{f}^s$ is the external solid force, $\mathbf{C}^s$ is the solid elasticity tensor, $\boldsymbol{\upepsilon}$ is the infinitesimal strain tensor, and $\mathbf{u}$ is the solid displacement. The elastic modulus $E(\rho)$ is interpolated using the modified SIMP approach where $E_\text{min}$ and $E_\text{max}$ are the lower and upper limits on $E$ representing the void and bulk elastic moduli, respectively. $p_E$ is the elastic modulus penalization parameter. The essential boundary conditions are defined as follows:
\begin{alignat}{3}
	\text{Fluid No-slip:} & \qquad \mathbf{v} = 0 & \quad \text{on } & \mathit{\Gamma}_{\mathbf{v}_0}, \label{eq:bc_noslip} \\
	\text{Fluid Inlet:} & \qquad \mathbf{v} = \mathbf{v}_\text{in} & \quad \text{on } & \mathit{\Gamma}_{\mathbf{v}_\text{in}}, \label{eq:bc_v_in} \\
	\text{Fluid Outlet:} & \qquad p = 0 & \quad \text{on } & \mathit{\Gamma}_{\mathbf{v}_\text{out}}, \label{eq:bc_p_out} \\
	\text{Structural Prescribed Displacement:} & \qquad \mathbf{u} = 0 & \quad \text{on } & \mathit{\Gamma}_{\mathbf{u}_0}.
\end{alignat}

Note that the fluid no-slip boundary condition in Eq. \ref{eq:bc_noslip} is only defined on the external domain boundaries aside from the inlet and outlet such that  $ \mathit{\partial \Omega}_f = \mathit{\Gamma}_{\mathbf{v}_0} \cup \mathit{\Gamma}_{\mathbf{v}_\text{in}} \cup \mathit{\Gamma}_{\mathbf{v}_\text{out}} $. The volume force term appended to the fluid momentum in Eq. \ref{eq:ns} automatically enforces a no-slip condition wherever needed within the solid domain and/or its fluid-structure interface. The natural boundary condition in an FSI problem is defined as follows:
\begin{equation}
	\boldsymbol{\upsigma}^f \cdot \mathbf{n}^f = \boldsymbol{\upsigma}^s \cdot \mathbf{n}^s \qquad \text{on } \mathit{\Gamma}_\text{FSI}
	\label{eq:trac_eq}
\end{equation}

\noindent where $\mathbf{n}^f$ and $\mathbf{n}^s$ are the normals to the fluid and solid surfaces, respectively. After calculating the fluidic forces per each finite element and before the global assembly, the forces must be multiplied by the pressure-coupling filtering function defined as\citep[p.~607]{Yoon2010}:
\begin{equation}
	\mathit{\Upsilon}(\rho_i) = \mathit{\Upsilon}_\text{min} + (\mathit{\Upsilon}_\text{max} - \mathit{\Upsilon}_\text{min}) \ \rho_i^{\, p_\mathit{\Upsilon}}
\end{equation}

\noindent where $\mathit{\Upsilon}_\text{max}$ and $\mathit{\Upsilon}_\text{min}$ are assumed 1 and 0 respectively. $p_\mathit{\Upsilon}$ is the pressure-coupling filter function interpolation parameter. The main intention of using the pressure-coupling filtering function is to prevent severe deformation of low density elements and potential singularities in solving the structural finite element system.

While the fluid-structure interface $\mathit{\Gamma}_\text{FSI}$ is explicitly defined in pure fluid-structure interactions (i.e. no porous media), this is not the case in density-based TO problems where all solids are porous to some degree. Hence, to apply the traction coupling condition in this formulation properly, $\mathit{\Gamma}_\text{FSI}$ must be taken as all edges of all solid finite elements. More discussion on this condition and how to apply it in a finite element form follows in subsection \ref{ssec:traction}.

In the next section, we transform the strong form of the governing equations into the discretized finite element form and extend the traction equilibrium condition described generally in Eq. \ref{eq:trac_eq} to include viscous stresses.

\section{Finite Element Formulations}
\label{sec:fin_elem_form}

\subsection{Fluid and Solid Governing Equations}
\label{ssec:govern_eqns}

In this work, conformal or matching discretization is assumed at the fluid-structure interface to ensure simple, yet accurate, force coupling. For non-conformal meshing, other mapping techniques can be used to enforce momentum conservation at the interface \cite{Farhat1998a}. The discretized finite element form of Navier-Stokes equations are obtained by multiplying the strong form of the continuity and momentum equations (Eqs. \ref{eq:contin} and \ref{eq:ns}) with suitable weight functions and integrating the resulting expressions over the appropriate computational domain. The divergence theorem (i.e. integration by parts) is typically implemented to reduce the continuity requirements on the interpolation/shape functions used to approximate the solution by moving the differentiation to the weight functions \cite{Reddy2010}. 

In this work, we implement the \textit{standard Galerkin method of weighted residuals} where the interpolation/shape functions themselves are used as test/weight functions. The resulting finite element model is of the \textit{velocity-pressure (or mixed)} type where both velocity and pressure are unknowns and are solved for simultaneously. We utilize a \textit{quadrilateral meshing} which is regular and structured wherever possible. To satisfy the \textit{Ladyzhenskaya-Babuska-Brezzi} condition (cf. \ \cite[p.~176]{Reddy2010}), \textit{\textcolor{black}{Q2Q1} Lagrangian} finite elements (i.e. 9 velocity nodes and 4 pressure nodes) are used with a low to moderate Reynolds number to avoid using stabilization techniques that artificially - but not necessarily accurately - dampen the discontinuities. The continuity equation (Eq. \ref{eq:contin}) is typically weighted by the pressure shape function $\boldsymbol{\Phi}$ while the momentum equation (Eq. \ref{eq:ns}) is typically weighted by the velocity shape function $\boldsymbol{\Psi}$, both are in column vector form. The resulting finite element system in 2D \textit{on the elemental level} is as follows:
\begin{equation}
	\underbrace{\begin{bmatrix}
			2 \mathbf{K}_{11} + \mathbf{K}_{22} + \mathbf{C(v)} & \mathbf{K}_{12} & - \mathbf{Q}_1 \\[4pt]
			\mathbf{K}_{21} & \mathbf{K}_{11} + 2 \mathbf{K}_{22} + \mathbf{C(v)} & - \mathbf{Q}_2 \\[4pt]
			- \mathbf{Q}_1^T & - \mathbf{Q}_2^T & \mathbf{0}
	\end{bmatrix}}_{\text{Conservation of Momentum and Mass}}
	\begin{Bmatrix}
		\mathbf{\hat{v}}_1 \\[4pt]
		\mathbf{\hat{v}}_2 \\[4pt]
		\mathbf{\hat{p}}
	\end{Bmatrix} +
	%%%
	\underbrace{\begin{bmatrix}
			\mathbf{A} & \mathbf{0} & \mathbf{0} \\[4pt]
			\mathbf{0} & \mathbf{A} & \mathbf{0} \\[4pt]
			\mathbf{0} & \mathbf{0} & \mathbf{0}
	\end{bmatrix}}_{\substack{\text{Brinkman} \\ \text{Penalization}}}
	\begin{Bmatrix}
		\mathbf{\hat{v}}_1 \\[4pt]
		\mathbf{\hat{v}}_2 \\[4pt]
		\mathbf{\hat{p}}
	\end{Bmatrix} = 
	\begin{Bmatrix}
		\mathbf{F}^f_1 \\[4pt]
		\mathbf{F}^f_2 \\[4pt]
		\mathbf{0}
	\end{Bmatrix}
	\label{eq:main_fea}
\end{equation}

The coefficient matrices in the finite element form are defined as ($\int_{-1}^{+1} \int_{-1}^{+1}$ and $\dd{\xi} \dd{\eta}$ are implied)\footnote{Summation is implied on repeated indices in $\mathbf{C(v)}$ but not in $\mathbf{K}_{ij}$.}:
\begingroup
\allowdisplaybreaks
\begin{align}
	\mathbf{K}_{ij} & = \mu \, \pdv{\mathbf{\Psi}}{x_i} \pdv{\mathbf{\Psi}}{x_j} ^T |\mathbf{J}|, \label{eq:2} \\
	%%%
	\mathbf{C(v)} & = \rho \, \mathbf{\Psi} \left[ \left( \mathbf{\Psi}^T \mathbf{\hat{v}}_i \right) \pdv{\mathbf{\Psi}}{x_i} ^T \right] |\mathbf{J}|, \\
	%%%
	\mathbf{Q}_i & = \pdv{\mathbf{\Psi}}{x_i} \, \mathbf{\Phi}^T |\mathbf{J}|, \\
	%%%
	\mathbf{A} & = \alpha(\rho) \mathbf{\Psi} \, \mathbf{\Psi}^T |\mathbf{J}|.
\end{align} 
\endgroup

\noindent where $\mathbf{\hat{v}}_1, \mathbf{\hat{v}}_2, \text{and } \mathbf{\hat{p}}$ are the nodal velocities in $x$ and $y$ and nodal pressures, respectively. $|\mathbf{J}|$ is the Jacobian determinant and $\xi$ and $\eta$ are the natural coordinates. $\mathbf{F}^f_1$ and $\mathbf{F}^f_2$ are externally applied nodal fluid forces which are not used in this work since the fluidic boundary conditions (Eqs. \ref{eq:bc_noslip} to \ref{eq:bc_p_out}) are implemented directly by setting nodal velocities/pressures to their appropriate values. Appropriate global assembly of Eq. \ref{eq:main_fea} is implemented and the resulting nonlinear system is solved using the undamped Newton-Raphson method \cite[p.~190]{Reddy2010}.

As for the Navier-Cauchy equations, their discretized finite element form is obtained using the \textit{principle of virtual displacement} following the procedure in \ \cite[p.~153]{Bathe2014finite}. To simplify the force coupling, we utilize \textit{Quadratic Lagrange} finite elements (i.e. 9 displacement nodes) for the solid domain. In 2D, the resulting system \textit{on the elemental level} is as follows:
\begin{align}
	\big[ \mathbf{k}^s \big] \big\{ \mathbf{\hat{u}}\big\} - \big\{ \mathit{\Upsilon}(\rho) \ \mathbf{F}^s \big\} = \mathbf{0}
\end{align}

\noindent where $\mathbf{\hat{u}}$ is the nodal displacements and $\mathbf{F}^s$ is the nodal solid forces. The solid stiffness matrix $\mathbf{k}^s$ for \textit{plane strain} conditions is obtained following the procedure in \ \cite[p.~164]{Bathe2014finite} with the only difference being the dependence of the elastic modulus on the design variable as detailed in Eq. \ref{eq:elastic_mod}, and appropriate global assembly is also performed to construct the global solid stiffness matrix. Isoparametric elements are used for both fluid and solid domains, so the geometry is described by the same shape functions describing the velocities/displacements.

\subsection{Traction Coupling}
\label{ssec:traction}

In literature, the traction equilibrium condition was derived considering hydrostatic stresses only (cf. \ \cite[p.~601]{Yoon2010} and \ \cite[p.~989]{Lundgaard2018}), which is suitable for some cases with negligible viscous stresses. In this work, we extend the traction equilibrium condition to include the viscous stresses as well. The fluidic forces applied on the solid structure are calculated directly from Eq. \ref{eq:trac_eq} as follows:
\begin{equation}
	\mathbf{f}^s = \int_{\mathit{\Gamma}_\text{FSI}} \boldsymbol{\upsigma}^f \cdot \mathbf{n}^f \dd{\mathit{\Gamma}}
\end{equation}

\noindent which can be put into the finite element form through multiplying by the velocity shape - and test - functions vector $\boldsymbol{\Psi}$ then discretizing:
\begin{equation}
	\mathbf{F}^s_i = \int_{\mathit{\Gamma}_\text{FSI}} \boldsymbol{\Psi} \left( \sigma^f_{ij} n^f_j \right) \dd{\mathit{\Gamma}}
	\label{eq:surf_force_general}
\end{equation}

After transformation from global to natural coordinates, Eq. \ref{eq:surf_force_general} can be further detailed as follows:
\begin{align}
	\mathbf{F}^s_x & = \int_{-1}^{+1} \boldsymbol{\Psi} \left[ \left( - \mathbf{\Phi}^T \mathbf{\hat{p}} + 2 \mu \pdv{\mathbf{\Psi}}{x} ^T \mathbf{\hat{v}}_1 \right) n^f_x + \mu \left( \pdv{\boldsymbol{\Psi}}{y} ^T \mathbf{\hat{v}}_1 + \pdv{\boldsymbol{\Psi}}{x} ^T \mathbf{\hat{v}}_2 \right) n^f_y \right] |\mathbf{J}_s| \dd{\xi_s},	\label{eq:surf_force_x_org}
	\\
	\mathbf{F}^s_y & = \int_{-1}^{+1} \boldsymbol{\Psi} \left[ \left( - \mathbf{\Phi}^T \mathbf{\hat{p}} + 2 \mu \pdv{\mathbf{\Psi}}{y} ^T \mathbf{\hat{v}}_2 \right) n^f_y + \mu \left( \pdv{\boldsymbol{\Psi}}{y} ^T \mathbf{\hat{v}}_1 + \pdv{\boldsymbol{\Psi}}{x} ^T \mathbf{\hat{v}}_2 \right) n^f_x \right] |\mathbf{J}_s| \dd{\xi_s}.
	\label{eq:surf_force_y_org}
\end{align}

\noindent where the normal vector components $n^f_x$ and $n^f_y$ can be calculated from elemental nodal coordinates, and $\dd{\xi_s}$ designates the natural coordinate along the finite element edge of interest in 2D. Note that $\mathbf{J}_s$ is the surface Jacobian that is different from the typical area Jacobian $\mathbf{J}$ (cf. \ \cite[p.~97]{Reddy2010} or \ \cite[p.~356]{Bathe2014finite}). 

In previous works on TOFSI problems (and in works on FSI with porous media in general), the divergence theorem is typically used to transform the surface integral into a volume integral as follows (cf. \ \cite[p.~598]{Yoon2010} and \ \cite[p.~989]{Lundgaard2018}):
\begin{equation}
	\oiint\limits_S \left( \mathbf{\boldsymbol{\Psi}} \sigma^f_{ij} \right) n^f_j \dd{\mathit{\Gamma}} = \iiint\limits_V \boldsymbol{\nabla} \cdot (\boldsymbol{\Psi} \sigma^f_{ij}) \dd{V} = \iiint\limits_V \left( \pdv{\boldsymbol{\Psi}}{x_j} \sigma^f_{ij} + \boldsymbol{\Psi} \pdv{\sigma^f_{ij}}{x_j} \right) \dd{V}
	\label{eq:div_theorm}
\end{equation}

By applying the divergence theorem in Eq. \ref{eq:div_theorm} to Eq. \ref{eq:surf_force_general}, the surface integral form in Eqs. \ref{eq:surf_force_x_org} and \ref{eq:surf_force_y_org} is transformed into a volume integral form as follows:
\begingroup
\allowdisplaybreaks
\begin{align}
	\begin{split}
		\mathbf{F}^s_x = \int_{-1}^{+1} \int_{-1}^{+1} \Bigg[ & \pdv{\boldsymbol{\Psi}}{x} \left( - \mathbf{\Phi}^T \mathbf{\hat{p}} + 2 \mu \pdv{\mathbf{\Psi}}{x} ^T \mathbf{\hat{v}}_1 \right) +  \pdv{\boldsymbol{\Psi}}{y} \left( \mu  \pdv{\boldsymbol{\Psi}}{y} ^T \mathbf{\hat{v}}_1 + \mu \pdv{\boldsymbol{\Psi}}{x} ^T \mathbf{\hat{v}}_2 \right) \\ 
		& + \boldsymbol{\Psi} \left( - \pdv{\mathbf{\Phi}}{x} ^T \mathbf{\hat{p}} \, + 
		2 \mu \pdv{\boldsymbol{\Psi}}{x,x} \mathbf{\hat{v}}_1 + \mu \pdv{\boldsymbol{\Psi}}{y,y} ^T \mathbf{\hat{v}}_1 + \mu \pdv{\boldsymbol{\Psi}}{x,y} ^T \mathbf{\hat{v}}_2 \right) \Bigg] |\mathbf{J}| \dd{\xi} \dd{\eta},
		\label{eq:vol_force_x_org}
	\end{split} \\
	\begin{split}
		\mathbf{F}^s_y = \int_{-1}^{+1} \int_{-1}^{+1} \Bigg[ & \pdv{\boldsymbol{\Psi}}{x} \left( \mu \pdv{\boldsymbol{\Psi}}{y} ^T \mathbf{\hat{v}}_1 + \mu  \pdv{\boldsymbol{\Psi}}{x} ^T \mathbf{\hat{v}}_2 \right) 
		+ \pdv{\boldsymbol{\Psi}}{y} \left( - \mathbf{\Phi}^T \mathbf{\hat{p}} + 2 \mu \pdv{\mathbf{\Psi}}{y} ^T \mathbf{\hat{v}}_2 \right) \\
		& + \boldsymbol{\Psi} \left( \mu \pdv{\boldsymbol{\Psi}}{y,x} ^T \mathbf{\hat{v}}_1 + 
		\mu \pdv{\boldsymbol{\Psi}}{x,x} ^T \mathbf{\hat{v}}_2 - \pdv{\mathbf{\Phi}}{y} ^T \mathbf{\hat{p}} + 2 \mu \pdv{\boldsymbol{\Psi}}{y,y} \mathbf{\hat{v}}_2 \right) \Bigg] |\mathbf{J}| \dd{\xi} \dd{\eta}.
		\label{eq:vol_force_y_org}
	\end{split}
\end{align}
\endgroup

Note the appearance of velocity first derivatives in the force coupling equations using surface integrals (i.e. Eqs. \ref{eq:surf_force_x_org} and \ref{eq:surf_force_y_org}). In the volume integral form of the force coupling (i.e. Eqs. \ref{eq:vol_force_x_org} and \ref{eq:vol_force_y_org}), velocity second derivatives and pressure first derivatives appear due to the use of the divergence theorem. The appearance of these derivatives and their effect on traction equilibrium when using total stress definitions is the topic of discussion in section \ref{sec:total_traction_section}.

\subsection{Sensitivity Analysis}
\label{ssec:snst_anls_general}

Typically, in TO problems, the number of design variables greatly exceeds the number of constraints. Hence, it is better from a computational perspective to use the \textit{adjoint} method in calculating the sensitivities as opposed to the \textit{direct} method (cf. \cite{Maute2003} for a discussion on the adjoint method in an FSI context). Given an objective function (or a constraint) $f$, the sensitivity of $f$ w.r.t. the design variable $\rho_i$ is calculated as follows:
\begingroup
\allowdisplaybreaks
\begin{align}
	& \dv{f}{\rho_i} = \pdv{f}{\rho_i} +	\pdv{f}{\mathbf{r}} ^T \pdv{\mathbf{r}}{\rho_i}, \\
	%%%
	& \dv{\mathbb{R}}{\rho_i} = \mathbf{0}, \quad \rightarrow \quad \pdv{\mathbb{R}}{\rho_i} + \pdv{\mathbb{R}}{\mathbf{r}} \pdv{\mathbf{r}}{\rho_i} = \mathbf{0}, \\
	%%%
	& \pdv{\mathbf{r}}{\rho_i} = - \pdv{\mathbb{R}}{\mathbf{r}} ^{-1}
	\pdv{\mathbb{R}}{\rho_i}, \\
	%%%
	& \dv{f}{\rho_i} = \pdv{f}{\rho_i} -	\pdv{f}{\mathbf{r}} ^T \pdv{\mathbb{R}}{\mathbf{r}} ^{-1} \pdv{\mathbb{R}}{\rho_i}, \\
	%%%
	& \dv{f}{\rho_i} = \pdv{f}{\rho_i} -
	\left( \pdv{\mathbb{R}}{\mathbf{r}} ^{-T}	\pdv{f}{\mathbf{r}} \right)^T		\pdv{\mathbb{R}}{\rho_i}, \\
	%%%
	& \boldsymbol{\uplambda} = \pdv{\mathbb{R}}{\mathbf{r}} ^{-T}
	\pdv{f}{\mathbf{r}}, \quad \rightarrow \quad \pdv{\mathbb{R}}{\mathbf{r}} ^T \boldsymbol{\uplambda} = \pdv{f}{\mathbf{r}}, \\
	%%%
	& \dv{f}{\rho_i} = \pdv{f}{\rho_i} - \boldsymbol{\uplambda}^T
	\pdv{\mathbb{R}}{\rho_i}.
\end{align}
\endgroup

\noindent where $\mathbf{r}$ is a vector of all state variables (i.e. solid displacements and fluid velocities and pressures), $\mathbb{R}$ represents the set of all (i.e. solid and fluid) finite element equations,  and $\boldsymbol{\uplambda}$ is the adjoint vector which is distinctive for each objective function or constraint.

It can be seen that the derivatives of the force coupling equations w.r.t. the fluid state variables are needed as a part of the $\partial \mathbb{R} / \partial \mathbf{r}$ term. For the case of hydrostatic coupling implemented earlier in \ \cite{Yoon2010} and \  \cite{Lundgaard2018}, only the derivatives w.r.t. the nodal pressures are needed. These derivatives are evaluated for the surface integral form directly from Eqs. \ref{eq:surf_force_x_org} and \ref{eq:surf_force_y_org} as follows:
\begin{align}
	\pdv{\mathbf{F}^s_x}{\mathbf{\hat{p}}} & = - \int_{-1}^{+1} \mathbf{\Psi} \ \mathbf{\Phi}^T n_x \ |\mathbf{J}_s| \dd{\xi_s}, \\
	\pdv{\mathbf{F}^s_y}{\mathbf{\hat{p}}} & = - \int_{-1}^{+1} \mathbf{\Psi} \ \mathbf{\Phi}^T n_y \ |\mathbf{J}_s| \dd{\xi_s}.
\end{align}

\noindent and for the volume integral form, they are evaluated directly from Eqs. \ref{eq:vol_force_x_org} and \ref{eq:vol_force_y_org} as follows:
\begin{align}
	\pdv{\mathbf{F}^s_x}{\mathbf{\hat{p}}} & = - \int_{-1}^{+1} \int_{-1}^{+1} \bigg( \pdv{\mathbf{\Psi}}{x} \ \mathbf{\Phi}^T + \mathbf{\Psi} \pdv{\mathbf{\Phi}}{x} ^T \bigg) |\mathbf{J}| \dd{\xi} \dd{\eta}, \\
	\pdv{\mathbf{F}^s_y}{\mathbf{\hat{p}}} & = - \int_{-1}^{+1} \int_{-1}^{+1} \bigg( \pdv{\mathbf{\Psi}}{y} \ \mathbf{\Phi}^T + \mathbf{\Psi} \pdv{\mathbf{\Phi}}{y} ^T \bigg) |\mathbf{J}| \dd{\xi} \dd{\eta}.
\end{align}

\section{Description of Test Problems}
\label{sec:tst_prblm}

Before we discuss the details of the total stress equilibrium equations in section \ref{sec:total_traction_section}, we first describe the setup of the numerical example to be used in the following sections.

\color{black}

The \textit{original} column in a channel problem was first discussed in a TOFSI context in \ \cite[p.~610]{Yoon2010} and has been used later as a benchmark problem in a number of works on TOFSI. In \ \cite{Lundgaard2018}, Lundgaard et al. developed a \textit{modified} version\footnote{\color{black} Lundgaard et al. \ \cite{Lundgaard2018} called this problem ``the wall".}, where they \textbf{(i)} increased the relative size of the design domain with respect to the whole computational domain, \textbf{(ii)} enlarged the geometry from the micro to the macro scale, and \textbf{(iii)} generally strengthened the fluid-structure dependency. We utilize the modified version of this problem as the numerical example of choice in this work. In the remainder of this work, we drop the designation ``modified" and simply refer to this problem as the ``column in a channel" test problem.

\begin{figure}[b]
	\centering
	\captionsetup{font={color=black}}
	\includegraphics[height=0.3\textwidth]{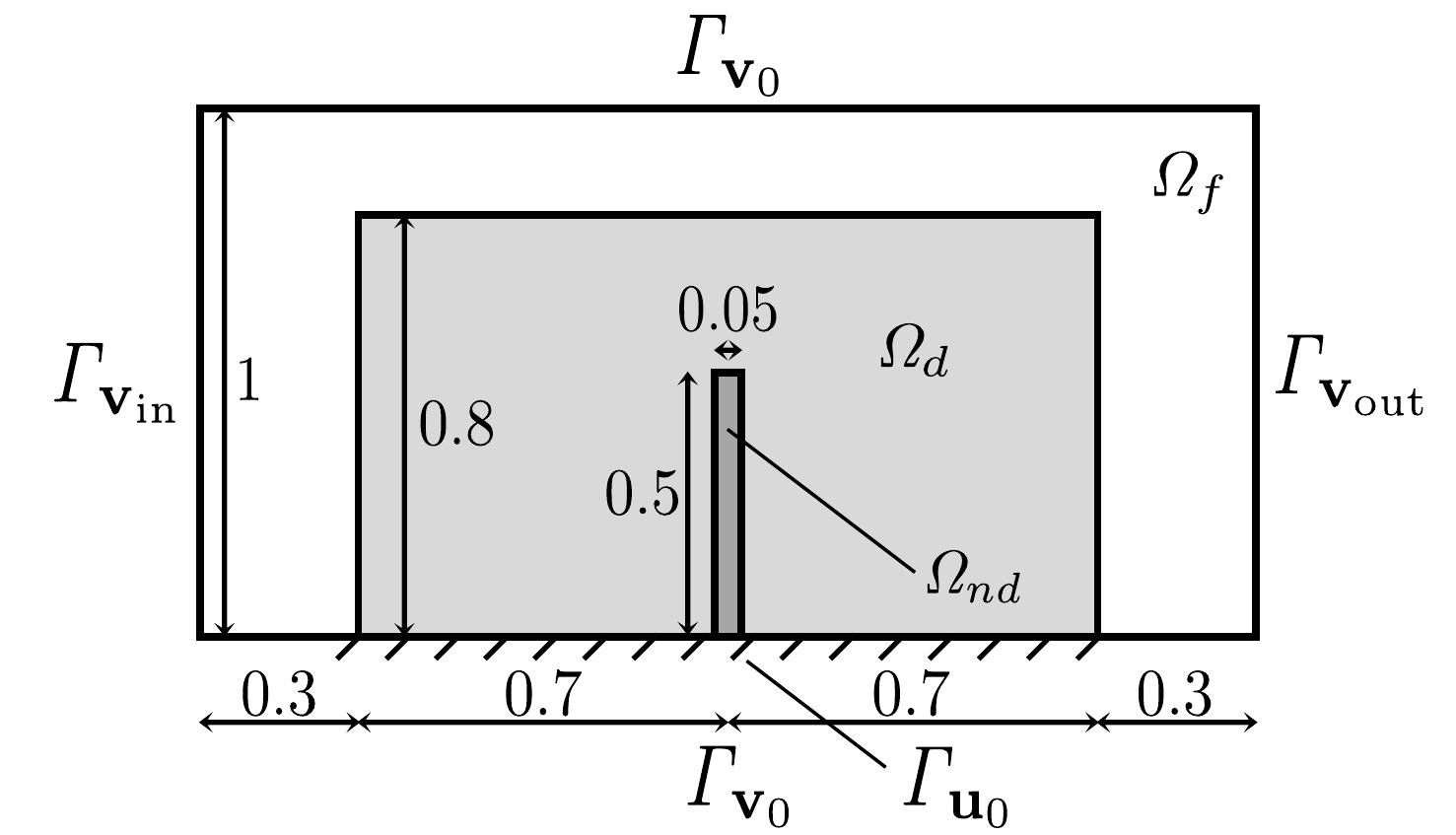}
	\caption{Diagram of the modified column in a channel problem by Lundgaard et al. \ \citep[p.~975]{Lundgaard2018}.}
	\label{fig:column_in_a_channel_sigmund_diagram}
\end{figure}

As detailed in Fig. \ref{fig:column_in_a_channel_sigmund_diagram}, the problem features a 1.4 $\times$ 0.8 m rectangular design space (light gray) placed inside a 2 $\times$ 1 m rectangular channel. To avoid trivial solutions, a 0.5 $\times$ 0.05 m non-design, elastic, solid column (dark gray) is placed within the design space to force the optimizer into a more sophisticated solution than a simple bump at the bottom of the channel. Note that even though the non-design column is not included in the design domain, it is still governed by the structural equations and is included in the calculation of the objective function.

The top and bottom surfaces of the channel $\mathit{\Gamma}_{\mathbf{v}_0}$ have a no-slip condition applied. A parabolic laminar flow profile is applied at the inlet $\mathit{\Gamma}_{\mathbf{v}_\text{in}}$ on the left with a maximum velocity of 1 m/s, and a zero pressure condition is applied at the outlet $\mathit{\Gamma}_{\mathbf{v}_\text{out}}$ on the right. The bottom surface of the design and non-design domains $\mathit{\Gamma}_{\mathbf{u}_0}$ is fixed to a ground structure. The fluid density $\rho_f$ is assumed 1  kg/m$^3$. The maximum and minimum Brinkman penalization limits are taken as $\alpha_\text{max} = 1e+9$ and $\alpha_\text{min} = 0$ as suggested in \ \citep[p.~974]{Lundgaard2018} and further discussed in \ \cite{abdelhamid2023calculation}. The solid elastic moduli parameters $E_\text{max}$ and $E_\text{min}$ are assumed $1e+5$ Pa and $1e-5$ Pa, respectively, while Poisson's ratio $\nu$ is taken as 0.3. The fluid dynamic viscosity $\mu$ is to be varied to produce different Reynolds numbers which are calculated according to Eq. \ref{eq:reynolds_number} as follows:
\begin{equation}
	Re = \frac{v_\text{max} \cdot \rho_f \cdot L_c}{\mu}.
	\label{eq:reynolds_number}
\end{equation}

\noindent where $L_c$ is a characteristic length that is taken as the width of the inlet, i.e. 1 m, in agreement with \ \citep[p.~971]{Lundgaard2018}. A regular, structured mesh of 45,300 (i.e. 302 in $x$ and 150 in $y$) quadrilateral Q2Q1 finite elements is implemented.

\color{black}

\section{Total Stress Equilibrium in Porous Fluid-Structure Interaction Problems}
\label{sec:total_traction_section}

It's well-known that in the Galerkin finite element formulation with typical Lagrangian shape functions, the primitive variables (i.e. velocities and pressure in fluid flow) are continuous across element boundaries, however, no such continuity is enforced on their derivatives\cite{Gresho2000}. This discontinuity in derivatives poses a challenge to the calculation of relevant derived quantities that utilize these derivatives such as fluxes and - pertinent to this work - traction at the fluid-structure interface.

A number of techniques have been discussed in literature to overcome this issue. Bazilevs and Hughes \cite{Bazilevs2007} enforced the Dirichlet no-slip boundary condition weakly in the context of turbulent fluid flow instead of the typical strong nodal imposition. Their mean velocity profile results showed good agreement with comparable direct numerical simulation results. Oshima et al. \cite{OSHIMA1998} developed a consistent method as a post-processing step where the fluxes or stresses are treated as unknowns which are solved for weakly after the typical Galerkin finite element solution is obtained. Zienkiewicz and Zhu \cite{Zienkiewicz1992, Zienkiewicz1992a} developed a superconvergent recovery technique for calculating the derivatives of finite element solutions at nodes. In this work, we implement the discrete version of this technique to overcome these discontinuities in the calculation of total stress coupling enforced at the fluid-structure interface. In particular, we apply the technique to recover the velocity derivatives w.r.t. coordinates $x$ and $y$.

\subsection{Recovery of Discontinuous Velocity Derivatives}

The superconvergent recovery technique by Zienkiewicz and Zhu \cite{Zienkiewicz1992, Zienkiewicz1992a} has its basis in the discovery of superconvergent points within each finite element where derivatives are known to have \textit{extraordinary} accuracy compared to the rest of the finite element (cf. \ \cite{Zienkiewicz1992, Zienkiewicz1992a} and references therein). These points correspond to reduced integration, so for the quadratic quadrilateral finite elements used in this work \textcolor{black}{(3 $\times$ 3 velocity nodes)}, these superconvergent points are the $2 \times 2$ Gaussian integration points with equal weights. The procedure we implement is described in detail in \ \cite[p.~1335]{Zienkiewicz1992} and in a fluid mechanics context in \ \cite[p.~952]{Gresho2000}. Once the recovered velocity derivatives at the nodes are calculated, they can be interpolated within each finite element using the same shape functions used to interpolate the primitive variables, i.e. velocities. Hence, this technique lends itself well to the volume integral form of force coupling in Eqs. \ref{eq:vol_force_x_org} and \ref{eq:vol_force_y_org} as the second derivatives of velocities are evaluated at the $3 \times 3$ Gaussian points within each finite element using the appropriate derivatives of the shape functions multiplied by the recovered nodal velocity derivatives as detailed later in Eqs. \ref{eq:vol_force_x_recov} and \ref{eq:vol_force_y_recov}.

We reiterate the recovery procedure here as it's pertinent to the sensitivity analysis discussion in subsection \ref{ssec:snst_anls_recovery}. We designate the recovered values with a bar on top, so the recovered nodal velocity derivatives are $\overline{\mathbf{\hat{v}}_{1,x}}$, $\overline{\mathbf{\hat{v}}_{1,y}}$, $\overline{\mathbf{\hat{v}}_{2,x}}$, and $\overline{\mathbf{\hat{v}}_{2,y}}$. The scalar recovered velocity derivative within a finite element $\overline{v_{i,j}}$ is related to the corresponding recovered nodal values of that particular element by interpolation using velocity shape functions as in:
\begin{equation}
	\overline{v_{i,j}} = \mathbf{\Psi}^T \overline{\mathbf{\hat{v}}_{i,j}}
\end{equation}

Each nodal value of the vector $\overline{\mathbf{\hat{v}}_{i,j}}$ is calculated within one or more\footnote{An average value is used when a certain node is shared between more than one recovery patch.} recovery patches where it's represented by a polynomial expansion of the following form:
\begin{equation}
	\overline{\hat{v}_{i,j}} = {\mathbf{P}_\text{O} (x,y)}^T \mathbf{a}_{i,j} \label{eq:p_o}
\end{equation}

\noindent where $\mathbf{P}$ contains the appropriate polynomial terms and is a function in coordinates $x$ and $y$ within the patch\footnote{In the original reference \ \cite{Zienkiewicz1992a}, $\mathbf{P}$ is a row vector, however in our implementation we define it as a column vector for consistency.}, and $\mathbf{a}_{i,j}$ - a vector of unknown parameters of the same size as $\mathbf{P}$ - is solved for within each patch through a least square fit minimization problem as follows:
\begin{equation}
	\sum_{k = 1}^{n} \mathbf{P}_\times (x_k,y_k) \ {\mathbf{P}_\times (x_k,y_k)}^T \ \mathbf{a}_{i,j} = \sum_{k = 1}^{n} \mathbf{P}_\times (x_k,y_k) \ v_{i,j} \label{eq:p_x}
\end{equation}

\noindent and $v_{i,j}$ is evaluated at the superconvergent points (i.e. 2 $\times$ 2 Gaussian points for a \textcolor{black}{Q2Q1} element) within each finite element in the recovery patch using the appropriate relation as in:
\begin{equation}
	v_{i,j} = \pdv{\mathbf{\Psi}}{x_j} ^T \mathbf{\hat{v}}_i
\end{equation}

The vector $\mathbf{P}$ is designated with a subscript $\square_\text{O}$ when evaluated at the to-be-recovered nodes as in Eq. \ref{eq:p_o} and with a subscript $\square_\times$ when evaluated at the superconvergent Gaussian points as in Eq. \ref{eq:p_x}. This distinction is critical in the sensitivity analysis discussion in subsection \ref{ssec:snst_anls_recovery}.

Note that within each patch, the four different derivatives - i.e. $\overline{\hat{v}_{1,x}}$,  $\overline{\hat{v}_{1,y}}$, $\overline{\hat{v}_{2,x}}$, and $\overline{\hat{v}_{2,y}}$ - are recovered using the same $\mathbf{P}$, however, a different $\mathbf{a}_{i,j}$ is calculated for each derivative, hence the subscript designation $\square_{i,j}$ for $\mathbf{a}_{i,j}$. A few remarks about our implementation of the superconvergent patch recovery technique are in order:
\begin{enumerate}[a.]
	\item Assuming $\mathbf{P} = \{1, \ x, \ y, \ xy, \ x^2, \ y^2, \ xy^2, \ x^2y, \ x^2y^2\}^T$ for \textcolor{black}{Q2Q1} elements and depending on the geometric scale of the problem, singularities might arise in solving the linear system $\mathbf{a} = \mathbf{A}^{-1} \mathbf{b}$ (where $\mathbf{A} = \mathbf{P} \ \mathbf{P}^T$ and $\mathbf{b} = \mathbf{P} \ v_{i,j}$) due to how the matrix $\mathbf{A}$ is constructed. For instance at point $(x,y) = (1,1)$, all terms of $\mathbf{P}$ are equal to 1. To overcome this issue, we elected to scale all $x$ and $y$ coordinates within each patch. We arrived at the following scaling relations using - may we say - `educated' trial and error:
	\begin{gather}
		x' = x \cdot s_x + t_x, \quad y' = y \cdot s_y + t_y, \\
		s_x = \frac{1}{ \max(\mathbf{\hat{x}}) - \min(\mathbf{\hat{x}}) }, \quad	t_x = - \min(\mathbf{\hat{x}}) \cdot s_x + 1.0, \\
		s_y = \frac{1}{ \max(\mathbf{\hat{y}}) - \min(\mathbf{\hat{y}}) }, \quad  t_y = - \min(\mathbf{\hat{y}}) \cdot s_y - 0.5.
	\end{gather}
	\noindent where $\mathbf{\hat{x}}$ and $\mathbf{\hat{y}}$ are the nodal coordinates of a random element within each recovery patch, usually the first element to be considered in the programming {\tt for} loop within each patch. It's possible that these scaling relations might need to be tuned based on the geometry of each problem as this was not rigorously tested for different problems.
	
	\item Depending on the regularity of the mesh and the location of the to-be-recovered nodes, some patches may contain less or more than 4 elements. In particular, for patches with less than 4 elements, the degree of the polynomial - hence the number of terms - used to construct $\mathbf{P}$ has to be reduced to avoid singularities. From our numerical experiments, patches containing two elements (e.g. the center patch node is at an external boundary) may use 4 polynomial terms, patches containing three elements (e.g. in the middle of an irregular mesh) may use 6 polynomial terms, and patches containing 4 or more elements may use all the 9 polynomial terms.
\end{enumerate}

\begin{figure}[b]
	\captionsetup{font={color=black}}
	\captionsetup[subfigure]{font+={color=black}}
	\subfloat[Original $\dfrac{\partial v_1}{\partial x}$]{\includegraphics[width=0.25\textwidth]{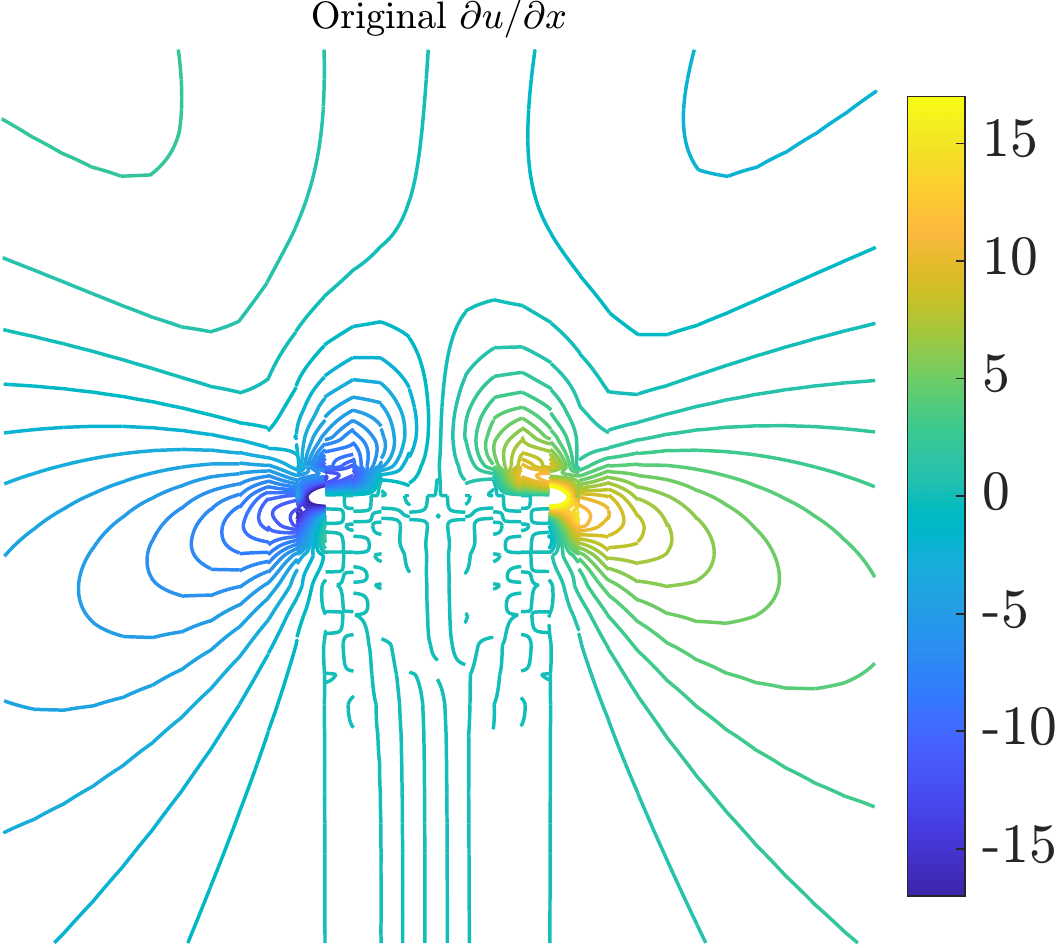}}
	\subfloat[Recovered $\dfrac{\partial v_1}{\partial x}$]{\includegraphics[width=0.25\textwidth]{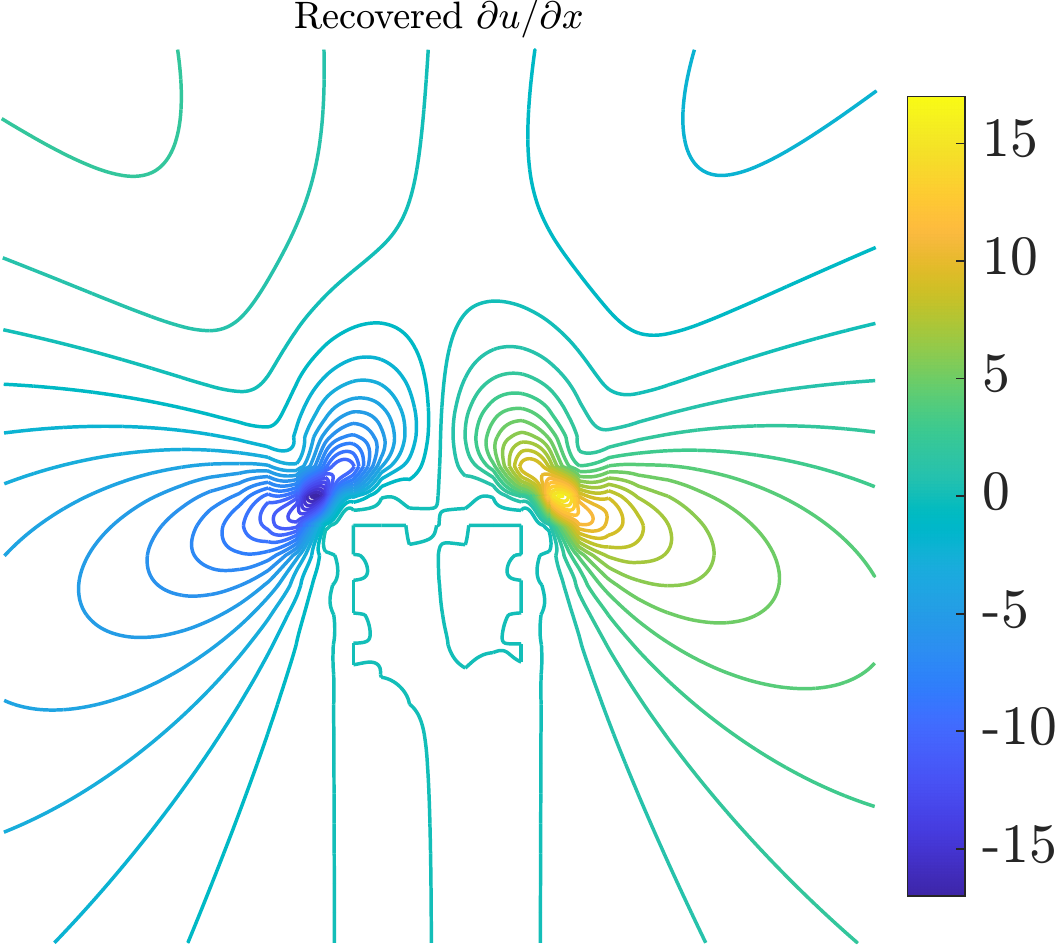}}
	\subfloat[Original $\dfrac{\partial v_2}{\partial y}$]{\includegraphics[width=0.25\textwidth]{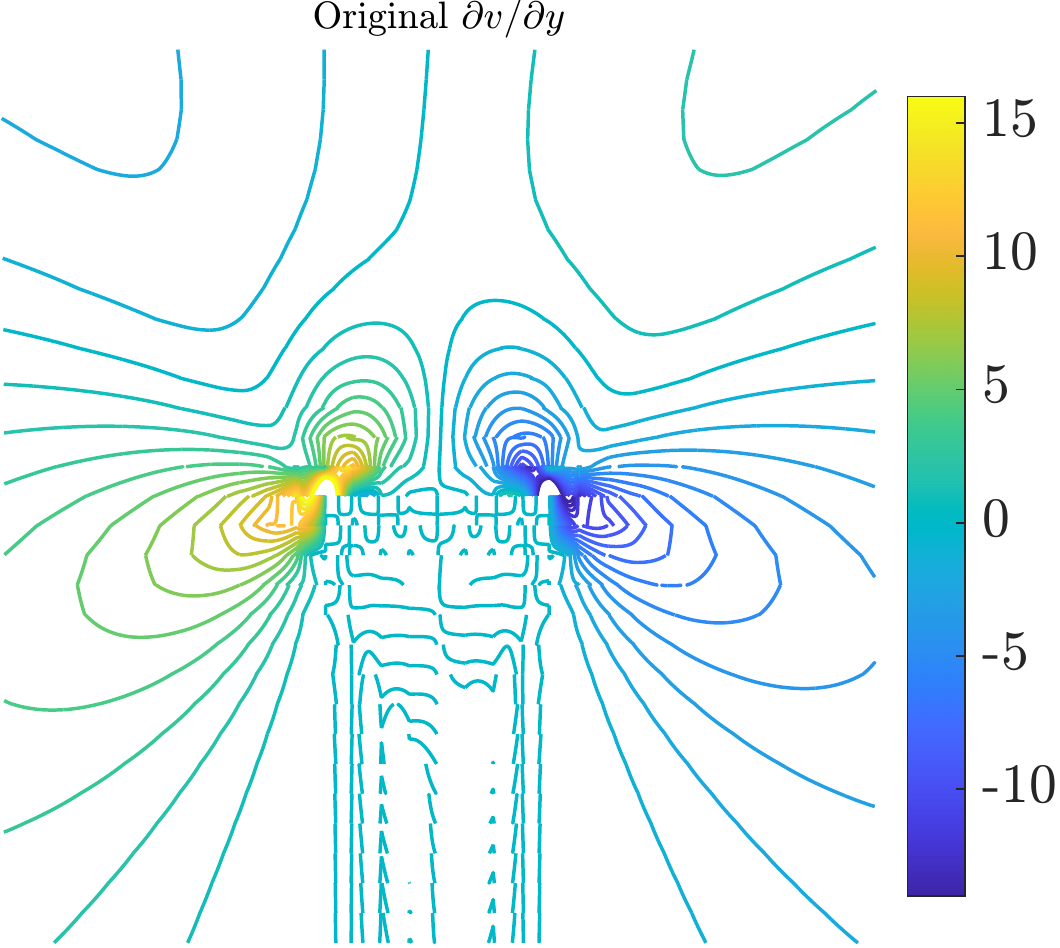}}
	\subfloat[Recovered $\dfrac{\partial v_2}{\partial y}$]{\includegraphics[width=0.25\textwidth]{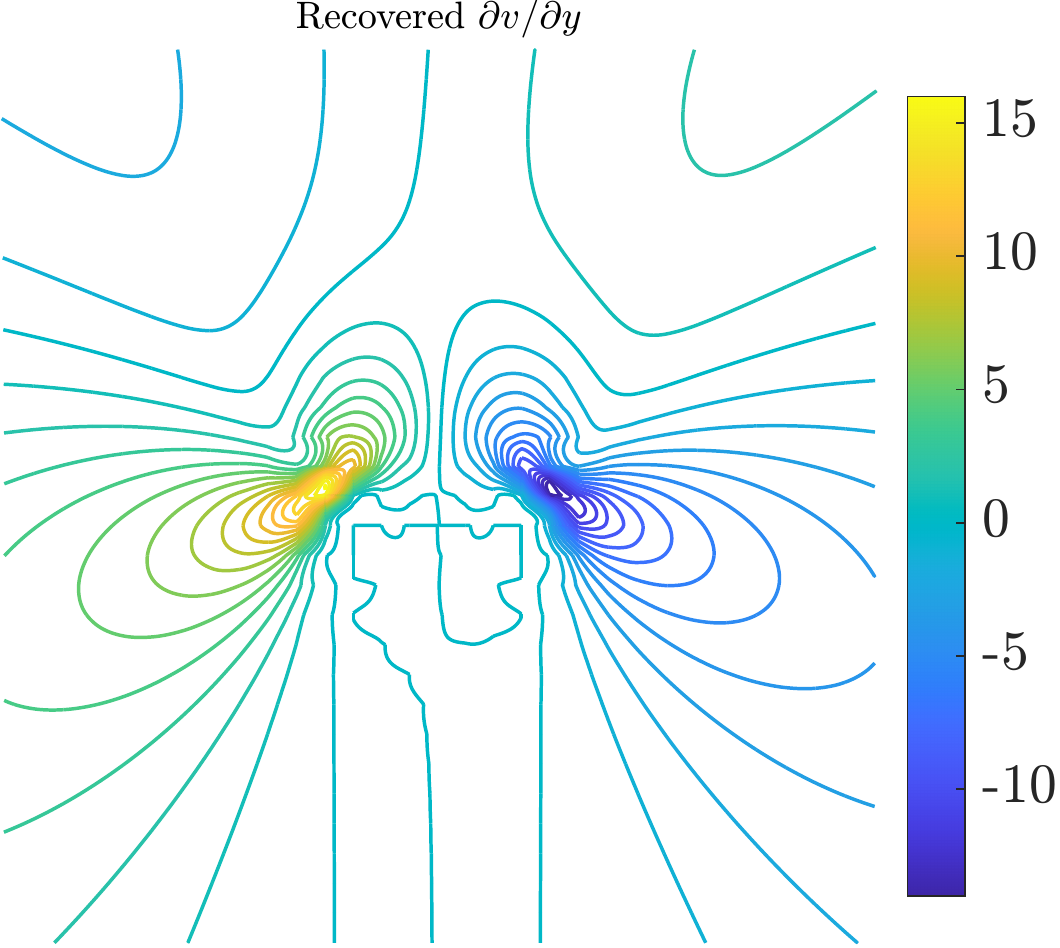}}
	\caption{Effect of the superconvergent patch recovery technique on velocity derivatives.}
	\label{fig:orgnl_vs_rcvrd_velocity_drvtv}
\end{figure}

To illustrate the effect of using the superconvergent patch recovery technique, Fig. \ref{fig:orgnl_vs_rcvrd_velocity_drvtv} shows the original vs. recovered $\partial v_1 / \partial x$ and $\partial v_2 / \partial y$ derivatives at a square area with coordinates \textcolor{black}{ $x = [ 0.9 \ \text{m}:1.1 \ \text{m} ]$ and $y = [ 0.4 \ \text{m}:0.6 \ \text{m} ]$ }of the \textcolor{black}{column} in a channel test problem, defined in section \ref{sec:tst_prblm}. In solving the fluid analysis problem and extracting these derivatives, we utilized \textcolor{black}{the default values described in section \ref{sec:tst_prblm} in addition to the following parameters; $\mu = 1$ Pa$\cdot$s, $p_\alpha = 5.25e-6$, $p_E = 1$, and $p_\mathit{\Upsilon} = 1$. An initial guess of 0.1 is used for all elements in the design domain (i.e. $\rho = 0.1$ in $\mathit{\Omega}_d$).} It can be clearly seen that using the superconvergent patch recovery technique restores the continuity of the velocity derivatives across the elemental boundaries. However, \textcolor{black}{the following phenomena} are observed:
\begin{enumerate}[a.]
	\item The polynomial degree representing the velocity derivatives is raised from bilinear in the original derivatives to biquadratic in the recovered ones.
	\item Some peaks in the original velocity derivatives are smoothed out in the recovered derivatives as evident at the two top corners of the non-design column, i.e. at \textcolor{black}{$(x,y) = (0.975,0.5) \text{ and } (1.025, 0.5)$ m.}
	\textcolor{black}{\item Due to the non-design column having almost zero velocities, the velocity derivatives within are almost vanishing and don't make a lot of sense.}
\end{enumerate}

Next, we move on to rederive the force coupling equations using the recovered velocity derivatives. The modified version of the surface integral form of force coupling with recovered derivatives is as follows:
\begin{align}
	\mathbf{F}^s_x & = \int_{-1}^{+1} \boldsymbol{\Psi} \left[ \left( - \mathbf{\Phi}^T \mathbf{\hat{p}} + 2 \mu \mathbf{\Psi}^T \overline{\mathbf{\hat{v}}_{1,x}} \right) n^f_x + \mu \left( \boldsymbol{\Psi}^T \overline{\mathbf{\hat{v}}_{1,y}} + \boldsymbol{\Psi}^T \overline{\mathbf{\hat{v}}_{2,x}} \right) n^f_y \right] |\mathbf{J}_s| \dd{\xi_s},	\label{eq:surf_force_x_recov}
	\\
	\mathbf{F}^s_y & = \int_{-1}^{+1} \boldsymbol{\Psi} \left[ \left( - \mathbf{\Phi}^T \mathbf{\hat{p}} + 2 \mu \mathbf{\Psi}^T \overline{\mathbf{\hat{v}}_{2,y}} \right) n^f_y + \mu \left( \boldsymbol{\Psi}^T \overline{\mathbf{\hat{v}}_{1,y}} + \boldsymbol{\Psi}^T \overline{\mathbf{\hat{v}}_{2,x}} \right) n^f_x \right] |\mathbf{J}_s| \dd{\xi_s}.
	\label{eq:surf_force_y_recov}
\end{align}

\noindent and the modified version of the volume integral form is as follows:
\begingroup
\allowdisplaybreaks
\begin{align}
	\begin{split}
		\mathbf{F}^s_x = \int_{-1}^{+1} \int_{-1}^{+1} & \Bigg[  \pdv{\boldsymbol{\Psi}}{x} \bigg( - \mathbf{\Phi} ^T \mathbf{\hat{p}} + 2 \mu \mathbf{\Psi} ^T \overline{\mathbf{\hat{v}}_{1,x}} \bigg)
		+ \pdv{\boldsymbol{\Psi}}{y} \bigg( \mu \boldsymbol{\Psi} ^T \overline{\mathbf{\hat{v}}_{1,y}} + \mu \boldsymbol{\Psi} ^T \overline{\mathbf{\hat{v}}_{2,x}} \bigg) \\ 
		& + \boldsymbol{\Psi} \bigg( - \pdv{\mathbf{\Phi}}{x} ^T \mathbf{\hat{p}} \, + 
		2 \mu \pdv{\boldsymbol{\Psi}}{x} ^T \overline{\mathbf{\hat{v}}_{1,x}} + \mu \pdv{\boldsymbol{\Psi}}{y} ^T \overline{\mathbf{\hat{v}}_{1,y}} + \mu \pdv{\boldsymbol{\Psi}}{y} ^T \overline{\mathbf{\hat{v}}_{2,x}} \bigg) \Bigg] |\mathbf{J}| \dd{\xi} \dd{\eta},
		\label{eq:vol_force_x_recov}
	\end{split} \\
	\begin{split}
		\mathbf{F}^s_y = \int_{-1}^{+1} \int_{-1}^{+1} & \Bigg[ \pdv{\boldsymbol{\Psi}}{x} \bigg( \mu \boldsymbol{\Psi} ^T \overline{\mathbf{\hat{v}}_{1,y}} + \mu \boldsymbol{\Psi} ^T \overline{\mathbf{\hat{v}}_{2,x}} \bigg) + \pdv{\boldsymbol{\Psi}}{y} \bigg( - \mathbf{\Phi}^T \mathbf{\hat{p}} + 2 \mu \mathbf{\Psi} ^T \overline{\mathbf{\hat{v}}_{2,y}} \bigg) \\
		& + \boldsymbol{\Psi} \bigg( \mu \pdv{\boldsymbol{\Psi}}{x} ^T \overline{\mathbf{\hat{v}}_{1,y}} + 
		\mu \pdv{\boldsymbol{\Psi}}{x} ^T \overline{\mathbf{\hat{v}}_{2,x}} - \pdv{\mathbf{\Phi}}{y} ^T \mathbf{\hat{p}} + 2 \mu \pdv{\boldsymbol{\Psi}}{y} ^T \overline{\mathbf{\hat{v}}_{2,y}} \bigg) \Bigg] |\mathbf{J}| \dd{\xi} \dd{\eta}.
		\label{eq:vol_force_y_recov}
	\end{split}
\end{align}
\endgroup

\subsection{Sensitivity Analysis using the Recovered Velocity Derivatives}
\label{ssec:snst_anls_recovery}
As detailed in the discussion on the sensitivity analysis (i.e. subsection \ref{ssec:snst_anls_general}), the derivatives of the solid governing equations w.r.t. the fluid state variables are needed as a part of the $\partial \mathbb{R} / \partial \mathbf{r}$ term. Before we discuss the sensitivity analysis any further, we share a thought on the discontinuities in the forces sensitivities.

Even though the derivatives of the forces w.r.t. velocities are independent of the nodal velocity values, there are still discontinuities at the elemental boundaries. By taking a closer look at any velocity derivative, take for instance $\pdv{v_1}/{x} = \pdv{\mathbf{\Psi}}/{x} ^T \mathbf{\hat{v}}_1$, evaluated at an elemental boundary from its two sides, the discontinuities arise due to a combination of \textbf{(i)} differing velocity components in $\mathbf{\hat{v}}_1$\footnote{Note that a velocity derivative, such as $\pdv{v_1}/{x} = \pdv{\mathbf{\Psi}}/{x} ^T \mathbf{\hat{v}}_1$, evaluated at an elemental boundary has non-zero components in $\pdv{\mathbf{\Psi}}/{x}$ not belonging to this particular boundary. Hence, $\pdv{v_1}/{x}$ includes contributions from nodes not belonging to this particular boundary.}, and \textbf{(ii)} differing shape functions derivatives, i.e. components of $\pdv{\mathbf{\Psi}}/{x}$. With the nodal velocities missing from the forces sensitivities, we focus on the discontinuities in the shape functions derivatives. The reason for discontinuities in the shape functions derivatives lies in the irregularity of the mesh and how this affects the Jacobian matrices of adjacent elements. In a highly regular mesh, these discontinuities almost disappear and the sensitivity analysis of the forces may be calculated directly by differentiating the original equations before recovering the velocity derivatives, i.e. by differentiating Eqs. \ref{eq:surf_force_x_org} and \ref{eq:surf_force_y_org} for the surface integral form and Eqs. \ref{eq:vol_force_x_org} and \ref{eq:vol_force_y_org} for the volume integral form. However, in the interest of generality in this work, we derive the sensitivity analysis equations after recovering the velocity derivatives.

The sensitivity analysis of the force equations w.r.t. velocities in the surface integral form is as follows:
\begingroup
\allowdisplaybreaks
\begin{align}
	\pdv{\mathbf{F}^s_x}{\hat{\mathbf{v}}_1} & = \int_{-1}^{+1} \boldsymbol{\Psi} \left[ \left( 2 \mu \mathbf{\Psi}^T \pdv{\overline{\mathbf{\hat{v}}_{1,x}}}{\hat{\mathbf{v}}_1} \right) n^f_x + \mu \left( \boldsymbol{\Psi}^T \pdv{\overline{\mathbf{\hat{v}}_{1,y}}}{\hat{\mathbf{v}}_1} \right) n^f_y \right] |\mathbf{J}_s| \dd{\xi_s},
	\\
	\pdv{\mathbf{F}^s_x}{\hat{\mathbf{v}}_2} & = \int_{-1}^{+1} \boldsymbol{\Psi} \left[  \mu \left( \boldsymbol{\Psi}^T \pdv{\overline{\mathbf{\hat{v}}_{2,x}}}{\hat{\mathbf{v}}_2} \right) n^f_y \right] |\mathbf{J}_s| \dd{\xi_s},
	\\
	\pdv{\mathbf{F}^s_y}{\hat{\mathbf{v}}_1} & = \int_{-1}^{+1} \boldsymbol{\Psi} \left[ \mu \left( \boldsymbol{\Psi}^T \pdv{\overline{\mathbf{\hat{v}}_{1,y}}}{\hat{\mathbf{v}}_1} \right) n^f_x \right] |\mathbf{J}_s| \dd{\xi_s},
	\\
	\pdv{\mathbf{F}^s_y}{\hat{\mathbf{v}}_2} & = \int_{-1}^{+1} \boldsymbol{\Psi} \left[ \left( 2 \mu \mathbf{\Psi}^T \pdv{\overline{\mathbf{\hat{v}}_{2,y}}}{\hat{\mathbf{v}}_2} \right) n^f_y + \mu \left( \boldsymbol{\Psi}^T \pdv{\overline{\mathbf{\hat{v}}_{2,x}}}{\hat{\mathbf{v}}_2} \right) n^f_x \right] |\mathbf{J}_s| \dd{\xi_s}.
\end{align}
\endgroup

\noindent while the volume integral form is as follows:
\begingroup
\allowdisplaybreaks
\begin{align}
	\begin{split}
		\pdv{\mathbf{F}^s_x}{\mathbf{\hat{\mathbf{v}}}_1} = \int_{-1}^{+1} \int_{-1}^{+1} & \Bigg[  \pdv{\boldsymbol{\Psi}}{x} \bigg( 2 \mu \mathbf{\Psi} ^T \pdv{\overline{\mathbf{\hat{v}}_{1,x}}}{\mathbf{\hat{\mathbf{v}}}_1} \bigg)
		+ \pdv{\boldsymbol{\Psi}}{y} \bigg( \mu \boldsymbol{\Psi} ^T \pdv{\overline{\mathbf{\hat{v}}_{1,y}}}{\mathbf{\hat{\mathbf{v}}}_1} \bigg) \\ 
		& + \boldsymbol{\Psi} \bigg( 2 \mu \pdv{\boldsymbol{\Psi}}{x} ^T \pdv{\overline{\mathbf{\hat{v}}_{1,x}}}{\mathbf{\hat{\mathbf{v}}}_1} + \mu \pdv{\boldsymbol{\Psi}}{y} ^T \pdv{\overline{\mathbf{\hat{v}}_{1,y}}}{\mathbf{\hat{\mathbf{v}}}_1} \bigg) \Bigg] |\mathbf{J}| \dd{\xi} \dd{\eta},
	\end{split} \\
	\begin{split}
		\pdv{\mathbf{F}^s_x}{\mathbf{\hat{\mathbf{v}}}_2} = \int_{-1}^{+1} \int_{-1}^{+1} & \Bigg[ \pdv{\boldsymbol{\Psi}}{y} \bigg( \mu \boldsymbol{\Psi} ^T \pdv{\overline{\mathbf{\hat{v}}_{2,x}}}{\mathbf{\hat{\mathbf{v}}}_2} \bigg) + \boldsymbol{\Psi} \bigg( \mu \pdv{\boldsymbol{\Psi}}{y} ^T \pdv{\overline{\mathbf{\hat{v}}_{2,x}}}{\mathbf{\hat{\mathbf{v}}}_2} \bigg) \Bigg] |\mathbf{J}| \dd{\xi} \dd{\eta},
	\end{split} \\
	\begin{split}
		\pdv{\mathbf{F}^s_y}{\mathbf{\hat{\mathbf{v}}}_1} = \int_{-1}^{+1} \int_{-1}^{+1} & \Bigg[ \pdv{\boldsymbol{\Psi}}{x} \bigg( \mu \boldsymbol{\Psi} ^T \pdv{\overline{\mathbf{\hat{v}}_{1,y}}}{\mathbf{\hat{\mathbf{v}}}_1} \bigg) + \boldsymbol{\Psi} \bigg( \mu \pdv{\boldsymbol{\Psi}}{x} ^T \pdv{\overline{\mathbf{\hat{v}}_{1,y}}}{\mathbf{\hat{\mathbf{v}}}_1} \bigg) \Bigg] |\mathbf{J}| \dd{\xi} \dd{\eta},
	\end{split} \\
	\begin{split}
		\pdv{\mathbf{F}^s_y}{\mathbf{\hat{\mathbf{v}}}_2} = \int_{-1}^{+1} \int_{-1}^{+1} & \Bigg[ \pdv{\boldsymbol{\Psi}}{x} \bigg( \mu \boldsymbol{\Psi} ^T \pdv{\overline{\mathbf{\hat{v}}_{2,x}}}{\mathbf{\hat{\mathbf{v}}}_2} \bigg) + \pdv{\boldsymbol{\Psi}}{y} \bigg( 2 \mu \mathbf{\Psi} ^T \pdv{\overline{\mathbf{\hat{v}}_{2,y}}}{\mathbf{\hat{\mathbf{v}}}_2} \bigg) \\
		& + \boldsymbol{\Psi} \bigg( \mu \pdv{\boldsymbol{\Psi}}{x} ^T \pdv{\overline{\mathbf{\hat{v}}_{2,x}}}{\mathbf{\hat{\mathbf{v}}}_2} + 2 \mu \pdv{\boldsymbol{\Psi}}{y} ^T \pdv{\overline{\mathbf{\hat{v}}_{2,y}}}{\mathbf{\hat{\mathbf{v}}}_2} \bigg) \Bigg] |\mathbf{J}| \dd{\xi} \dd{\eta}.
	\end{split}
\end{align}
\endgroup

The derivatives of a recovered nodal value  $\overline{\hat{v}_{i,j}}$ w.r.t. the nodal velocities of a specific element $\mathbf{\hat{\mathbf{v}}}_i^e$ are calculated, in vector form, as follows:
\begin{equation}
	\begin{split}
		\pdv{\overline{\hat{v}_{i,j}}}{\mathbf{\hat{v}}_i^e} & = \frac{1}{m} \sum_{l=1}^{m} \Bigg[  {\mathbf{P}_\text{O}}^T  \pdv{\mathbf{a}_l}{\mathbf{\hat{v}}_i} \Bigg] \\
		& = \frac{1}{m} \sum_{l=1}^{m} \Bigg[ {\mathbf{P}_\text{O}}^T \bigg[ \sum^{n} {\mathbf{P}_\times} \ {\mathbf{P}_\times}^T \bigg]^{-1} \bigg( \sum^{q} \pdv{\mathbf{\Psi}}{x_j} \Bigr|_{\times} {\mathbf{P}_\times}^T \bigg)_e \Bigg]_l.
	\end{split} \label{eq:drvtv_spr_vlct}
\end{equation}

\noindent with no summation on $i$ in the LHS term. $m$ is the number of patches sharing a particular recovered node, $n$ is the number of superconvergent points within a particular patch $l$, and $q$ is the number of superconvergent Gaussian points within a particular element $e$ which is 4 for \textcolor{black}{Q2Q1} elements.

Note that the use of superconvergent patch recovery extended the dependency between the adjacent elements' nodes, further complicating the sensitivity analysis. However, fortunately, Eq. \ref{eq:drvtv_spr_vlct} is independent of nodal velocity values, hence it need not be updated for new state variables. Nevertheless, in cases with mesh deformation, the equation needs to be solved after any changes in nodal $x$ and $y$ coordinates.

\subsection{Verification of the Sensitivity Analysis}
\label{ssec:verification_of_snst}

\begin{figure}[b]
	\captionsetup{font={color=black}}
	\centering
	\includegraphics[width=0.5\textwidth]{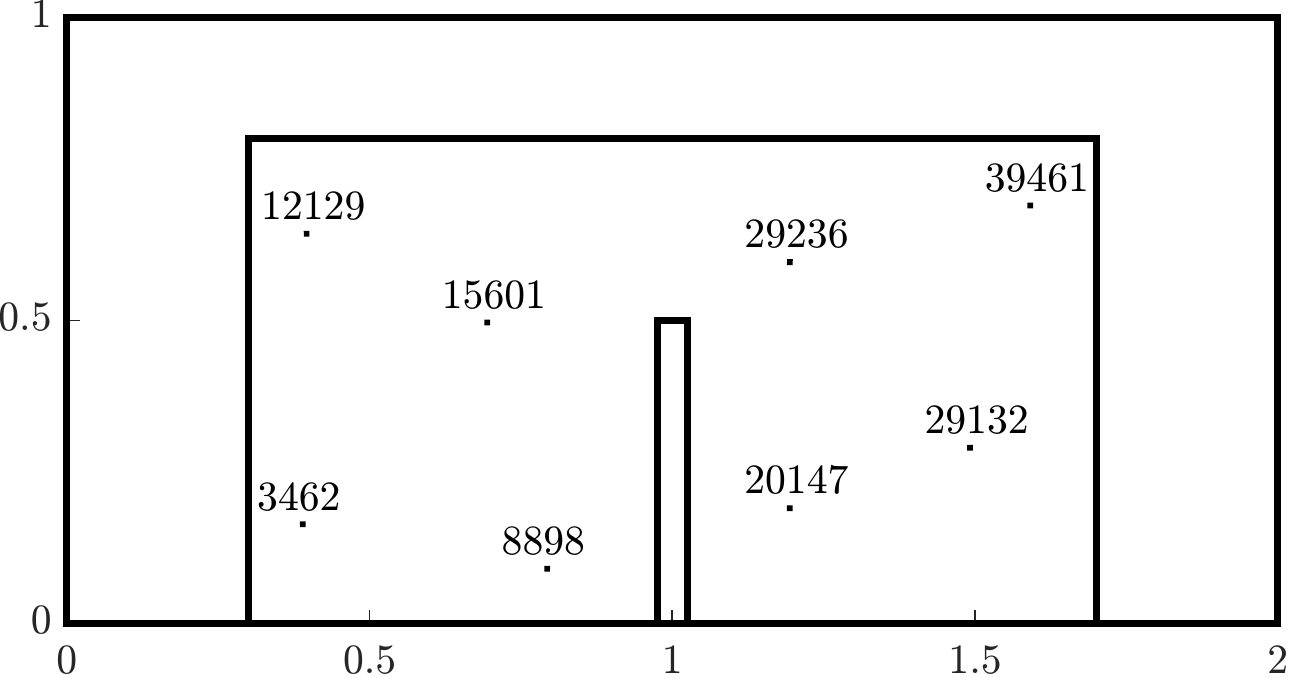}
	\caption{Location and number designations of the finite elements selected for sensitivity verification.}
	\label{fig:lctn_of_snst_elmnts}
\end{figure}

\begin{table}[t!]
	\centering
	\captionsetup{font={color=black}}
	\captionsetup{width=0.82\textwidth}
	\caption{Verification of analytically-derived sensitivities using the complex-step derivative approximation method.}
	\label{tab:snst_vrfctn}
	\color{black}
	\begin{tabular}{@{}crrrr@{}} \toprule
		\color{black}
		& \multicolumn{4}{c}{Pressure Coupling} \\ 
		& \multicolumn{4}{c}{\cbox{Total Stress Coupling}} \\ 
		\cmidrule{1-5}
		Element & Imaginary Step & Obj. Fun. Imaginary Part & Analytical Sensitivity & Normalized Error \\ \midrule
		12129 & $+1e-10$ & $+8.137062143103e-14$ & $+8.137062142673e-04$ & $+5.28e-09\%$ \\
		& \cbox{$+1e-10$} & \cbox{$+8.912673735197e-14$} & \cbox{$+8.912673735248e-04$} & \cbox{$-5.72e-10\%$} \\[6pt]
		3462 & $+1e-10$ & $-6.415904186362e-14$ & $-6.415904186347e-04$ & $+2.34e-10\%$ \\
		& \cbox{$+1e-10$} & \cbox{$-7.133040232407e-14$} & \cbox{$-7.133040232431e-04$} & \cbox{$-3.36e-10\%$} \\[6pt]
		15601 & $-1e-10$ & $-5.915480351527e-14$ & $+5.915480351356e-04$ & $+2.89e-09\%$ \\
		& \cbox{$-1e-10$} & \cbox{$-6.090842047842e-14$} & \cbox{$+6.090842047873e-04$} & \cbox{$-5.09e-10\%$} \\[6pt]
		8898 & $-1e-10$ & $+1.718044899909e-13$ & $-1.718044899907e-03$ & $+1.16e-10\%$ \\
		& \cbox{$-1e-10$} & \cbox{$+1.815291681449e-13$} & \cbox{$-1.815291681448e-03$} & \cbox{$+5.51e-11\%$} \\[6pt]
		29236 & $+1e-10$ & $+1.783493810518e-13$ & $+1.783493810553e-03$ & $-1.96e-09\%$ \\
		& \cbox{$+1e-10$} & \cbox{$+1.882944445386e-13$} & \cbox{$+1.882944445355e-03$} & \cbox{$+1.65e-09\%$} \\[6pt]
		20147 & $+1e-10$ & $-1.358101868978e-13$ & $-1.358101868974e-03$ & $+2.95e-10\%$ \\
		& \cbox{$+1e-10$} & \cbox{$-1.450049498019e-13$} & \cbox{$-1.450049498008e-03$} & \cbox{$+7.59e-10\%$} \\[6pt]
		39461 & $-1e-10$ & $-9.033751176788e-14$ & $+9.033751176807e-04$ & $-2.10e-10\%$ \\
		& \cbox{$-1e-10$} & \cbox{$-1.019007554562e-13$} & \cbox{$+1.019007554518e-03$} & \cbox{$+4.32e-09\%$} \\[6pt]
		29132 & $+1e-10$ & $-4.193677350665e-14$ & $-4.193677350532e-04$ & $+3.17e-09\%$ \\
		& \cbox{$+1e-10$} & \cbox{$-4.553959134022e-14$} & \cbox{$-4.553959133947e-04$} & \cbox{$+1.65e-09\%$} \\	\bottomrule
	\end{tabular}
\end{table}

In this subsection, we utilize the complex-step derivative approximation method (cf. \ \cite{Martins2003}) to verify the sensitivities calculated for the pressure as well as total stress coupling cases. The complex-step derivative approximation method is implemented by introducing a small imaginary step $ih$ to a particular design variable $\rho_i$ at a time, and calculating the approximate derivative of the objective function $f$ w.r.t. that particular design variable $\rho_i$ as follows:
\begin{equation}
	\pdv{f}{\rho_i} \approx \frac{\text{Im} \left[ f(\rho_i + ih) \right] }{h} 
\end{equation} 

\noindent where $\text{Im}(f)$ indicates the imaginary part of $f$.

Utilizing the \textcolor{black}{column} in a channel test problem defined in section \ref{sec:tst_prblm}, we randomly sampled eight finite elements within the design domain and calculated the sensitivities at these elements using the complex-step derivative approximation method against the analytically-derived sensitivities discussed in subsection \ref{ssec:snst_anls_recovery}. Figure \ref{fig:lctn_of_snst_elmnts} shows the location of the finite elements of interest within the design domain. For the results in table \ref{tab:snst_vrfctn}, we utilized \textcolor{black}{the parameters used in subsection 5.1 without filtering.} The error in table \ref{tab:snst_vrfctn} is calculated by subtracting the analytical sensitivity from the complex sensitivity then dividing over the analytical sensitivity. As can be noted, the errors are within a reasonable range attesting to the validity of the derived analytical sensitivity.

% \textcolor{black}{It can be noted that the finite elements at the fluid-structure interface have more pronounced errors compared to the internal ones.}

\section{Numerical Experiments}
\label{sec:numerical_exp}

In this section, we present some numerical results through solving the test problem described in section \ref{sec:tst_prblm}. Our main intention in the following discussion is to demonstrate the difference between the pressure vs. total stress force coupling. The \textbf{m}ethod of \textbf{m}oving \textbf{a}symptotes (MMA) is used to update the design variables \cite{Svanberg1987}. In this work, all the derivatives are derived analytically - by hand - and the entire TOFSI problem is coded and solved in MATLAB.

In general, TOFSI with pressure coupling uses less computational resources when compared with TOFSI with total stress coupling. As mentioned in the last paragraph in section \ref{ssec:snst_anls_recovery}, in the absence of mesh deformations, the derivatives of the fluidic forces w.r.t. the state variables are calculated only once in the first iteration. In the subsequent iterations, the pressure-coupling filtering function is updated with the new design variables and multiplied by the derivatives of the fluidic forces w.r.t. the state variables. Based on our coding practices, we estimate the increase in computational costs to be below \textcolor{black}{10\%}.

\color{black}
The column in a channel problem (described in section \ref{sec:tst_prblm}) is solved with the objective function of minimizing the compliance of the structure at a specific volume fraction as follows:
\begin{equation}
	\begin{aligned}
		\text{minimize:} & \quad f = \sum \mathbf{\hat{u}} ^T \mathbf{k}^s \mathbf{\hat{u}}, \\[4pt]
		\text{subject to:} & \quad \mathbf{0} \leq \boldsymbol{\uprho} \leq \boldsymbol{1} , \\[4pt]
		& \quad \sum_i V_i/V_0 \leq V_f .
	\end{aligned}
\end{equation}

\noindent where $V_i$ is the volume of finite element $i$, $V_0$ is the volume of the entire design domain, and the volume fraction $V_f$ is set to 0.1 following \ \citep[p.~976]{Lundgaard2018}. The volume fraction is further used as the initial guess for all the optimized designs. The problem is optimized for two Reynolds numbers; $Re = 1$ and 10, obtained by varying the fluid dynamic viscosity as discussed in section \ref{sec:tst_prblm} such that $\mu = 1$ and 0.1 Pa$\cdot$s, respectively.

\subsection{In Search of Discrete Designs}

\begin{figure}[b!]
	\centering
	\captionsetup{font={color=black}}
	\includegraphics[width=0.9\textwidth]{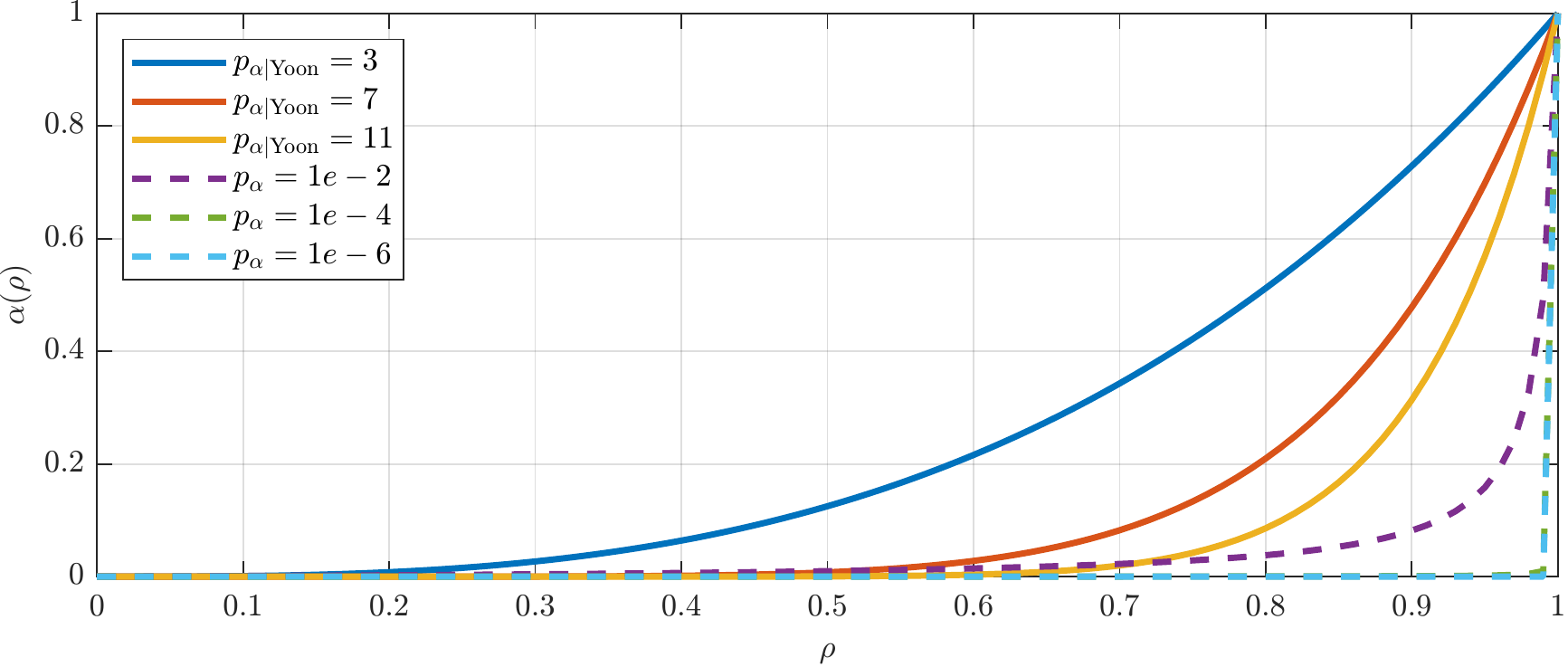}
	\caption{A comparison of two Brinkman penalization interpolation functions. Solid lines represent the interpolation function used by Yoon (Eq. 33 in \ \citep[p.~607]{Yoon2010}) while dashed lines represent Eq. \ref{eq:invrs_perm} in this work.  A unity $\alpha_\text{max}$ and zero $\alpha_\text{min}$ are used for demonstration purposes.}
	\label{fig:brnk_intrp_fnc}
\end{figure}

In \ \cite{Yoon2010}, Yoon used constant interpolation parameters to solve the original column in a channel problem. Specifically, he used $p_E = 3$, $p_\mathit{\Upsilon} = 3$, and $p_\alpha = 11$ for different Reynolds numbers, albeit the Brinkman penalization interpolation function was defined differently from Eq. \ref{eq:invrs_perm} and more similar to the typical SIMP equation. Figure \ref{fig:brnk_intrp_fnc} shows a comparison between the Brinkman interpolation function used by Yoon\cite{Yoon2010} and the interpolation function implemented in this work (Eq. \ref{eq:invrs_perm}), showing that the former is generally more sensitive to intermediate density elements. Later in \ \cite{Lundgaard2018}, Lundgaard et al. solved the modified column in a channel problem using the robust formulation\cite{Sigmund2009} with $p_E = 1$, $p_\mathit{\Upsilon} = 1$, and $p_\alpha$ calibrated for different Reynolds numbers. Lundgaard et al. \cite{Lundgaard2018} brought considerable attention to the effect of $p_\alpha$ on the optimized designs. 

Initially, we attempted to solve the problem using continuation on $p_E$, $p_\mathit{\Upsilon}$, and $p_\alpha$ in order to avoid the exacerbated computational costs associated with using the robust formulation. However, TOFSI problems proved to be notoriously non-convex and whatever computational costs were saved by avoiding the robust formulation, were consumed in the numerous trial-and-error attempts to calibrate the interpolation parameters. Moreover, the best combination we could find still contained considerable intermediate densities. Eventually, we settled on a combination of the robust formulation and continuation on $p_E$ and $p_\mathit{\Upsilon}$. From our perspective, obtaining close-to-discrete designs takes precedence over enhancing the objective function.

\begin{figure}[b!]
	\centering
	\captionsetup{font={color=black}}
	\captionsetup[subfigure]{font+={color=black}}
	\newcommand{\mywdth}{0.325\textwidth}
	\captionsetup[subfigure]{justification=centering}
	\subfloat[$\delta_{E|\mathit{\Upsilon}} = 1.00$, \\ \vspace{1pt} $p_\mathit{\Upsilon}$ = \{1.50, 2.00, 3.00, 4.00\}, \\ \vspace{1pt} $f = 0.014624, \ DM = 1.19 \%$.]{\includegraphics[width=\mywdth]{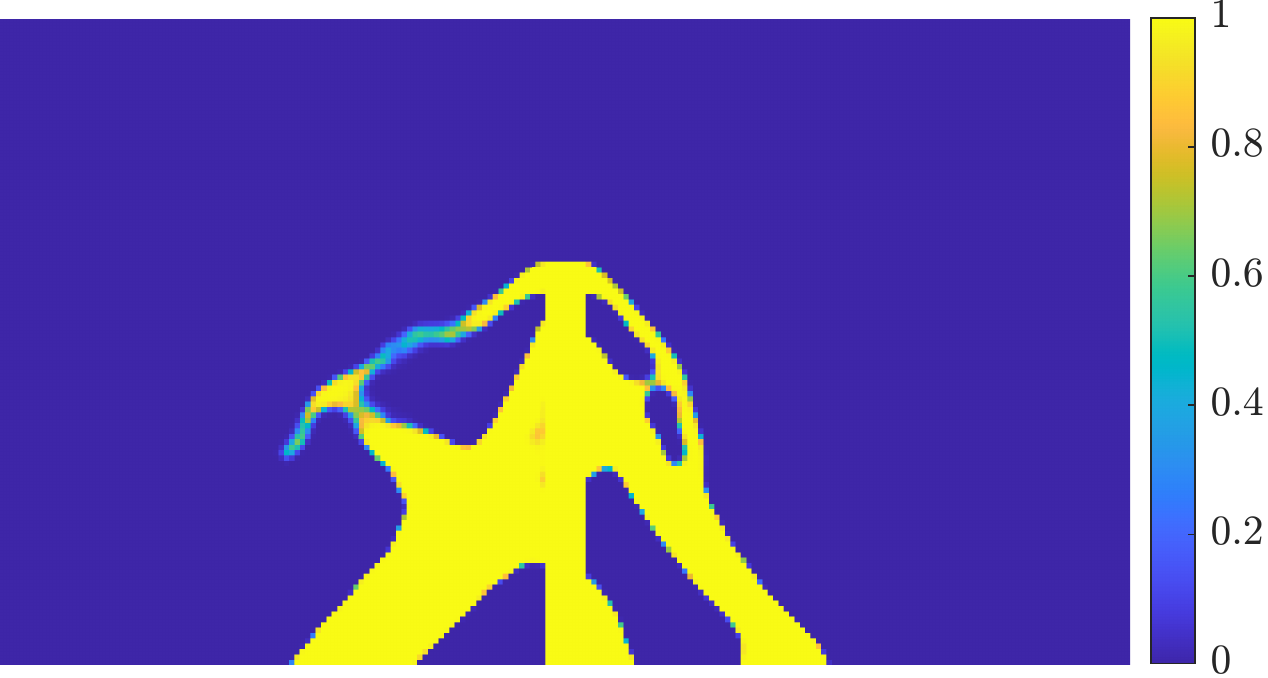}}
	\hspace{2pt}
	\subfloat[$\delta_{E|\mathit{\Upsilon}} = 1.25$, \\ \vspace{1pt} $p_\mathit{\Upsilon}$ = \{1.40, 1.80, 2.60, 3.40\}, \\ \vspace{1pt} $f = 0.014340, \ DM = 0.64 \%$.]{\includegraphics[width=\mywdth]{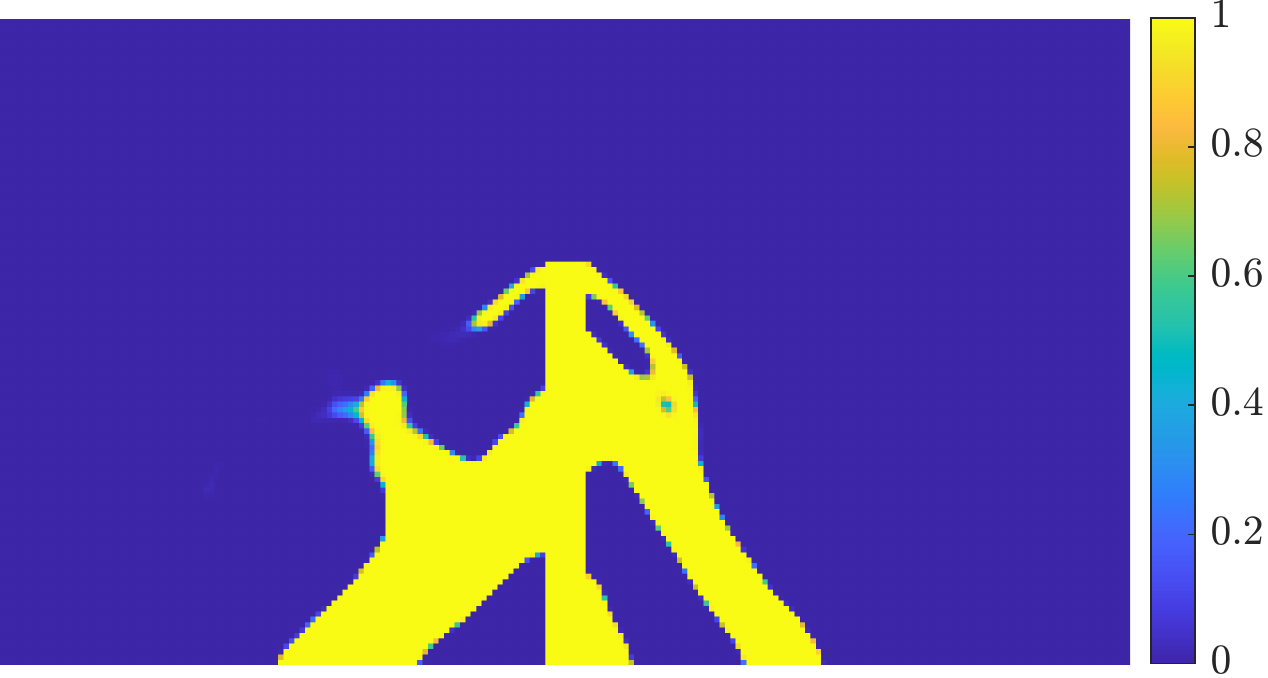}}
	\hspace{2pt}
	\subfloat[$\delta_{E|\mathit{\Upsilon}} = 1.50$, \\ \vspace{1pt} $p_\mathit{\Upsilon}$ = \{1.33, 1.67, 2.33, 3.00\}, \\ \vspace{1pt} $f = 0.011839, \ DM = 0.36 \%$.]{\includegraphics[width=\mywdth]{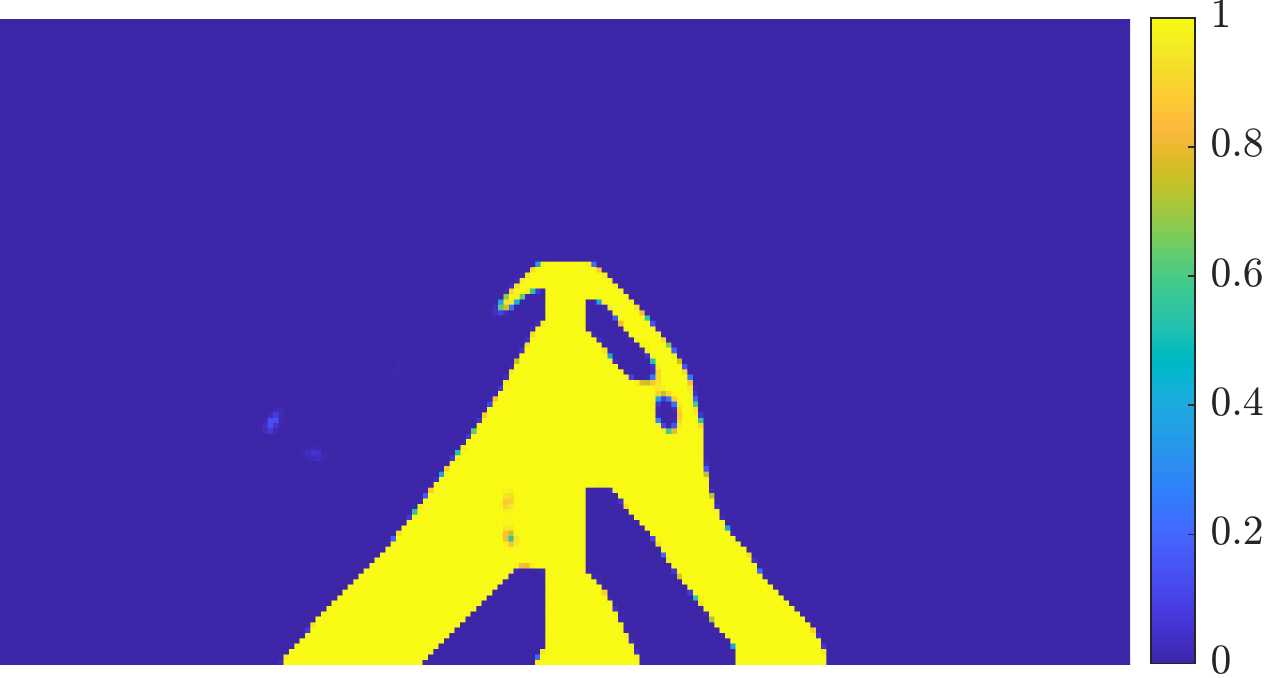}} \\
	\subfloat[$\delta_{E|\mathit{\Upsilon}} = 1.75$, \\ \vspace{1pt} $p_\mathit{\Upsilon}$ = \{1.29, 1.57, 2.14, 2.71\}, \\ \vspace{1pt} $f = 0.017344, \ DM = 0.70 \%$.]{\includegraphics[width=\mywdth]{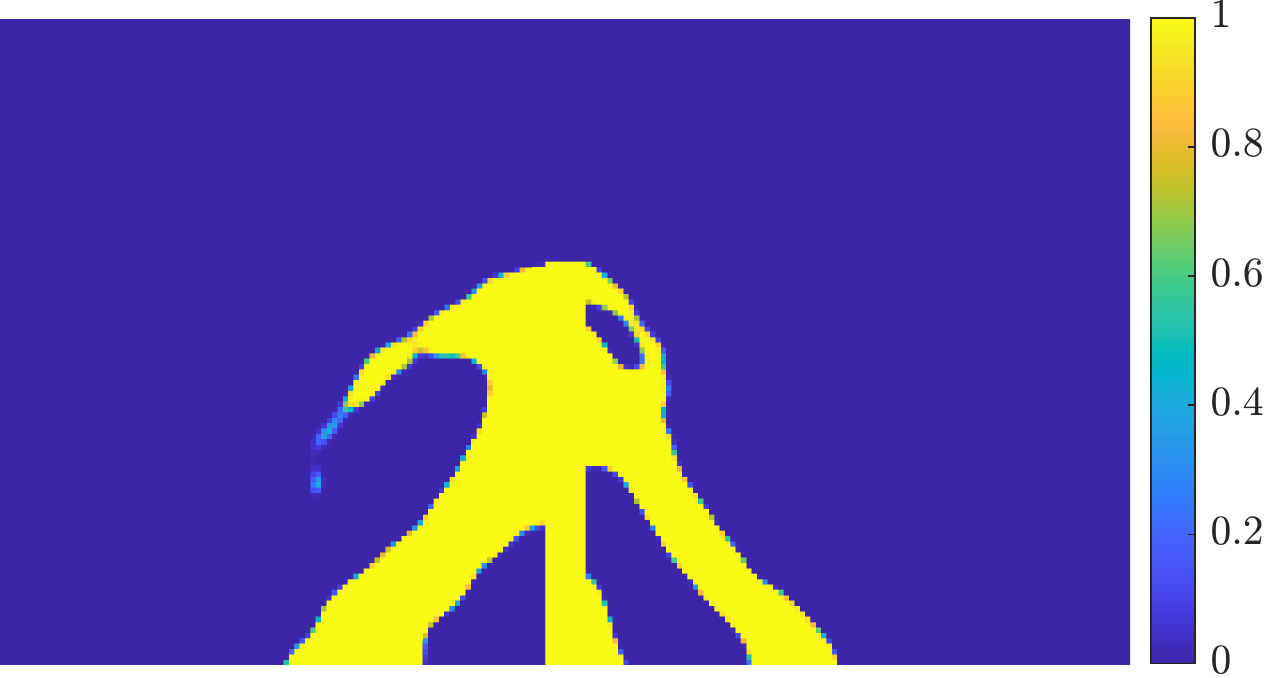}}
	\hspace{2pt}
	\subfloat[$\delta_{E|\mathit{\Upsilon}} = 2.00$, \\ \vspace{1pt} $p_\mathit{\Upsilon}$ = \{1.25, 1.50, 2.00, 2.50\}, \\ \vspace{1pt} $f = 0.011654, \ DM = 0.30 \%$. \label{ffig:5e}]{\includegraphics[width=\mywdth]{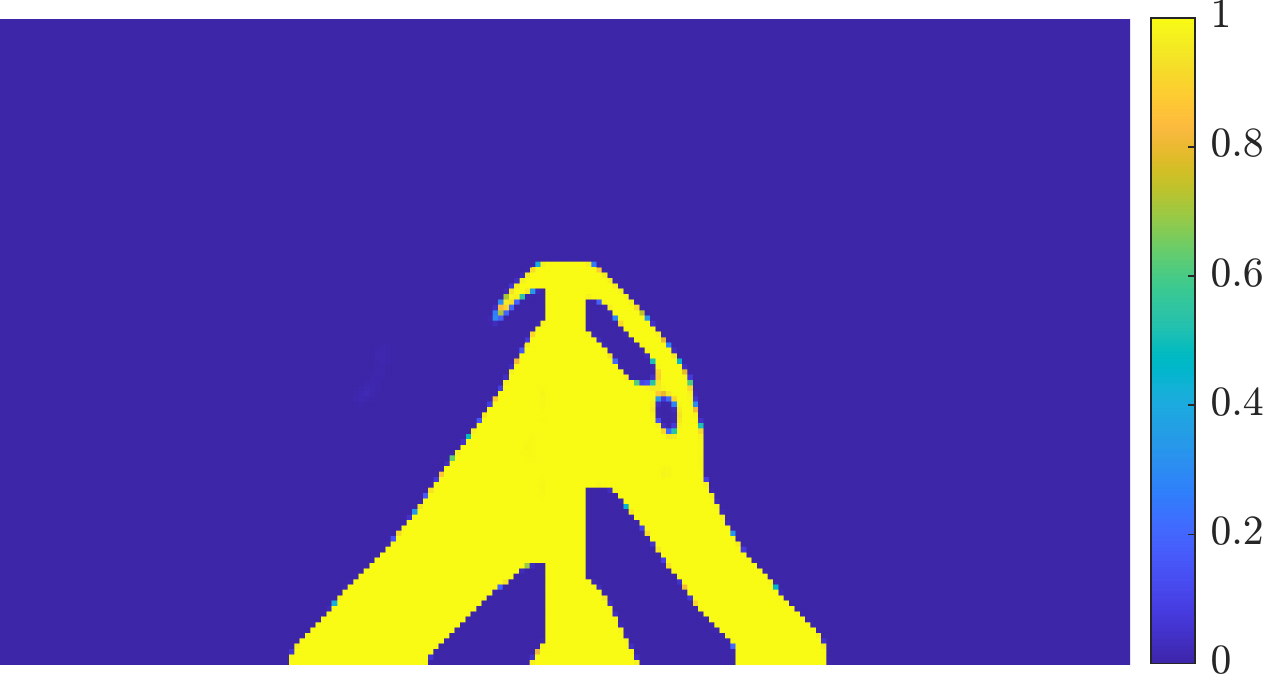}}
	\hspace{2pt}
	\subfloat[$\delta_{E|\mathit{\Upsilon}} = 2.25$, \\ \vspace{1pt} $p_\mathit{\Upsilon}$ = \{1.22, 1.44, 1.89, 2.33\}, \\ \vspace{1pt} $f = 0.011797, \ DM = 0.61 \%$.]{\includegraphics[width=\mywdth]{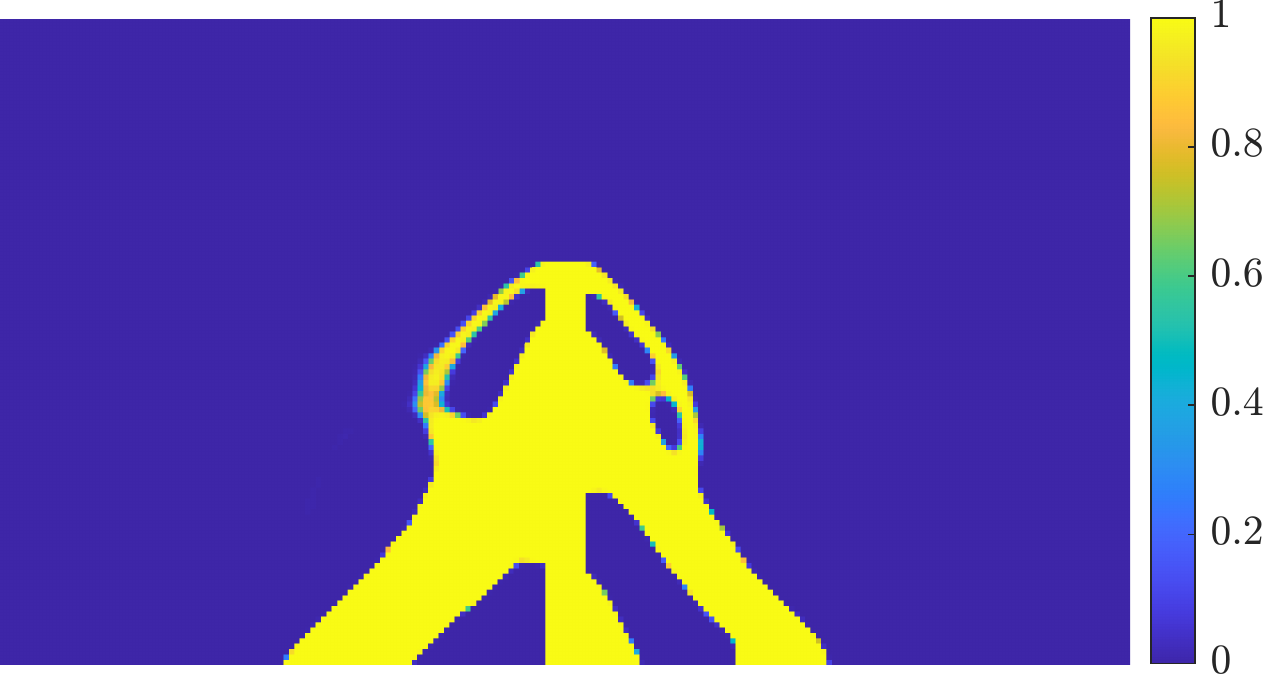}} \\
	\subfloat[$\delta_{E|\mathit{\Upsilon}} = 2.50$, \\ \vspace{1pt} $p_\mathit{\Upsilon}$ = \{1.20, 1.40, 1.80, 2.20\}, \\ \vspace{1pt} $f = 0.011683, \ DM = 0.33 \%$.]{\includegraphics[width=\mywdth]{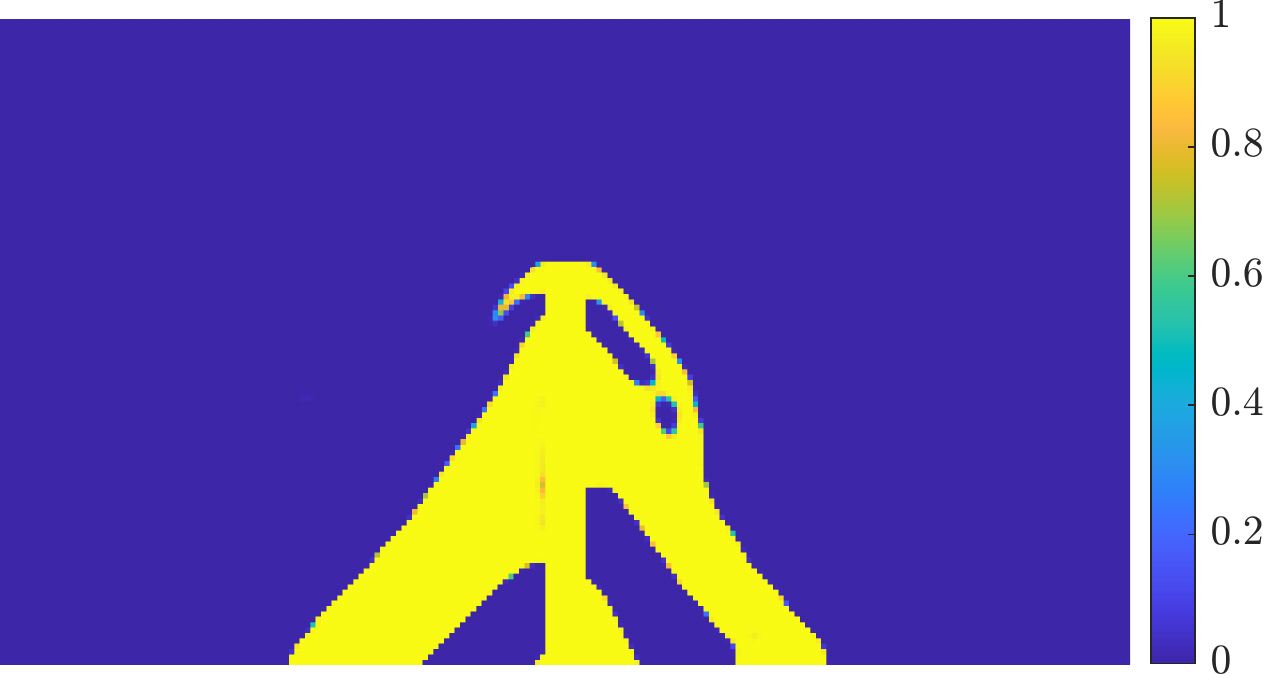}}
	\hspace{2pt}
	\subfloat[$\delta_{E|\mathit{\Upsilon}} = 2.75$, \\ \vspace{1pt} $p_\mathit{\Upsilon}$ = \{1.18, 1.36, 1.73, 2.09\}, \\ \vspace{1pt} $f = 0.011665, \ DM = 0.28 \%$.]{\includegraphics[width=\mywdth]{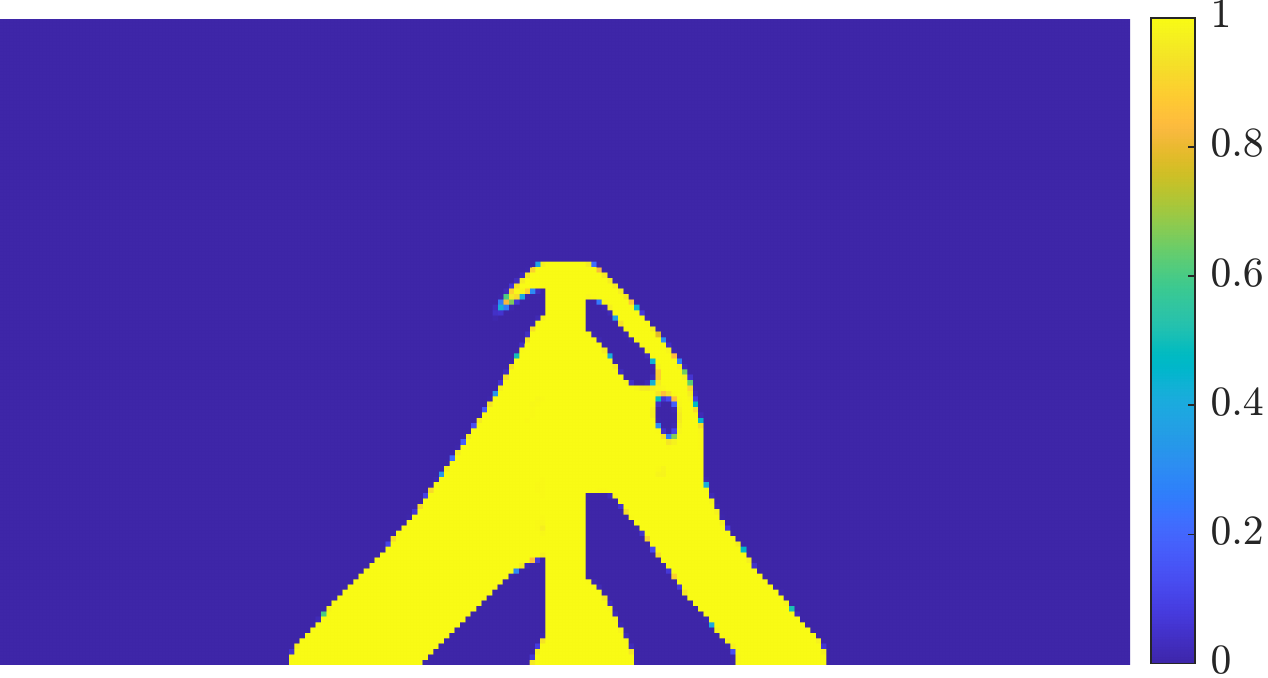}}
	\hspace{2pt}
	\subfloat[$\delta_{E|\mathit{\Upsilon}} = 3.00$, \\ \vspace{1pt} $p_\mathit{\Upsilon}$ = \{1.17, 1.33, 1.67, 2.00\}, \\ \vspace{1pt} $f = 0.011444, \ DM = 0.41
	\%$.]{\includegraphics[width=\mywdth]{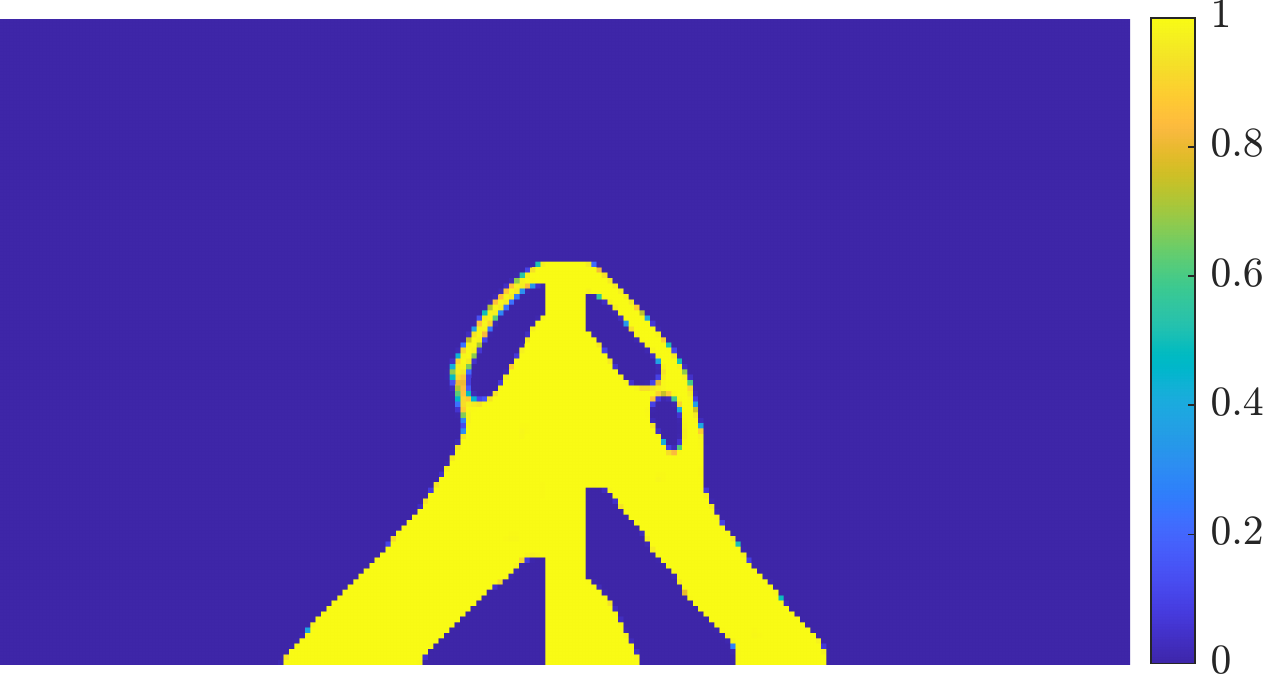}}
	\caption{Effect of $p_E$ vs $p_\mathit{\Upsilon}$ on optimized designs with  $p_\alpha = 5.25e-6$ at $Re = 1$ with pressure coupling. At iterations \{21, 41, 61, 81\}, $p_E$ is set to \{1.5, 2.0, 3.0, 4.0\} while $p_\mathit{\Upsilon}$ is calculated for different $\delta_{E|\mathit{\Upsilon}}$ as shown in the captions. $f$ is the objective function while $DM$ is the discreteness measure.}
	\label{fig:pe_vs_pu_5.25e-6_prs}
\end{figure}

To generate the following designs, the problem is started with $p_E = 1$, $p_\mathit{\Upsilon} = 1$, and $\beta = 4$. These parameters are updated such that at iterations \{21, 41, 61, 81\}, $p_E$ is raised by \{0.5, 0.5, 1.0, 1.0\}, $p_\mathit{\Upsilon}$ is raised by \{0.5/$\delta_{E|\mathit{\Upsilon}}$, 0.5/$\delta_{E|\mathit{\Upsilon}}$, 1.0/$\delta_{E|\mathit{\Upsilon}}$, 1.0/$\delta_{E|\mathit{\Upsilon}}$\}, and $\beta$ is set to \{8, 16, 32, 64\}. $\delta_{E|\mathit{\Upsilon}}$ is a parameter that is used to modify the penalization of $p_\mathit{\Upsilon}$ against that of $p_E$, hence it takes values greater than or equal to unity. For $\delta_{E|\mathit{\Upsilon}} = 1$, $p_E$ and $p_\mathit{\Upsilon}$ experience the same penalization. As $\delta_{E|\mathit{\Upsilon}}$ increases, $p_\mathit{\Upsilon}$ gets less penalized than $p_E$. For the robust formulation, we followed Eqs. 10-13 in \ \citep[p.~973]{Lundgaard2018} with projection thresholding values of $\eta_n = 0.5$, $\eta_d = 0.49$, and $\eta_e = 0.51$ for the nominal, dilated, and eroded designs\footnote{For some reason unknown to us, the projection thresholds of \{0.5, 0.3, 0.7\}  implemented in Lundgaard et al. \cite{Lundgaard2018} produced designs that didn't contain any structural details. Instead all material was concentrated around the non-design column as if to minimize the surface area, hence minimize the differences between the three designs.}, respectively. The relative filtering radius is set to 5.3, which is a bit high but is necessary to avoid thin members that tend to contain intermediate density elements even at high $\beta$ values.

To achieve this relatively low number of iterations per step, we used the default settings of the MMA code in its min/max form\citep[p.~3]{Svanberg2004} except for the following changes; \textbf{(i)} the move limit is set to 0.1, \textbf{(ii)} a positive offset of 1 is added to the objective function before inputting into the MMA code to increase the constraints violation, and \textbf{(iii)} the objective function (and its derivatives) is multiplied by 100 to raise the order of magnitude of the derivatives for the cases of $Re = 10$. These modifications generally increase the ``aggressiveness" of the MMA optimizer with a slight detriment to the optimality of the final design. In our numerical experiments, the maximum change in the design variables is generally above 0.05 and typically higher at the start and after each continuation update. The volume fraction constraint is applied on the dilated design with its value updated every 2 iterations to account for this increased aggressiveness, cf. Eq. 14 in \ \citep[p.~774]{Wang2011}.

\begin{figure}[b!]
	\newcommand{\mywdth}{0.325\textwidth}
	\centering
	\captionsetup{font={color=black}}
	\captionsetup[subfigure]{font+={color=black}}
	\captionsetup[subfigure]{justification=centering}
	\subfloat[$\delta_{E|\mathit{\Upsilon}} = 1.00$, \\ \vspace{1pt} $p_\mathit{\Upsilon}$ = \{1.50, 2.00, 3.00, 4.00\}, \\ \vspace{1pt} $f = 0.011405, \ DM = 0.42 \%$.]{\includegraphics[width=\mywdth]{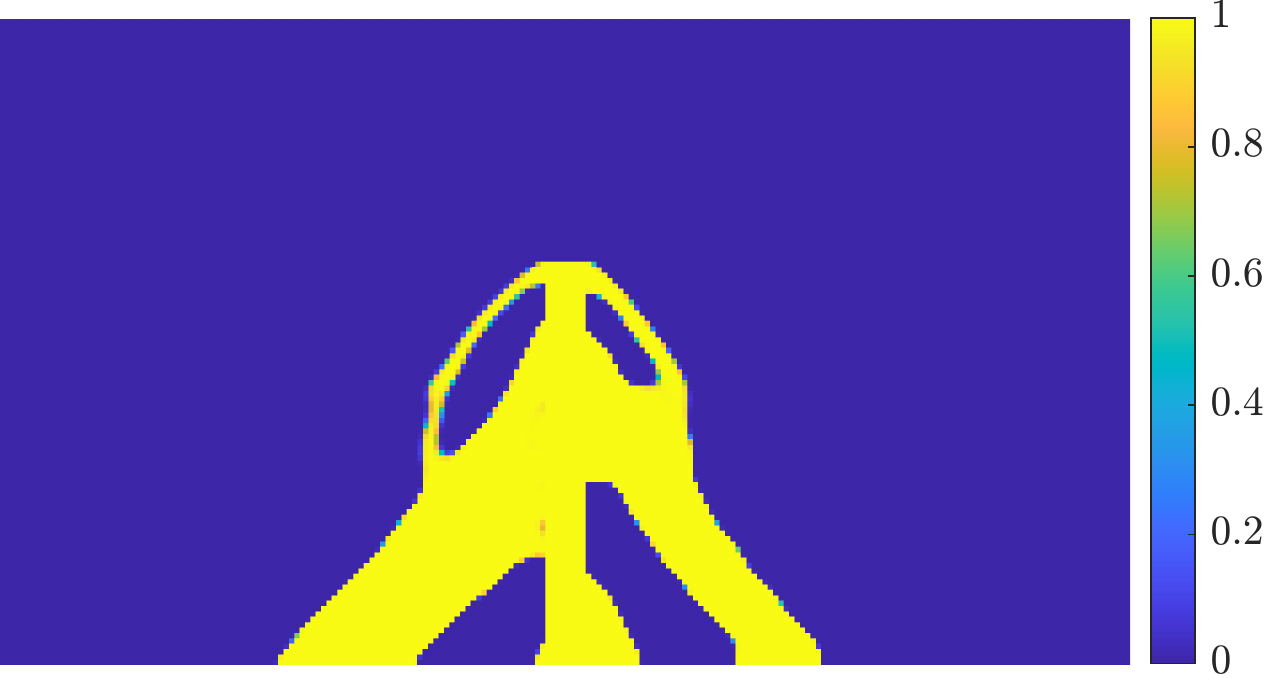}}
	\hspace{2pt}
	\subfloat[$\delta_{E|\mathit{\Upsilon}} = 1.50$, \\ \vspace{1pt} $p_\mathit{\Upsilon}$ = \{1.33, 1.67, 2.33, 3.00\}, \\ \vspace{1pt} $f = 0.011324, \ DM = 0.40 \%$.]{\includegraphics[width=\mywdth]{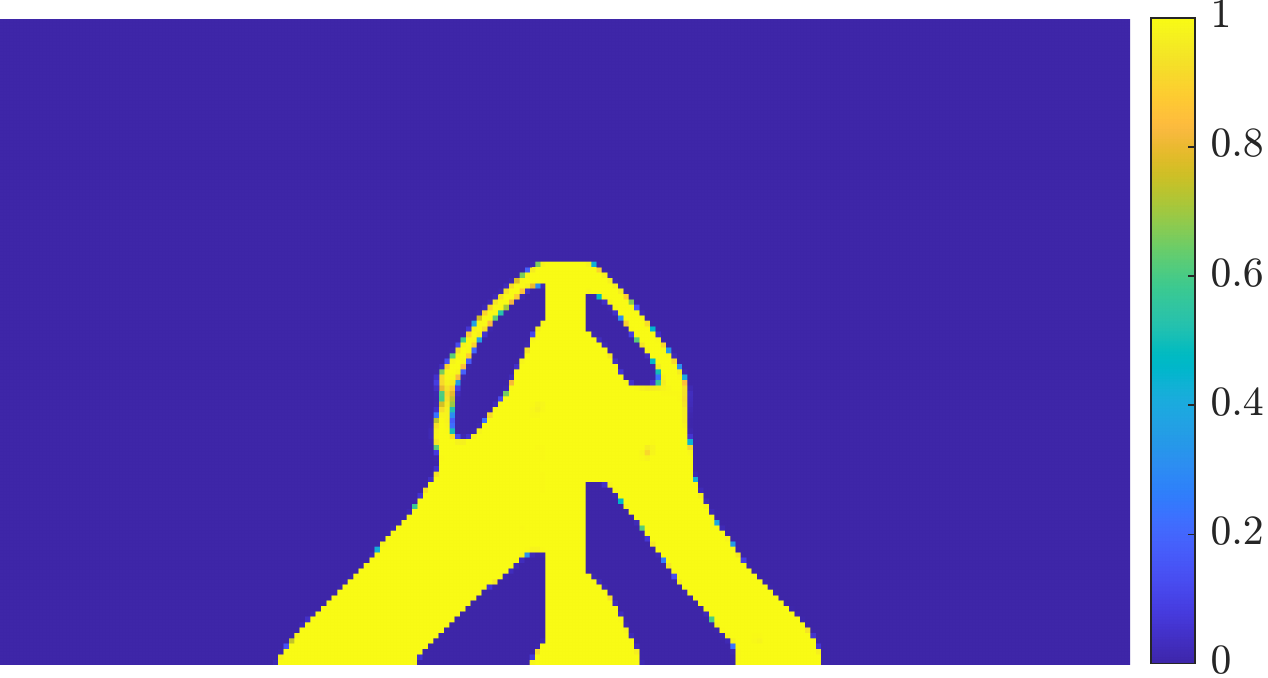}}
	\hspace{2pt}
	\subfloat[$\delta_{E|\mathit{\Upsilon}} = 2.00$, \\ \vspace{1pt} $p_\mathit{\Upsilon}$ = \{1.25, 1.50, 2.00, 2.50\}, \\ \vspace{1pt} $f = 0.011370, \ DM = 0.37 \%$. \label{ffig:6c}]{\includegraphics[width=\mywdth]{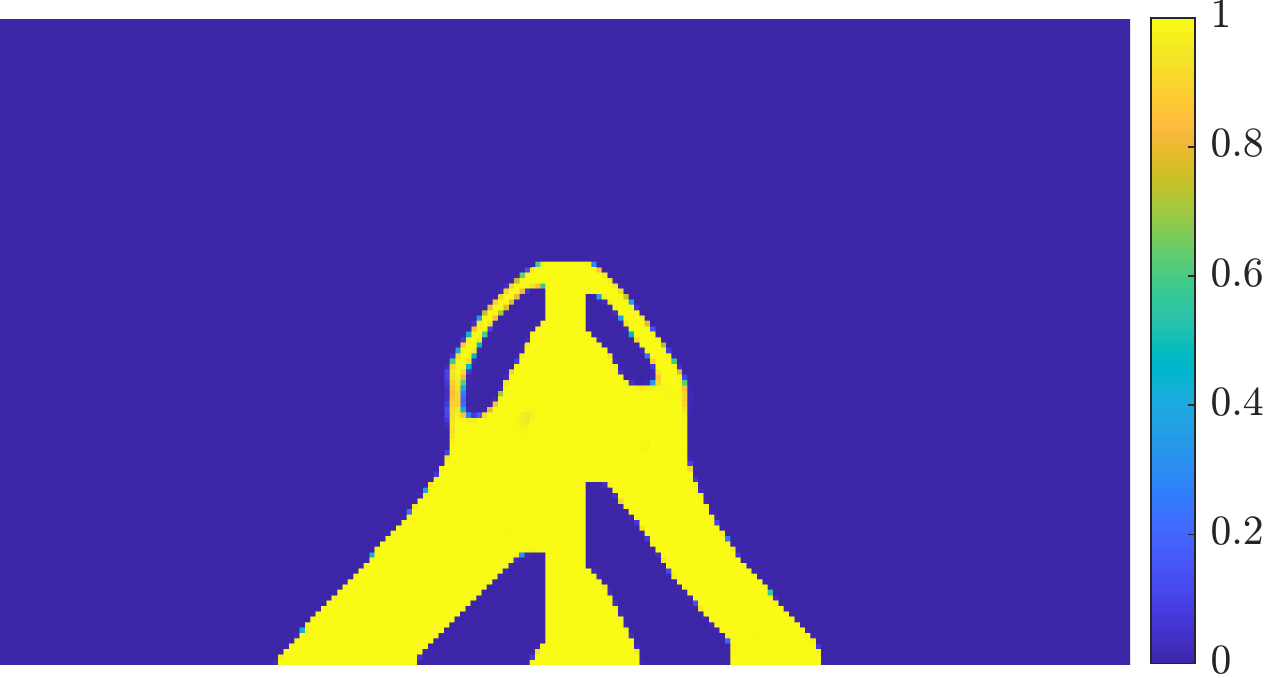}} \\
	\subfloat[$\delta_{E|\mathit{\Upsilon}} = 2.50$, \\ \vspace{1pt} $p_\mathit{\Upsilon}$ = \{1.20, 1.40, 1.80, 2.20\}, \\ \vspace{1pt} $f = 0.011394, \ DM = 0.34 \%$.]{\includegraphics[width=\mywdth]{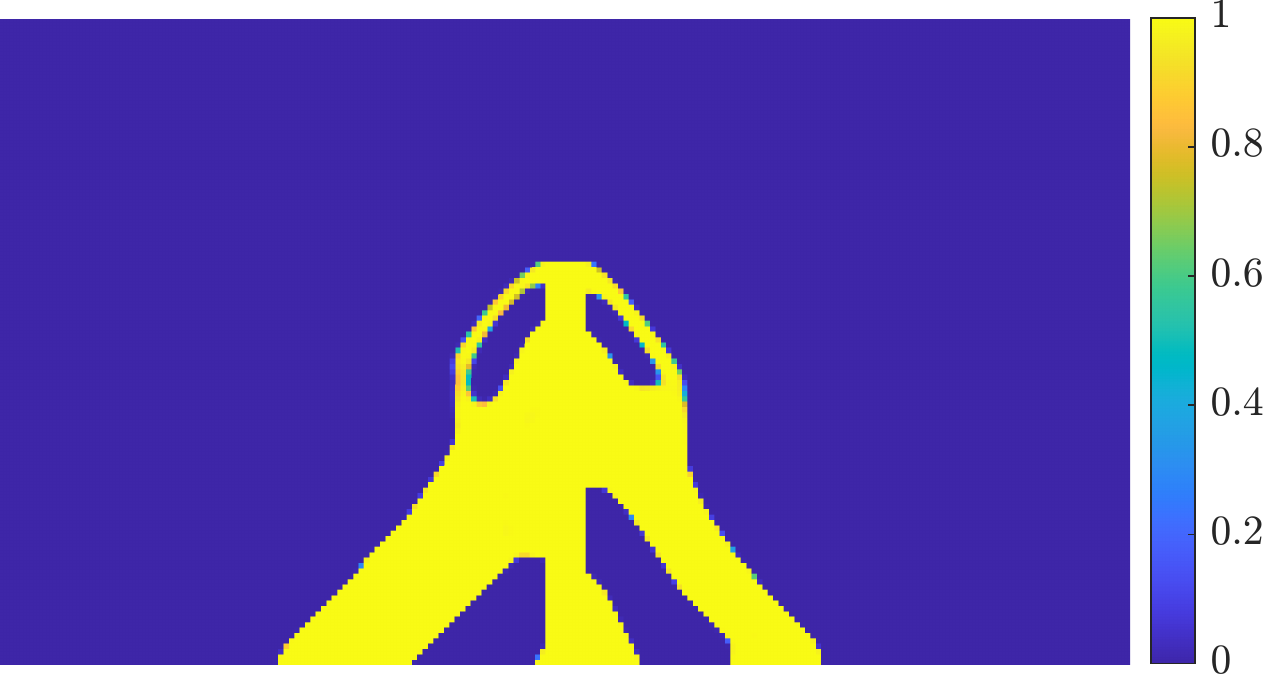}}
	\hspace{2pt}
	\subfloat[$\delta_{E|\mathit{\Upsilon}} = 3.00$, \\ \vspace{1pt} $p_\mathit{\Upsilon}$ = \{1.17, 1.33, 1.67, 2.00\}, \\ \vspace{1pt} $f = 0.011459, \ DM = 0.35 \%$.]{\includegraphics[width=\mywdth]{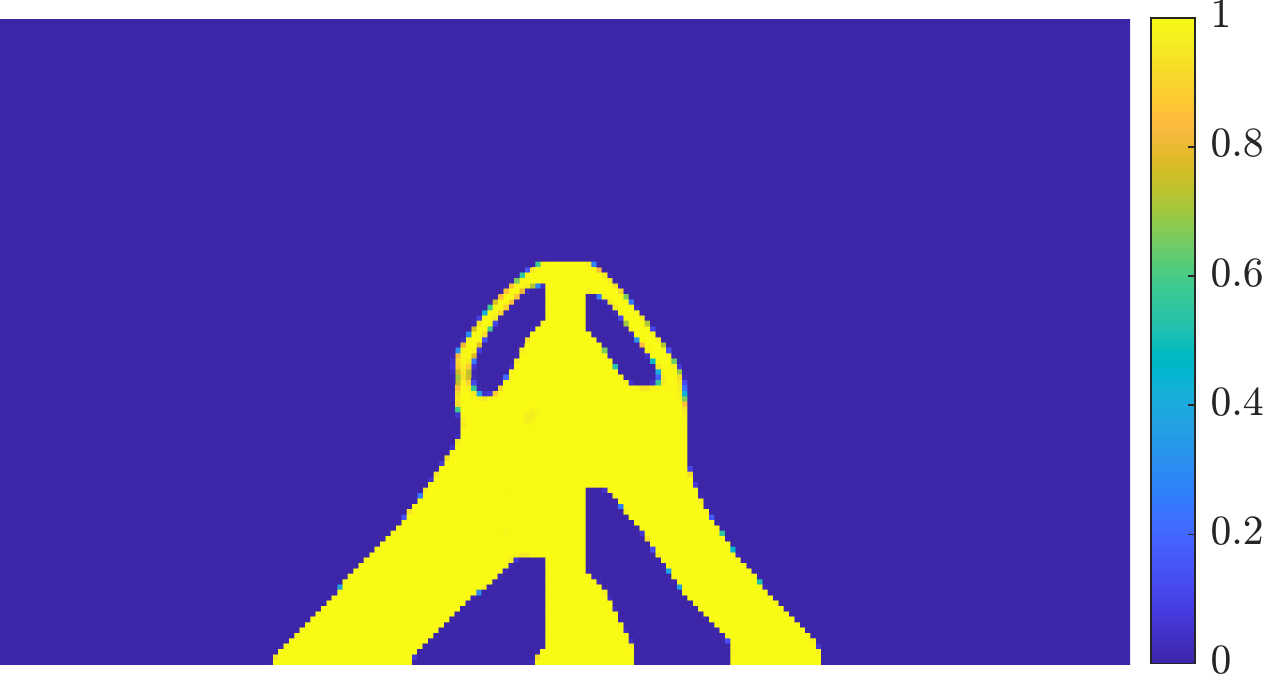}}
	\caption{Effect of $p_E$ vs $p_\mathit{\Upsilon}$ on optimized designs with  $p_\alpha = 2.5e-6$ at $Re = 1$ with pressure coupling. At iterations \{21, 41, 61, 81\}, $p_E$ is set to \{1.5, 2.0, 3.0, 4.0\} while $p_\mathit{\Upsilon}$ is calculated for different $\delta_{E|\mathit{\Upsilon}}$ as shown in the captions. $f$ is the objective function while $DM$ is the discreteness measure.}
	\label{fig:pe_vs_pu_2.5e-6_prs}
\end{figure}

\begin{figure}[t!]
	\newcommand{\mywdth}{0.325\textwidth}
	\captionsetup{font={color=black}}
	\captionsetup[subfigure]{font+={color=black}}
	\centering
	\captionsetup[subfigure]{justification=centering}
	\subfloat[$\delta_{E|\mathit{\Upsilon}} = 1.00$, \\ \vspace{1pt} $p_\mathit{\Upsilon}$ = \{1.50, 2.00, 3.00, 4.00\}, \\ \vspace{1pt} $f = 0.016401, \ DM = 0.53 \%$.]{\includegraphics[width=\mywdth]{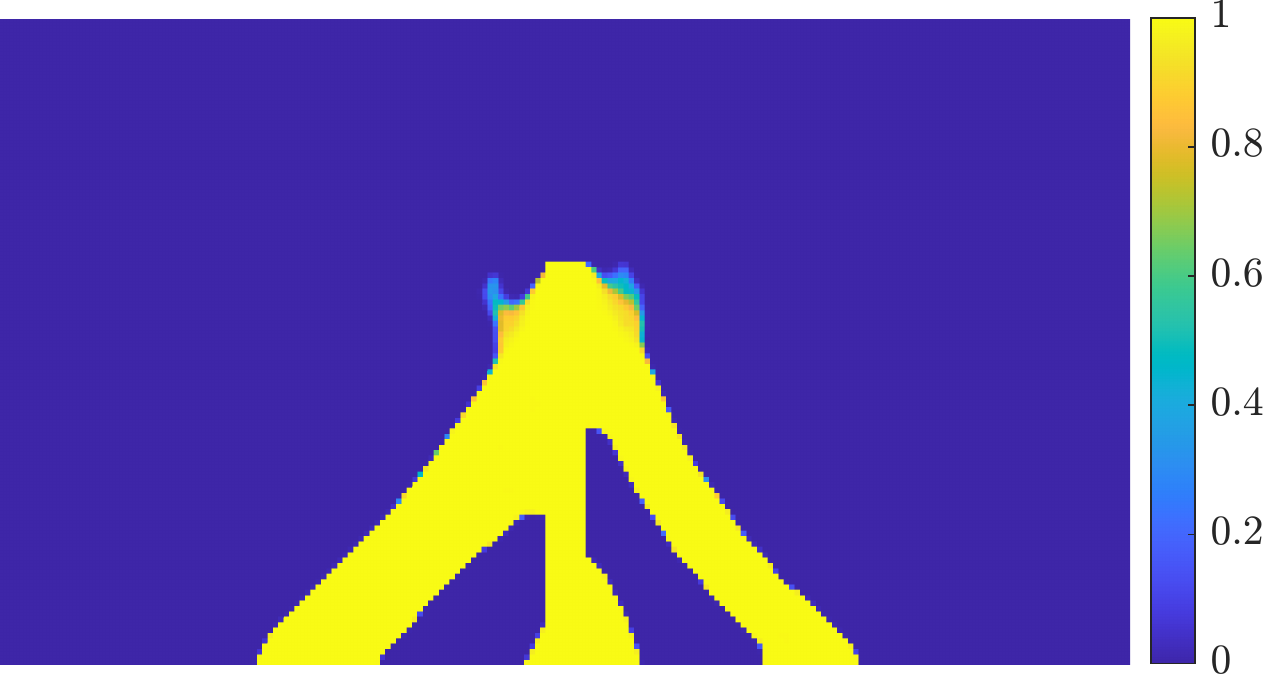}}
	\hspace{2pt}
	\subfloat[$\delta_{E|\mathit{\Upsilon}} = 1.50$, \\ \vspace{1pt} $p_\mathit{\Upsilon}$ = \{1.33, 1.67, 2.33, 3.00\}, \\ \vspace{1pt} $f = 0.016582, \ DM = 0.34 \%$.]{\includegraphics[width=\mywdth]{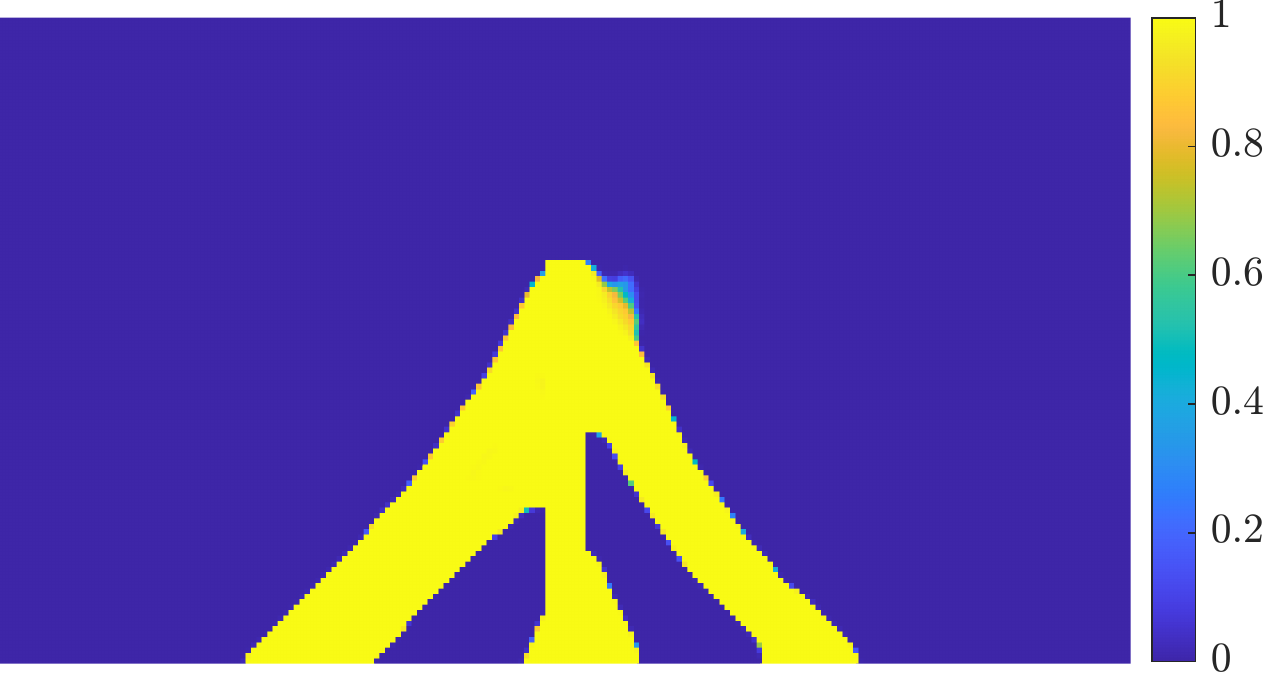}}
	\hspace{2pt}
	\subfloat[$\delta_{E|\mathit{\Upsilon}} = 2.00$, \\ \vspace{1pt} $p_\mathit{\Upsilon}$ = \{1.25, 1.50, 2.00, 2.50\}, \\ \vspace{1pt} $f = 0.016210, \ DM = 0.19 \%$. \label{ffig:7c}]{\includegraphics[width=\mywdth]{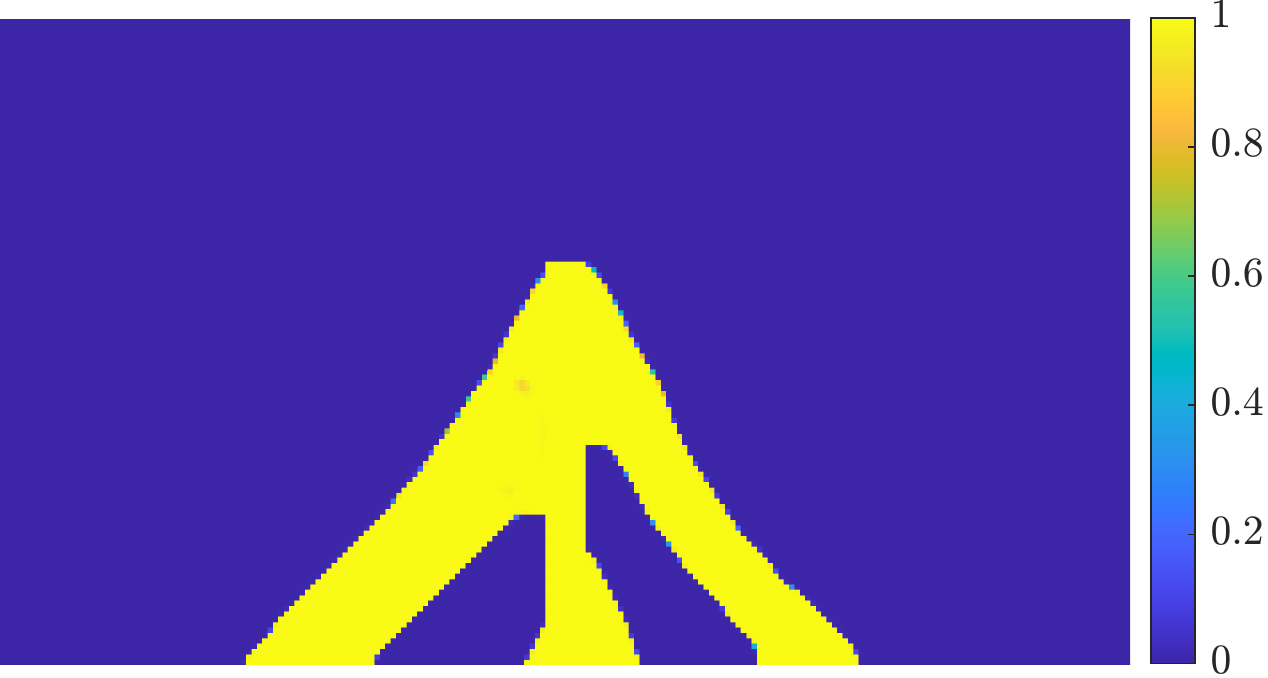}} \\
	\subfloat[$\delta_{E|\mathit{\Upsilon}} = 2.50$, \\ \vspace{1pt} $p_\mathit{\Upsilon}$ = \{1.20, 1.40, 1.80, 2.20\}, \\ \vspace{1pt} $f = 0.016413, \ DM = 0.17 \%$.]{\includegraphics[width=\mywdth]{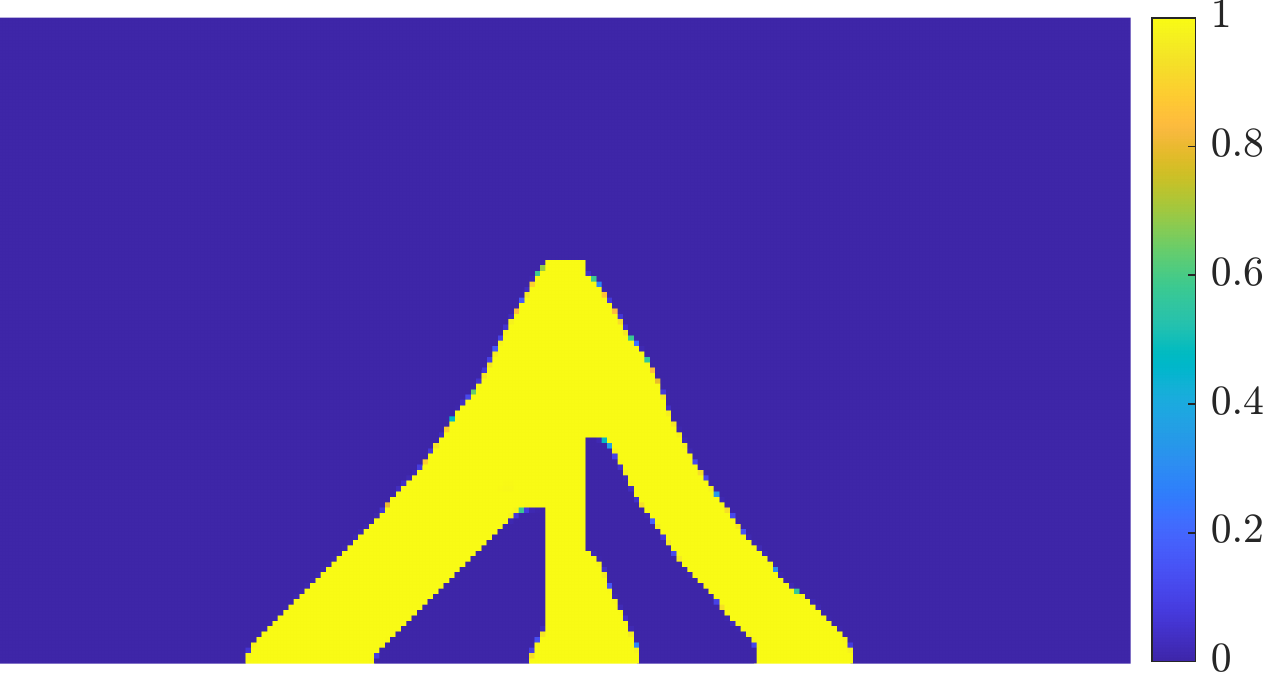}}
	\hspace{2pt}
	\subfloat[$\delta_{E|\mathit{\Upsilon}} = 3.00$, \\ \vspace{1pt} $p_\mathit{\Upsilon}$ = \{1.17, 1.33, 1.67, 2.00\}, \\ \vspace{1pt} $f = 0.016328, \ DM = 0.16 \%$.]{\includegraphics[width=\mywdth]{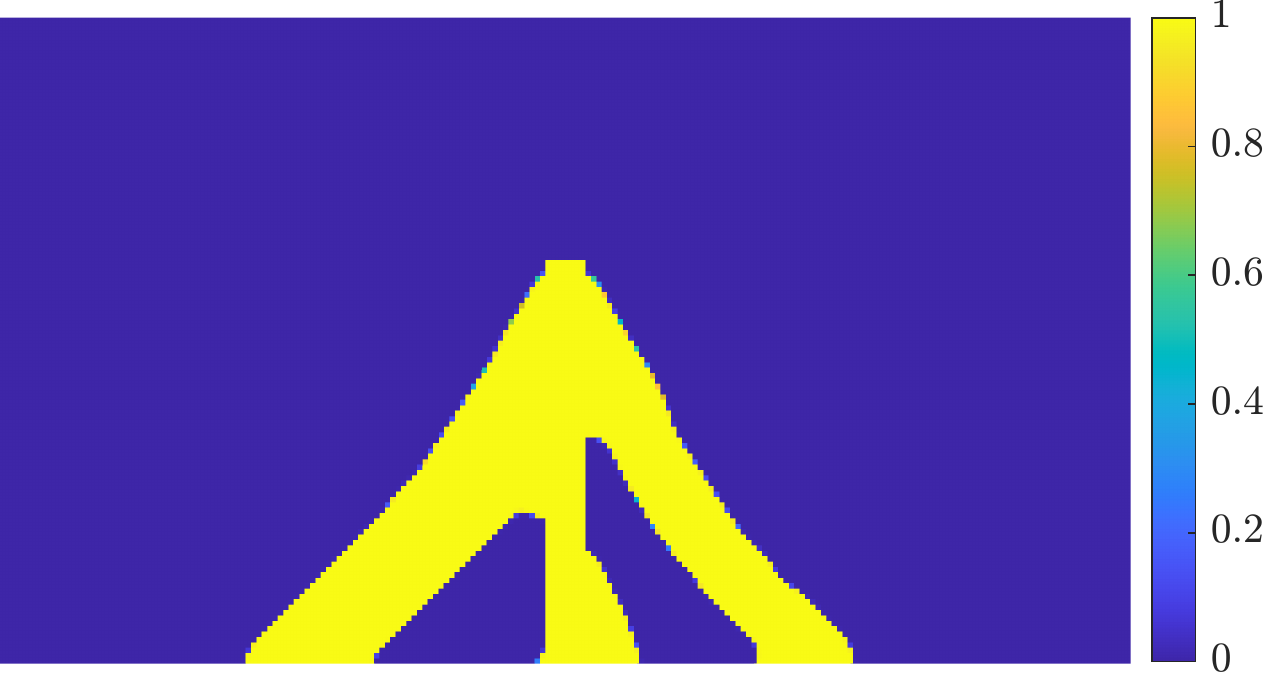}}
	\caption{Effect of $p_E$ vs $p_\mathit{\Upsilon}$ on optimized designs with  $p_\alpha = 2.5e-6$ at $Re = 1$ with total stress coupling. At iterations \{21, 41, 61, 81\}, $p_E$ is set to \{1.5, 2.0, 3.0, 4.0\} while $p_\mathit{\Upsilon}$ is calculated for different $\delta_{E|\mathit{\Upsilon}}$ as shown in the captions. $f$ is the objective function while $DM$ is the discreteness measure.}
	\label{fig:pe_vs_pu_2.5e-6_ttl}
\end{figure}

From our numerical experiments, we didn't notice any convexification effects in implementing continuation of $p_\alpha$, hence it's set to a single value throughout the entire optimization process. However, we noticed that the optimized designs are far more sensitive to the selection of $p_\alpha$ than to the selection of $p_E$ or $p_\mathit{\Upsilon}$. This is in agreement with the conclusion reached in \ \citep[p.~990]{Lundgaard2018}. To settle first the selection of $\delta_{E|\mathit{\Upsilon}}$, hence $p_\mathit{\Upsilon}$, we present the optimized designs for different $p_\mathit{\Upsilon}$ for two $p_\alpha$ values for the case of $Re = 1$. Figure \ref{fig:pe_vs_pu_5.25e-6_prs}\footnote{The discreteness measure $DM$ is calculated following Eq. 41 in \ \citep[p.~415]{Sigmund2007}.} shows the effect of different $\delta_{E|\mathit{\Upsilon}}$ on the optimized designs under pressure coupling for $p_\alpha = 5.25e-6$. It seems that, generally, close-to-discrete designs are better obtained for $\delta_{E|\mathit{\Upsilon}}$ somewhere between 2 and 3. Meaning that intermediate density elements are exposed to a relatively stronger fluid forces, hence driving their density values towards the discrete limits. In addition, the effect of $\delta_{E|\mathit{\Upsilon}}$ is far more significant when $p_\alpha$ is not properly selected. A value of $p_\alpha = 2.5e-6$ seemed more appropriate for the flow conditions of $Re = 1$ for both pressure and total stress coupling as can be seen in Figs. \ref{fig:pe_vs_pu_2.5e-6_prs} and \ref{fig:pe_vs_pu_2.5e-6_ttl}, respectively. It can be noted that the effect of $\delta_{E|\mathit{\Upsilon}}$ in the case of the more appropriate $p_\alpha = 2.5e-6$ is far less significant, yet still noticeable. Similarly to the case of $p_\alpha = 5.25e-6$, a value of $\delta_{E|\mathit{\Upsilon}}$ between 2 and 3 produces better designs for the case of $p_\alpha = 2.5e-6$ as well. It's critical to mention that the above discussion on the results presented in Figs. \ref{fig:pe_vs_pu_5.25e-6_prs}, \ref{fig:pe_vs_pu_2.5e-6_prs}, and \ref{fig:pe_vs_pu_2.5e-6_ttl} concerns the effect of the relative values of $p_\mathit{\Upsilon}$ vs $p_E$ rather than the absolute values of $p_\mathit{\Upsilon}$.

\begin{figure}[b!]
	\centering
	\captionsetup[subfigure]{justification=centering}
	\captionsetup{font+={color=black}}
	\captionsetup[subfigure]{font+={color=black}}
	\newcommand{\mywdth}{0.325\textwidth}
	\newcommand{\myvspcng}{-1pt}
	\subfloat[$p_\alpha = 1e-6, \ f = 0.00008945,$ \\ \vspace{1pt} $DM = 2.20 \%$. \label{fig:8a}]{\includegraphics[width=\mywdth]{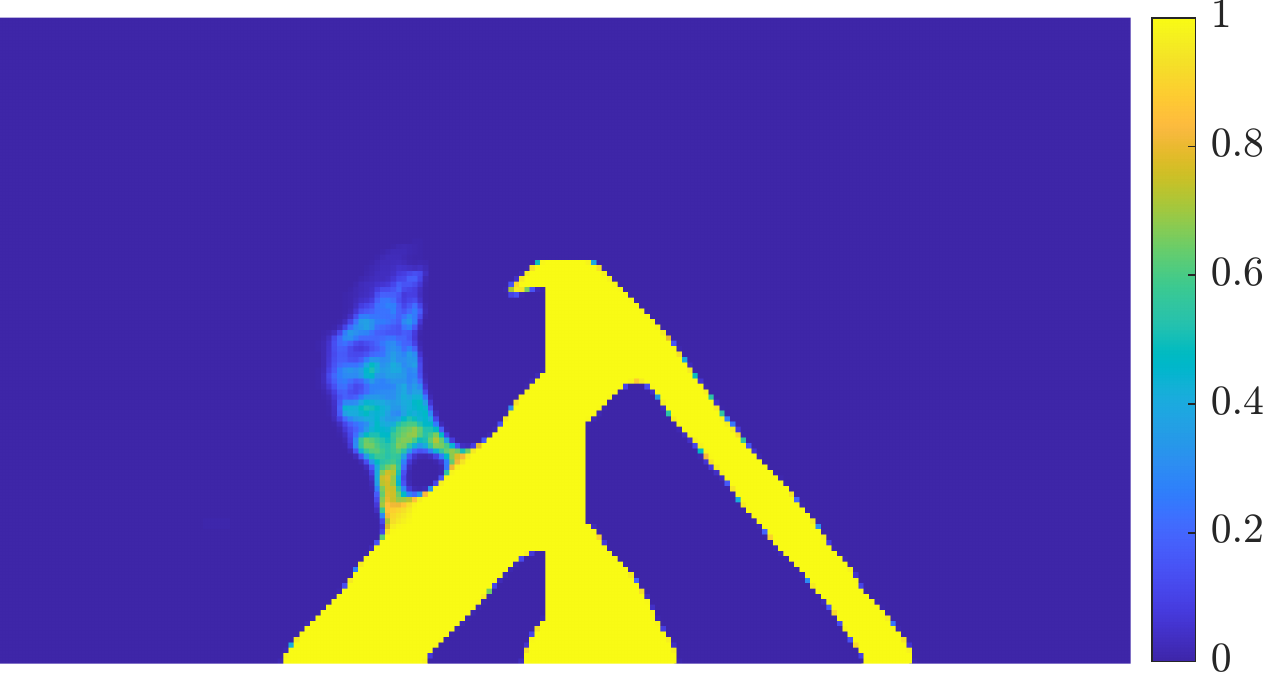}}
	\hspace{2pt}
	\subfloat[$p_\alpha = 9e-7, \ f = 0.00009682,$ \\ \vspace{1pt} $DM = 2.02 \%$.]{\includegraphics[width=\mywdth]{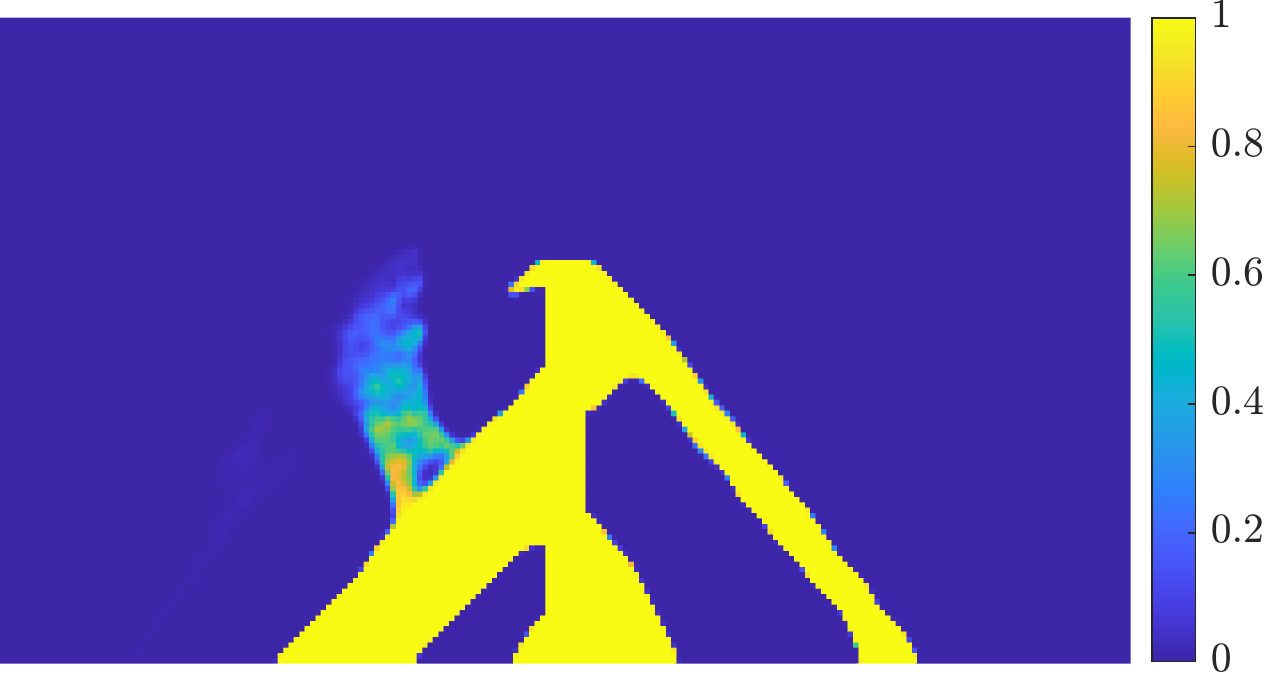}}
	\hspace{2pt}
	\subfloat[$p_\alpha = 8e-7, \ f = 0.00009529,$ \\ \vspace{1pt} $DM = 2.12
	\%$.]{\includegraphics[width=\mywdth]{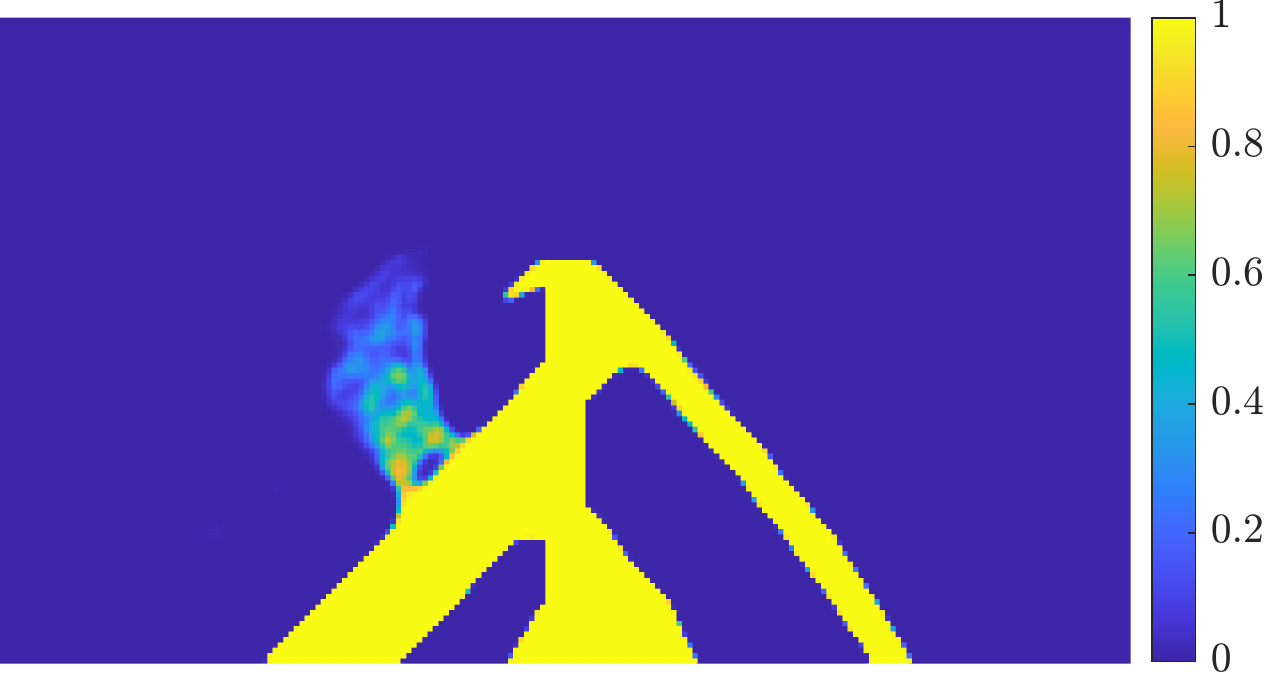}}\\[\myvspcng]
	\subfloat[$p_\alpha = 7e-7, \ f = 0.00011389,$ \\ \vspace{1pt} $DM = 1.09 \%$.]{\includegraphics[width=\mywdth]{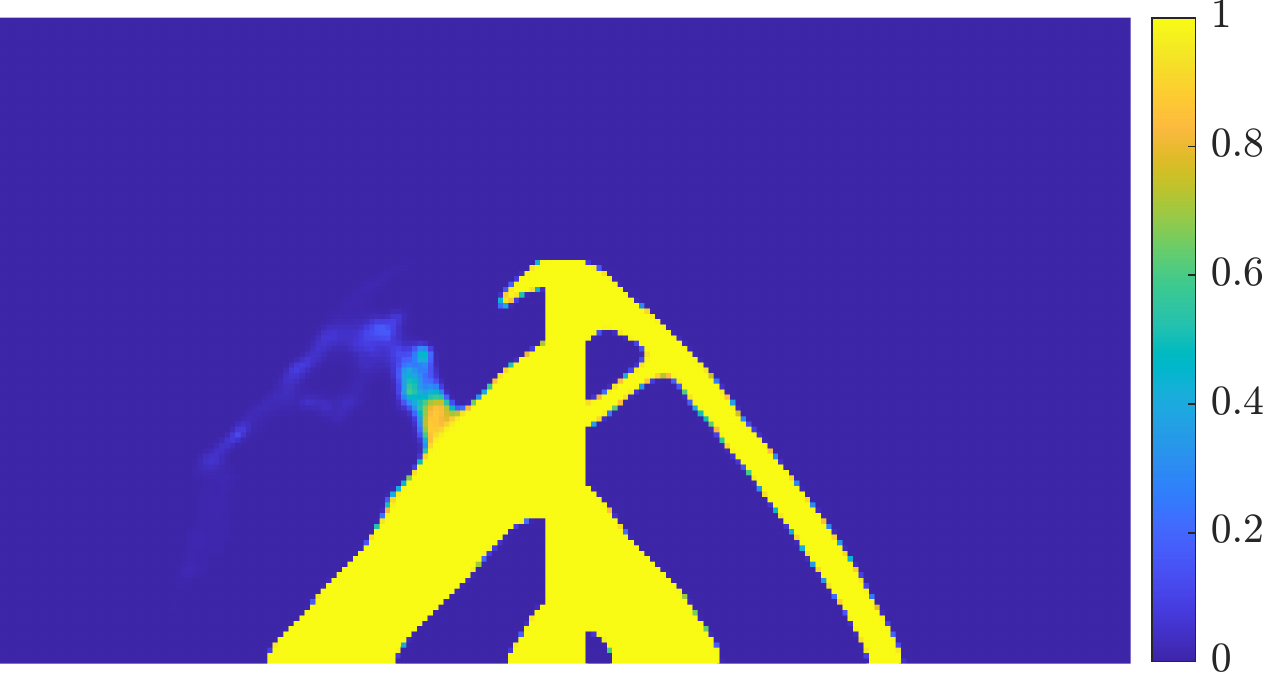}}
	\hspace{2pt}
	\subfloat[$p_\alpha = 6e-7, \ f = 0.00011513,$ \\ \vspace{1pt} $DM = 1.82 \%$. \label{fig:8e}]{\includegraphics[width=\mywdth]{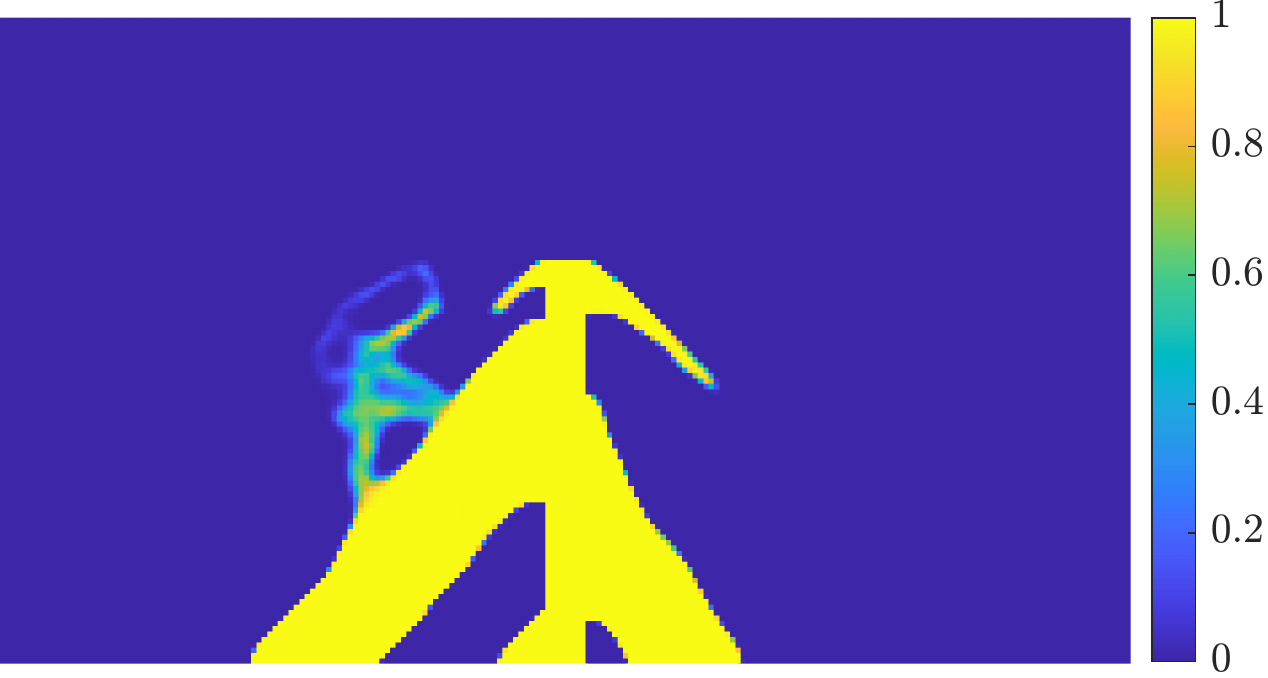}}
	\hspace{2pt}
	\subfloat[$p_\alpha = 5e-7, \ f = 0.00011391,$ \\ \vspace{1pt} $DM = 0.77	
	\%$.]{\includegraphics[width=\mywdth]{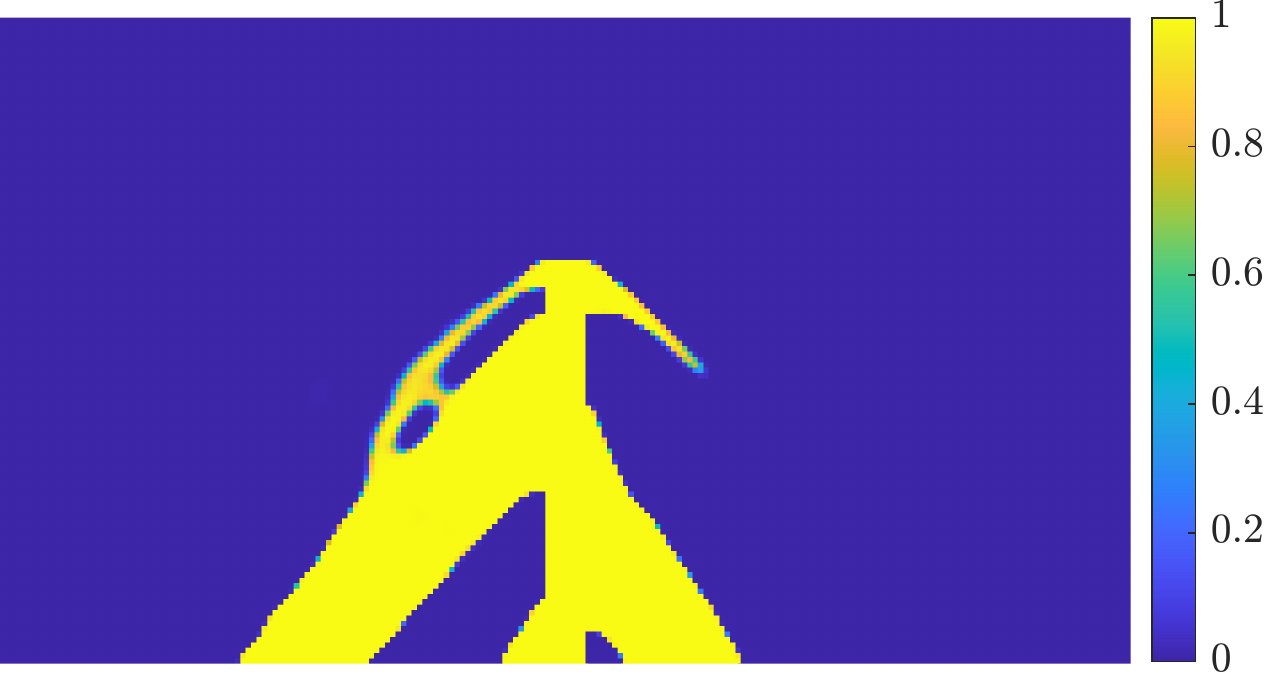}}\\[\myvspcng]
	\subfloat[$p_\alpha = 4e-7, \ f = 0.00010853,$ \\ \vspace{1pt} $DM = 0.42	
	\%$.]{\includegraphics[width=\mywdth]{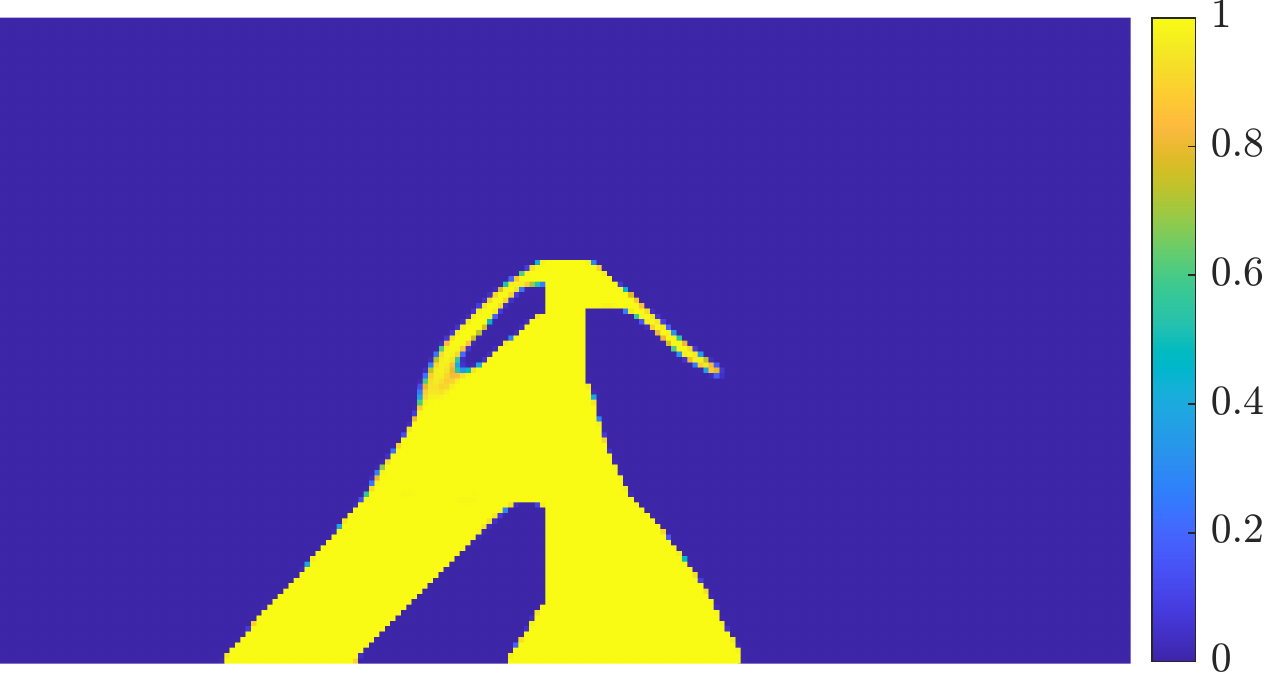}}
	\hspace{2pt}
	\subfloat[$p_\alpha = 3e-7, \ f = 0.00010678,$ \\ \vspace{1pt} $DM = 0.47
	\%$.]{\includegraphics[width=\mywdth]{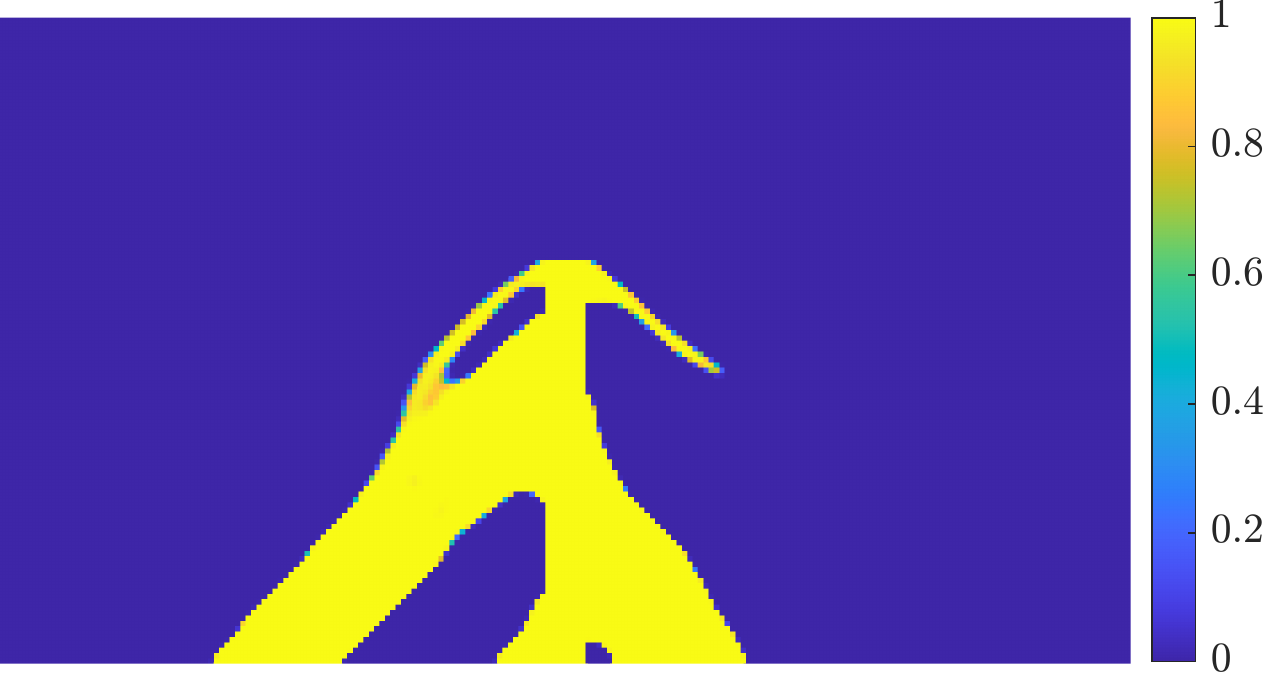}}
	\hspace{2pt}
	\subfloat[$p_\alpha = 2e-7, \ f = 0.00010525,$ \\ \vspace{1pt} $DM = 0.42
	\%$.]{\includegraphics[width=\mywdth]{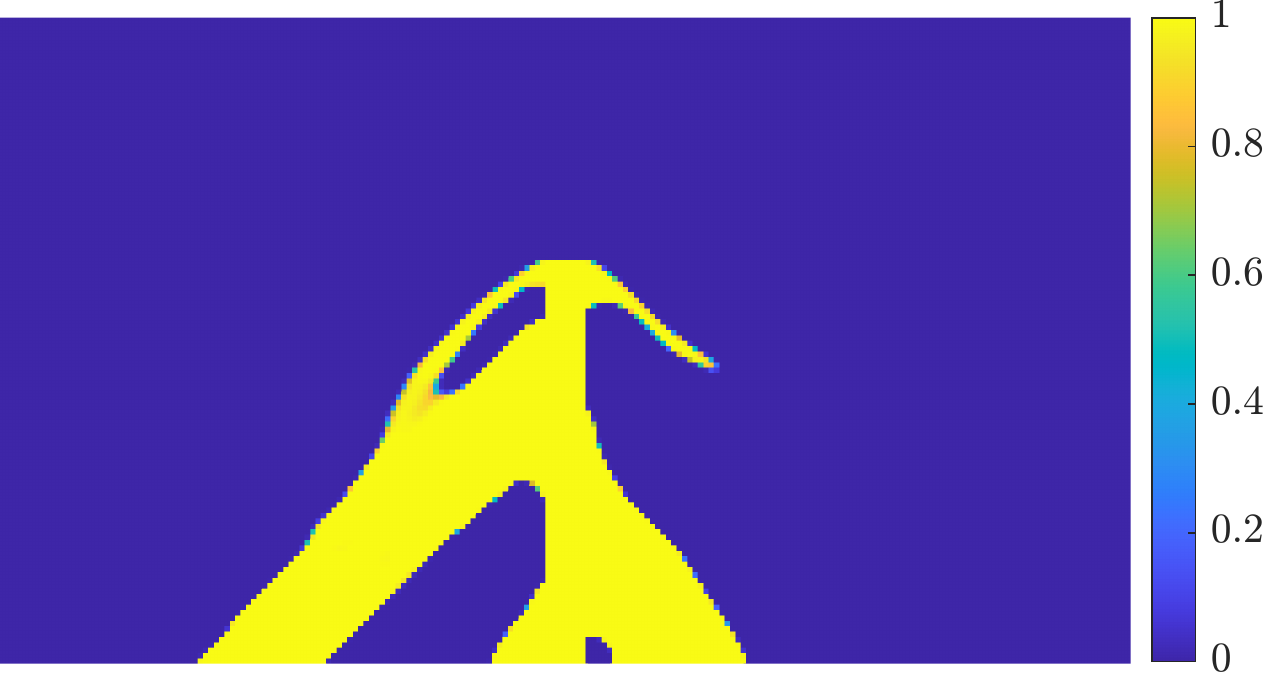}}\\[\myvspcng]
	\subfloat[$p_\alpha = 1e-7, \ f = 0.00010543,$ \\ \vspace{1pt} $DM = 0.41
	\%$.]{\includegraphics[width=\mywdth]{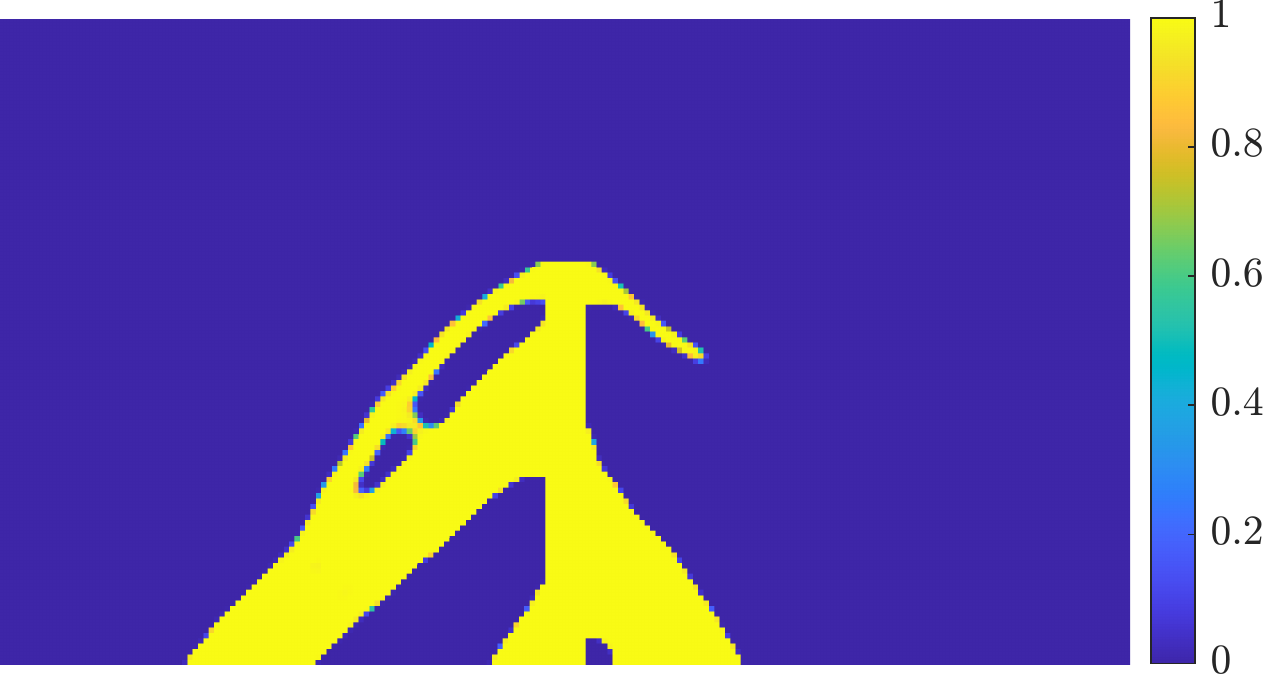}}
	\hspace{2pt}
	\subfloat[$p_\alpha = 9e-8, \ f = 0.00010585,$ \\ \vspace{1pt} $DM = 0.37
	\%$. \label{ffig:8k}]{\includegraphics[width=\mywdth]{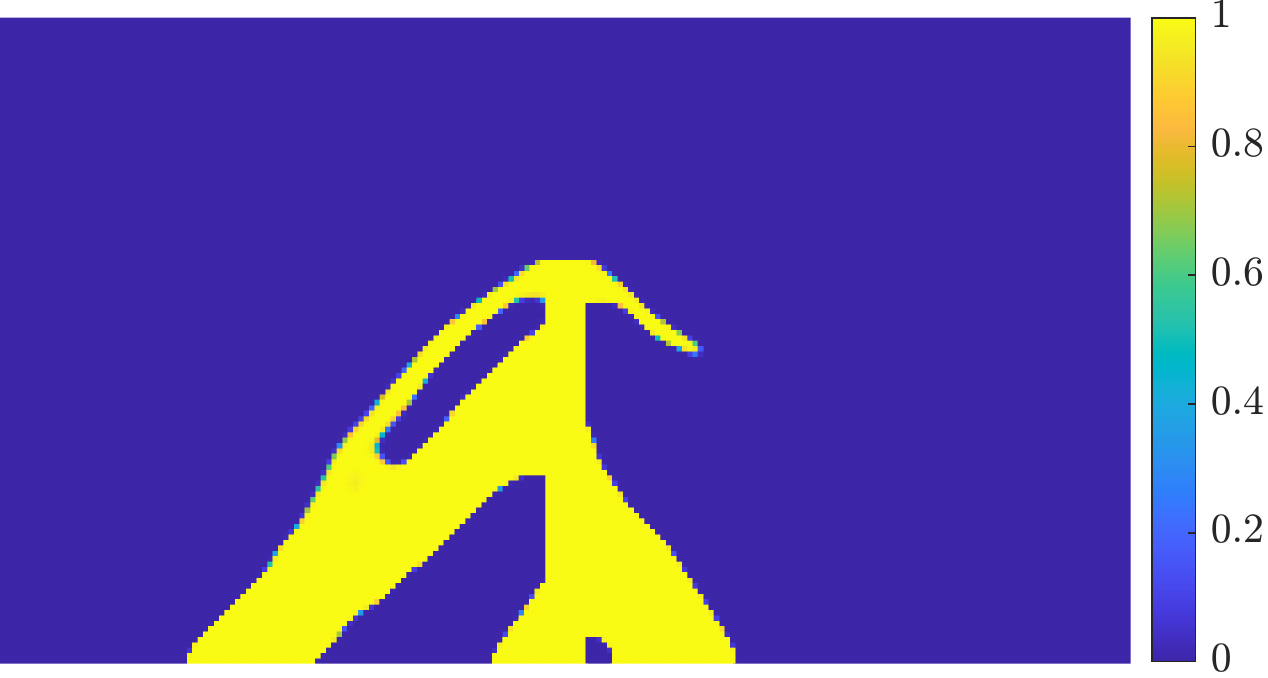}}
	\hspace{2pt}
	\subfloat[$p_\alpha = 8e-8, \ f = 0.00010346,$ \\ \vspace{1pt} $DM = 0.32
	\%$. \label{fig:8l}]{\includegraphics[width=\mywdth]{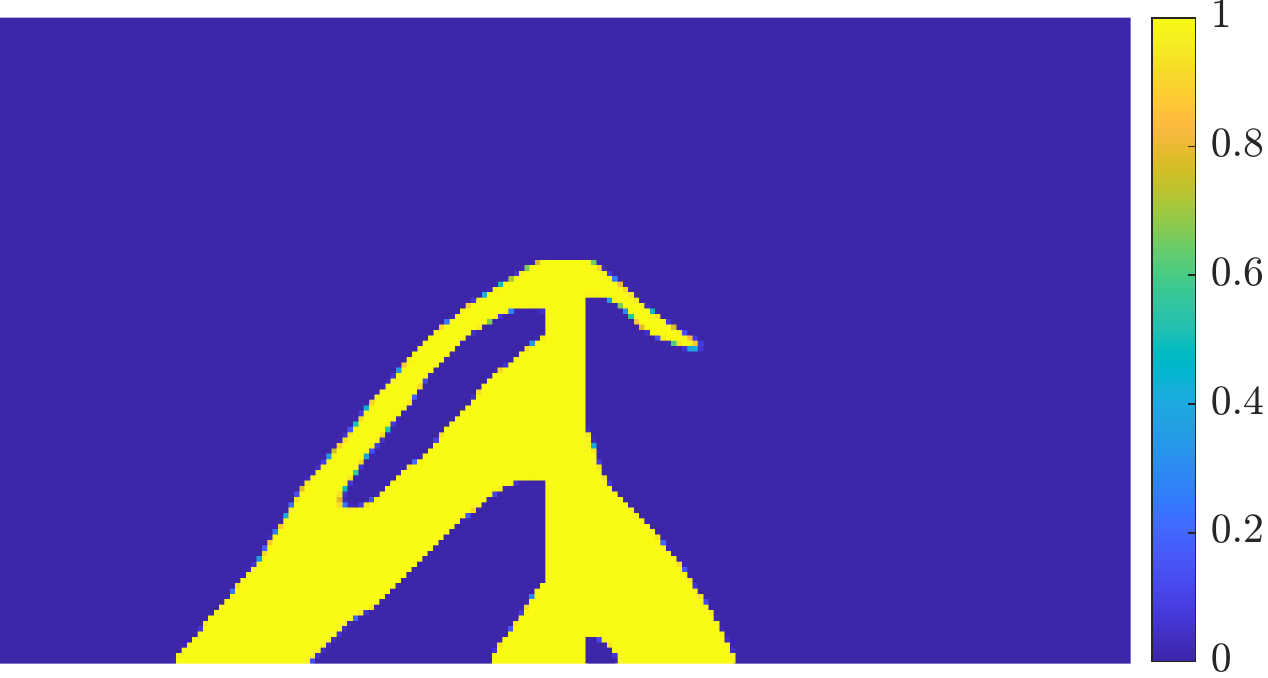}}\\[\myvspcng]
	\subfloat[$p_\alpha = 7e-8, \ f = 0.00010245,$ \\ \vspace{1pt} $DM = 0.32
	\%$.]{\includegraphics[width=\mywdth]{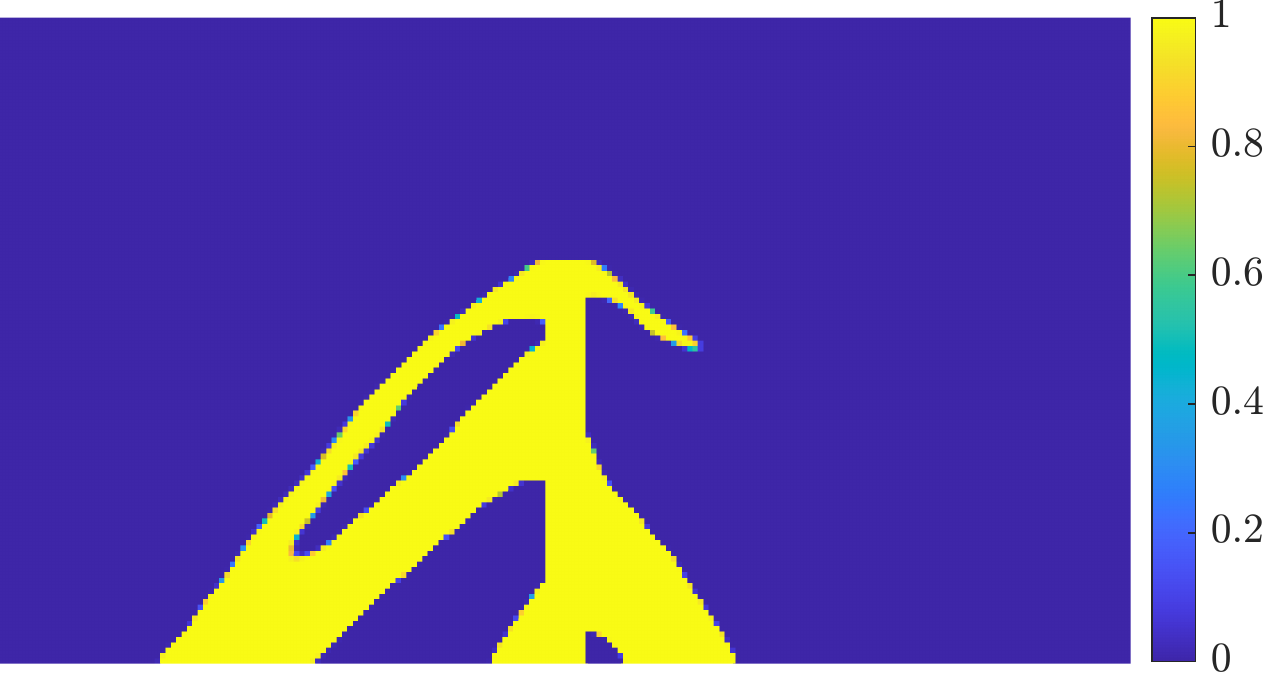}}
	\hspace{2pt}
	\subfloat[$p_\alpha = 6e-8, \ f = 0.00010107,$ \\ \vspace{1pt} $DM = 0.37
	\%$.]{\includegraphics[width=\mywdth]{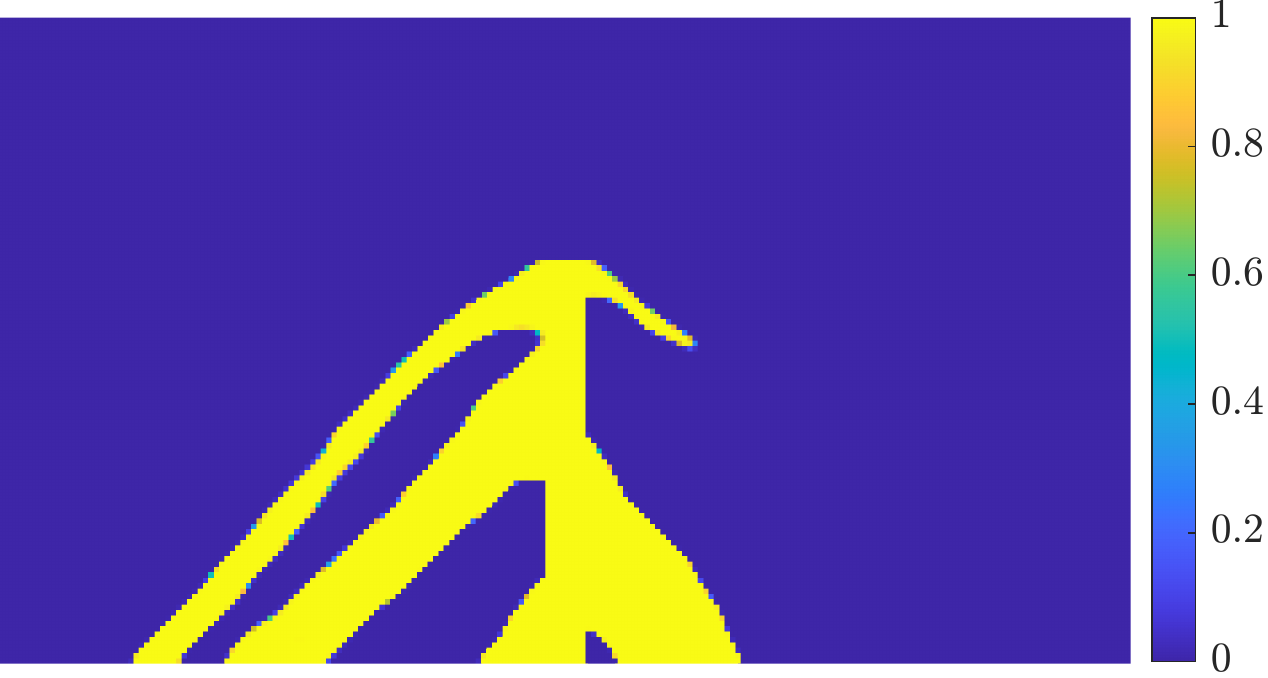}}
	\hspace{2pt}
	\subfloat[$p_\alpha = 5e-8, \ f = 0.00010066,$ \\ \vspace{1pt} $DM = 0.34
	\%$. \label{fig:8o}]{\includegraphics[width=\mywdth]{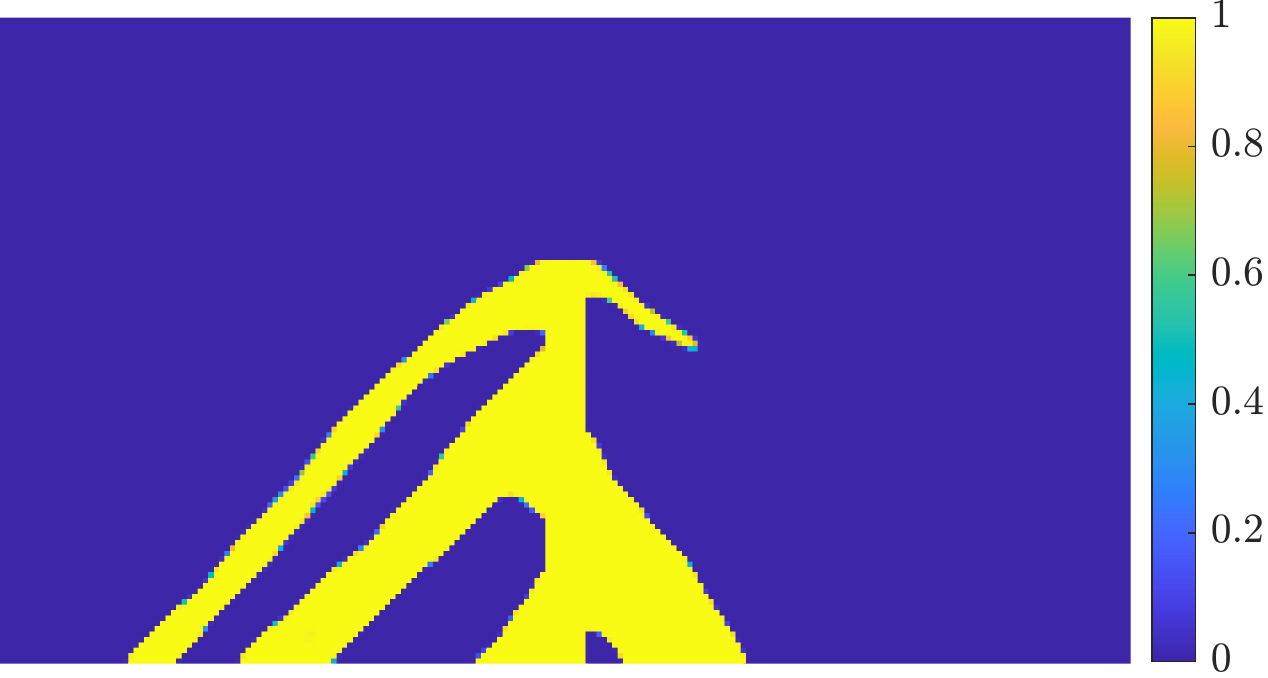}}
	\caption{Effect of $p_\alpha$ on optimized designs for the case $Re = 10$ with pressure coupling. At iterations \{21, 41, 61, 81\}, $p_E$ is set to \{1.5, 2.0, 3.0, 4.0\} while $p_\mathit{\Upsilon}$ is set to \{1.25, 1.50, 2.00, 2.50\}, i.e. $\delta_{E|\mathit{\Upsilon}} = 2$.}
	\label{fig:effect_of_p_alpha_re10_prs}
\end{figure}

\begin{figure}[b!]
	\centering
	\captionsetup{font+={color=black}}
	\captionsetup[subfigure]{font+={color=black}}
	\captionsetup[subfigure]{justification=centering}
	\newcommand{\mywdth}{0.325\textwidth}
	\newcommand{\myvspcng}{-1pt}
	\subfloat[$p_\alpha = 1e-6, \ f = 0.00013791,$ \\ \vspace{1pt} $DM = 2.68 \%$. \label{fig:9a}]{\includegraphics[width=\mywdth]{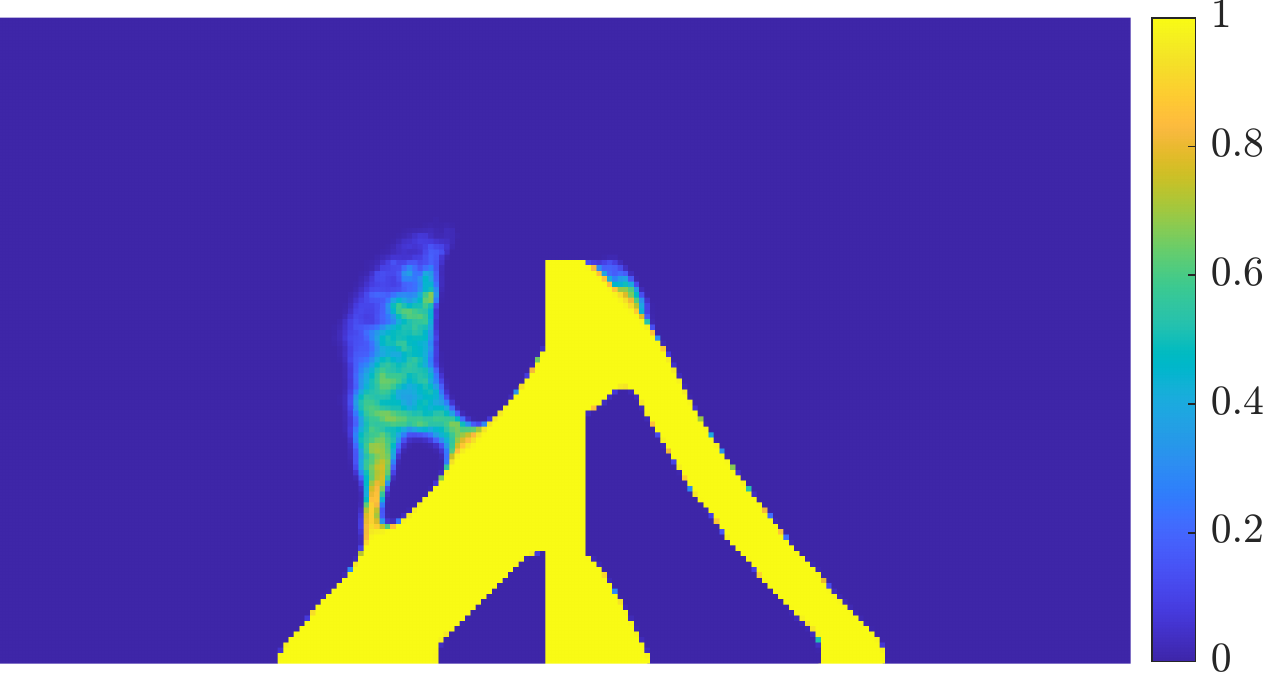}}
	\hspace{2pt}
	\subfloat[$p_\alpha = 9e-7, \ f = 0.00041615,$ \\ \vspace{1pt} $DM = 1.75 \%$.]{\includegraphics[width=\mywdth]{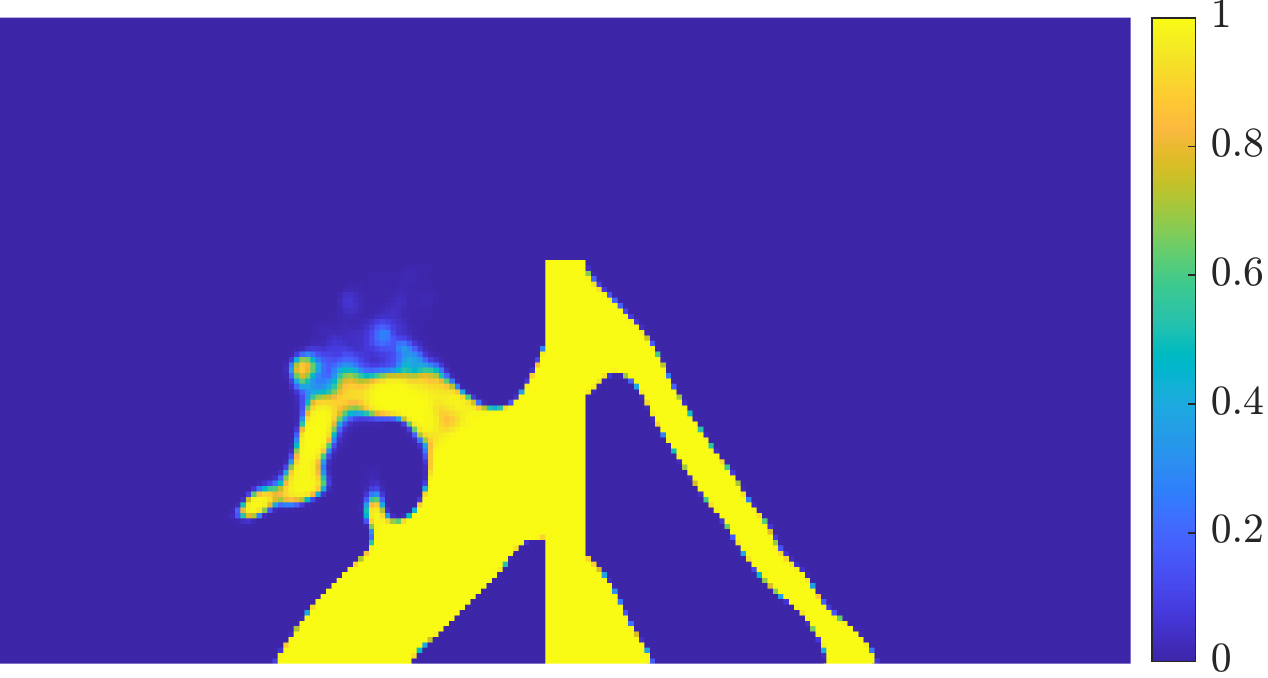}}
	\hspace{2pt}
	\subfloat[$p_\alpha = 8e-7, \ f = 0.00013941,$ \\ \vspace{1pt} $DM = 1.90 \%$.]{\includegraphics[width=\mywdth]{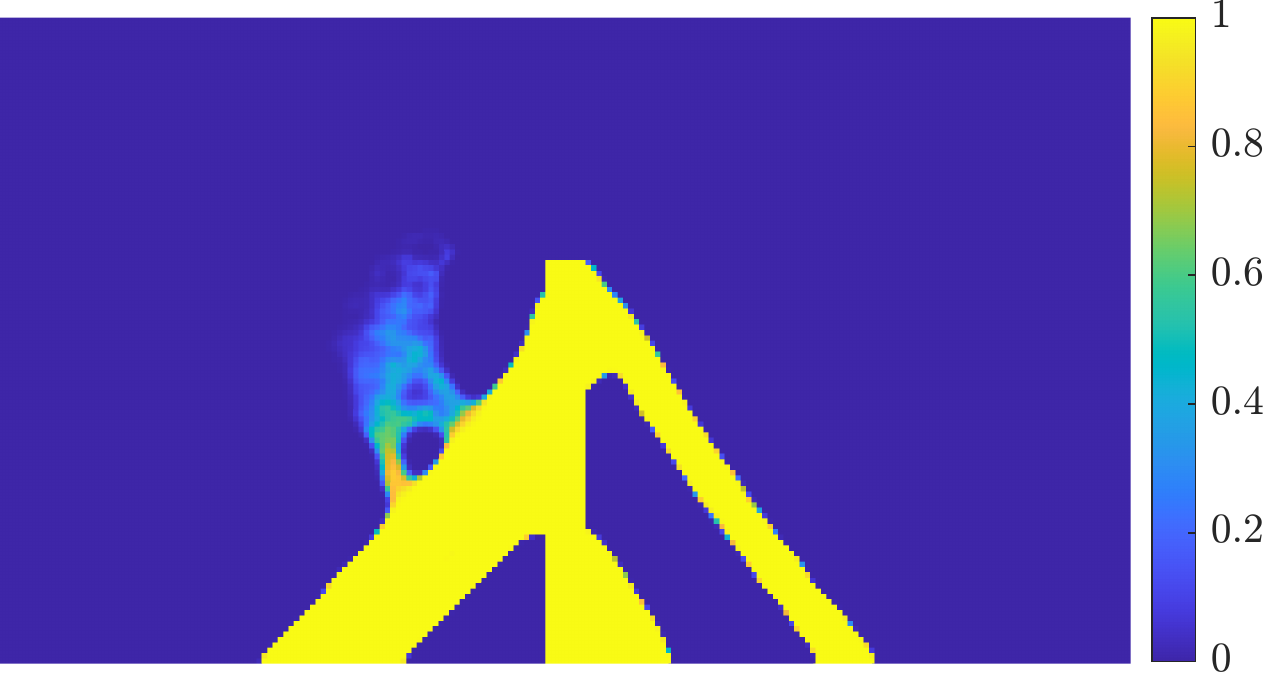}}\\[\myvspcng]
	\subfloat[$p_\alpha = 7e-7, \ f = 0.00013797,$ \\ \vspace{1pt} $DM = 1.92 \%$.]{\includegraphics[width=\mywdth]{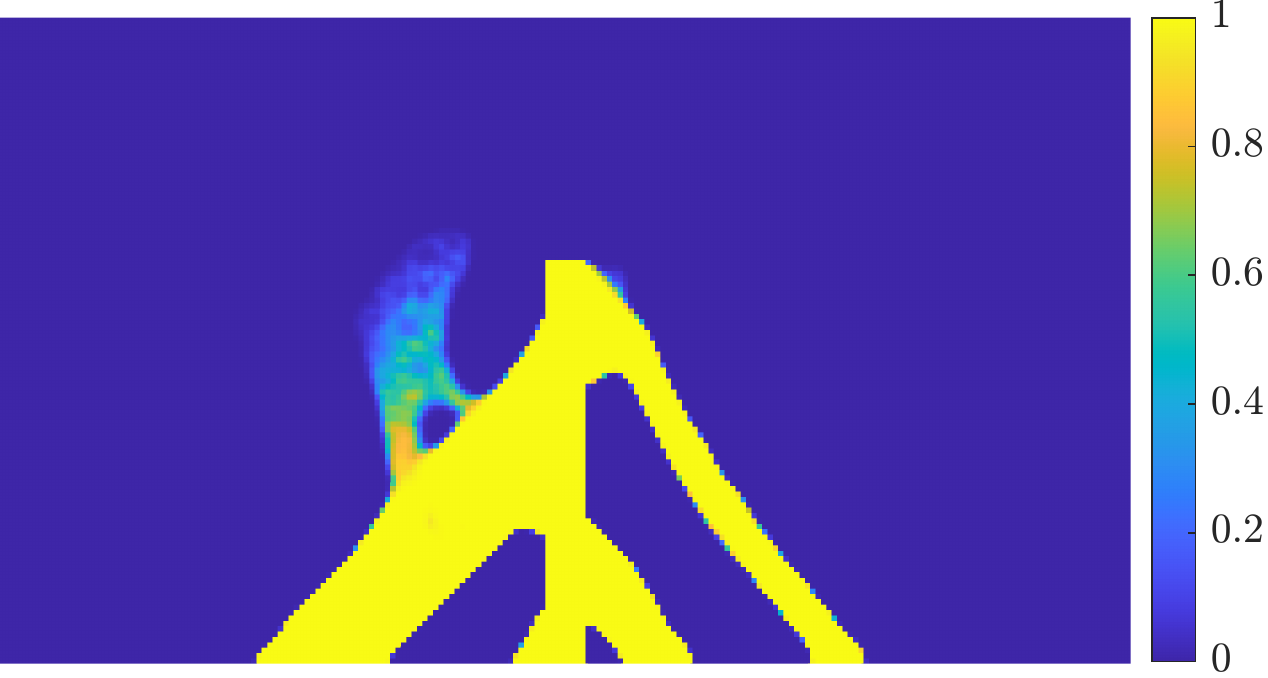}}
	\hspace{2pt}
	\subfloat[$p_\alpha = 6e-7, \ f = 0.00014292,$ \\ \vspace{1pt} $DM = 1.63 \%$.]{\includegraphics[width=\mywdth]{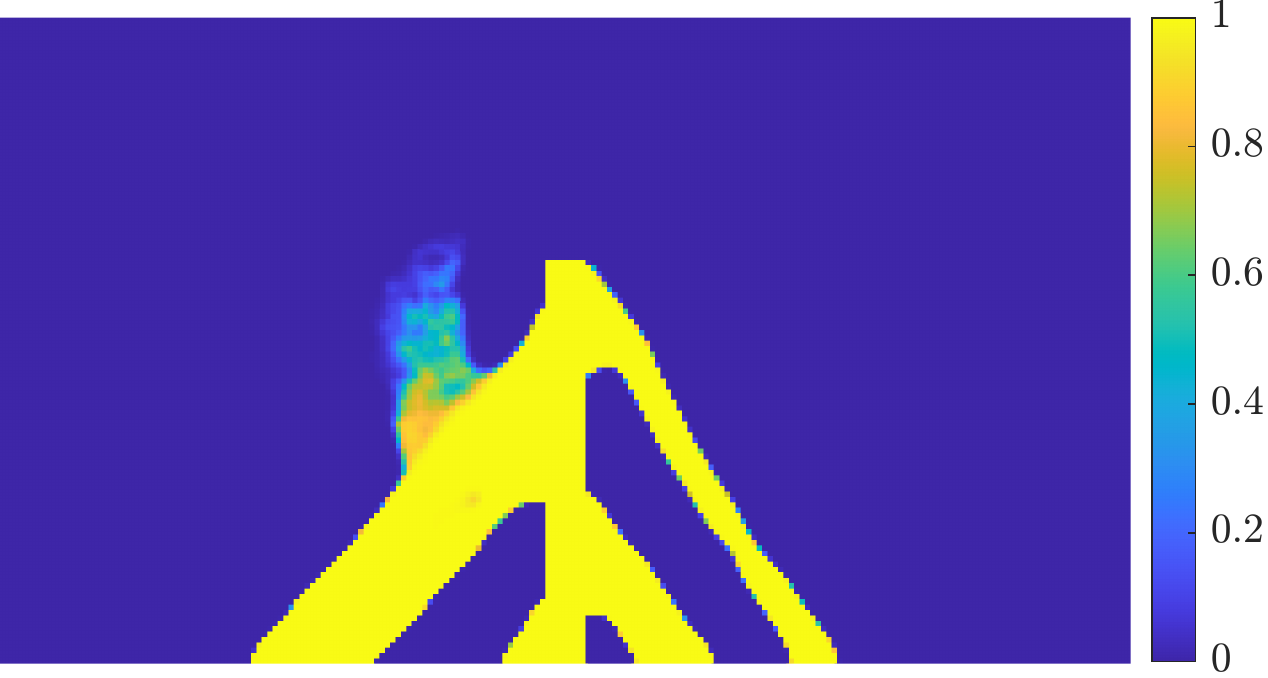}}
	\hspace{2pt}
	\subfloat[$p_\alpha = 5e-7, \ f = 0.00015141,$ \\ \vspace{1pt} $DM = 1.35 \%$.]{\includegraphics[width=\mywdth]{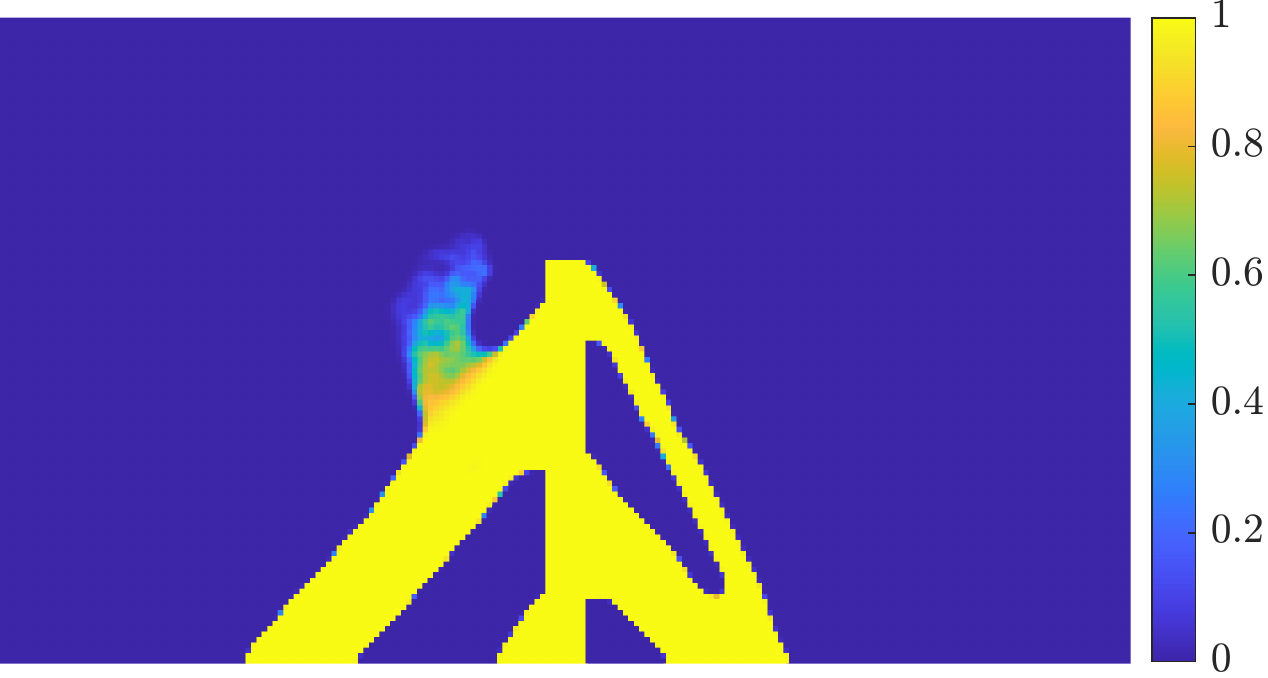}}\\[\myvspcng]
	\subfloat[$p_\alpha = 4e-7, \ f = 0.00016479,$ \\ \vspace{1pt} $DM = 0.93 \%$.]{\includegraphics[width=\mywdth]{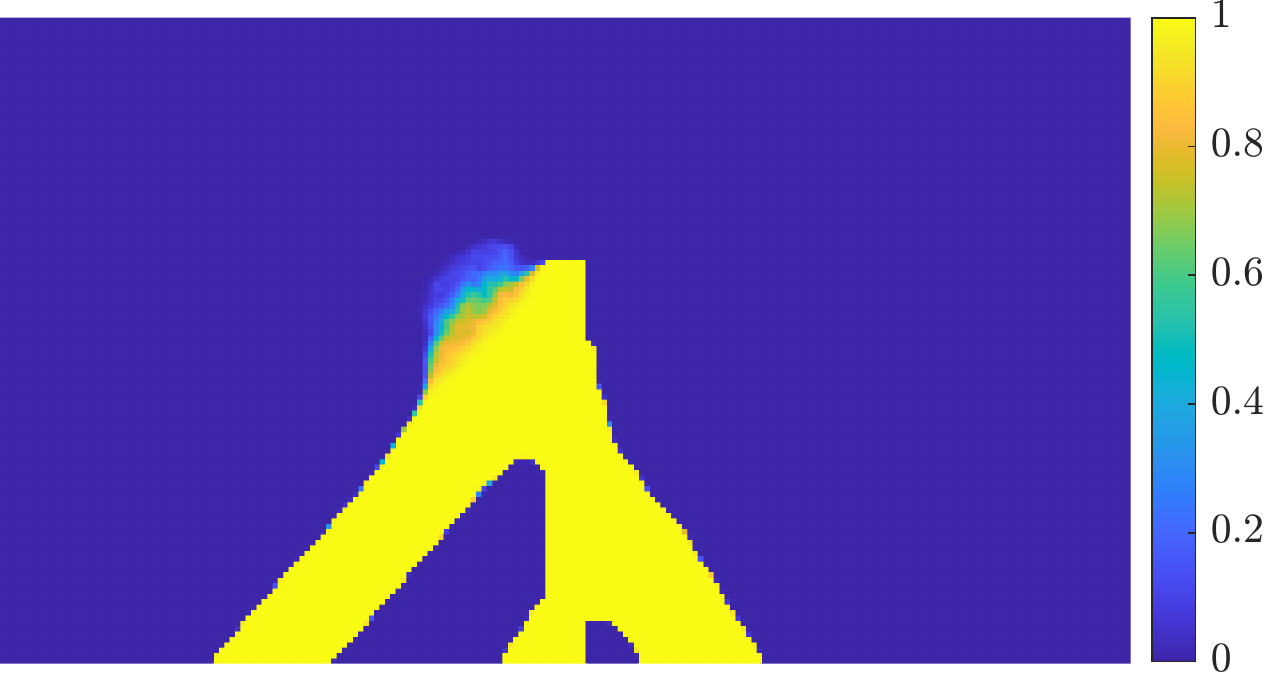}}
	\hspace{2pt}
	\subfloat[$p_\alpha = 3e-7, \ f = 0.00016824,$ \\ \vspace{1pt} $DM = 0.93 \%$.]{\includegraphics[width=\mywdth]{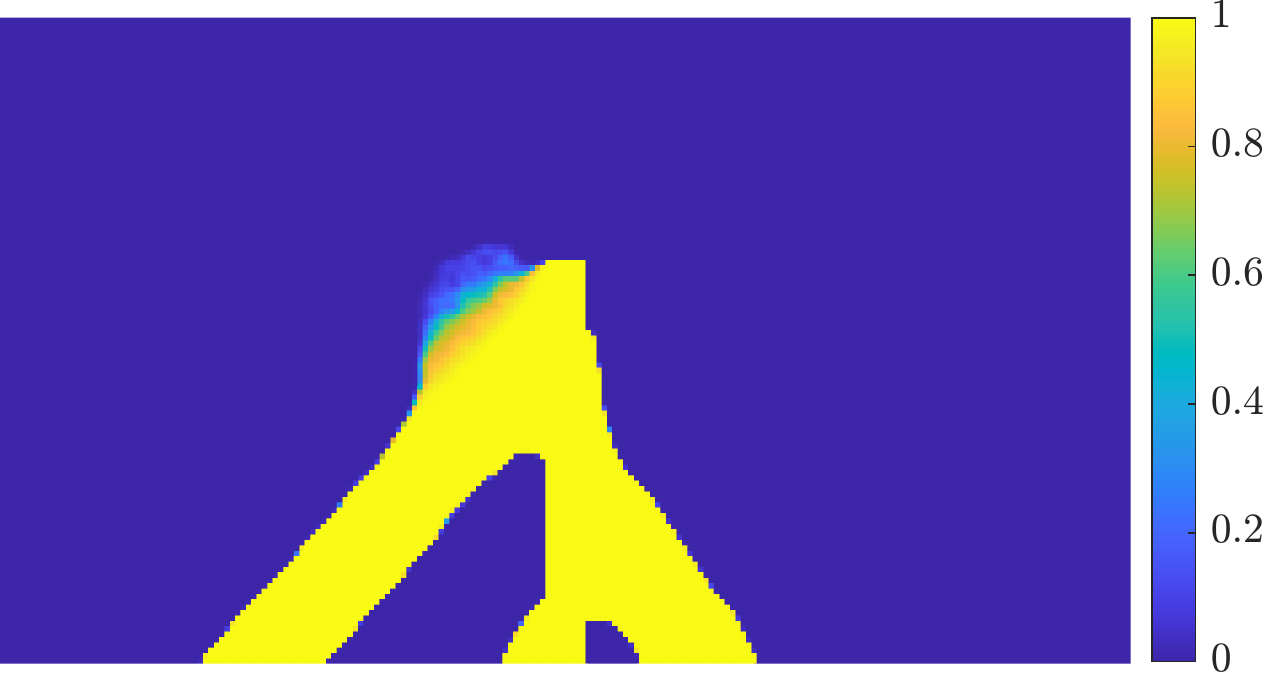}}
	\hspace{2pt}
	\subfloat[$p_\alpha = 2e-7, \ f = 0.00017229,$ \\ \vspace{1pt} $DM = 0.80 \%$. \label{fig:9i}]{\includegraphics[width=\mywdth]{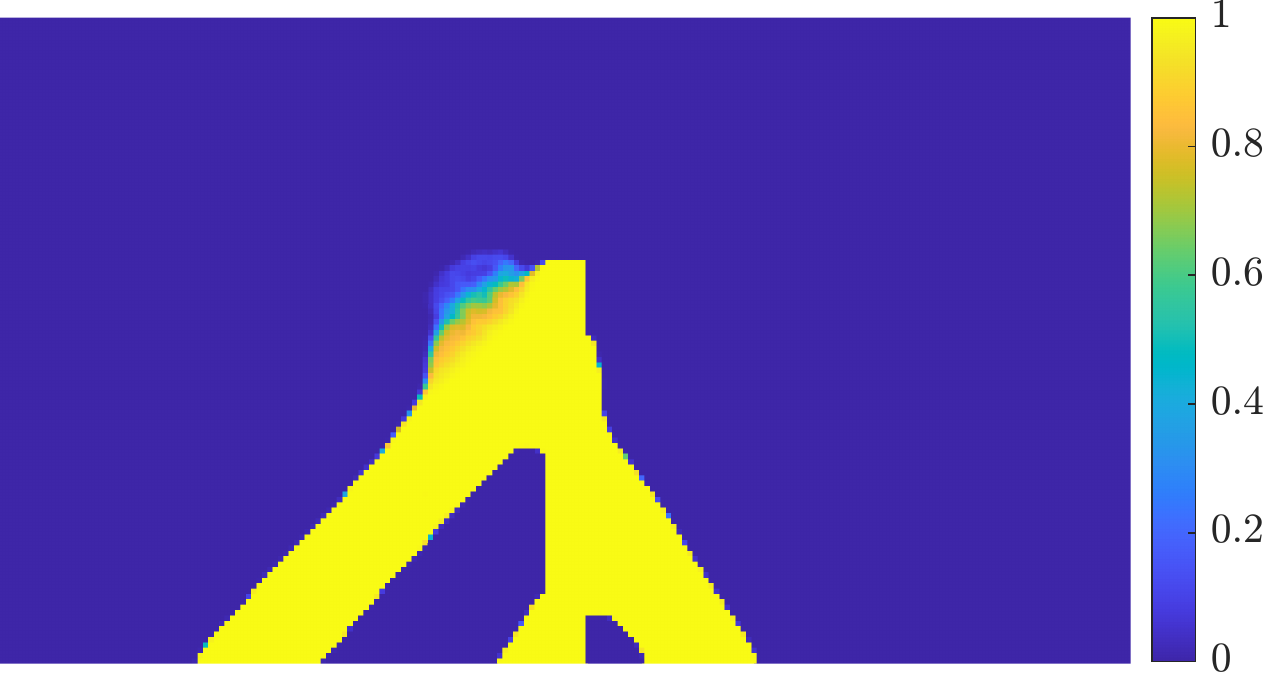}}\\[\myvspcng]
	\subfloat[$p_\alpha = 1e-7, \ f = 0.00017840,$ \\ \vspace{1pt} $DM = 0.55 \%$.]{\includegraphics[width=\mywdth]{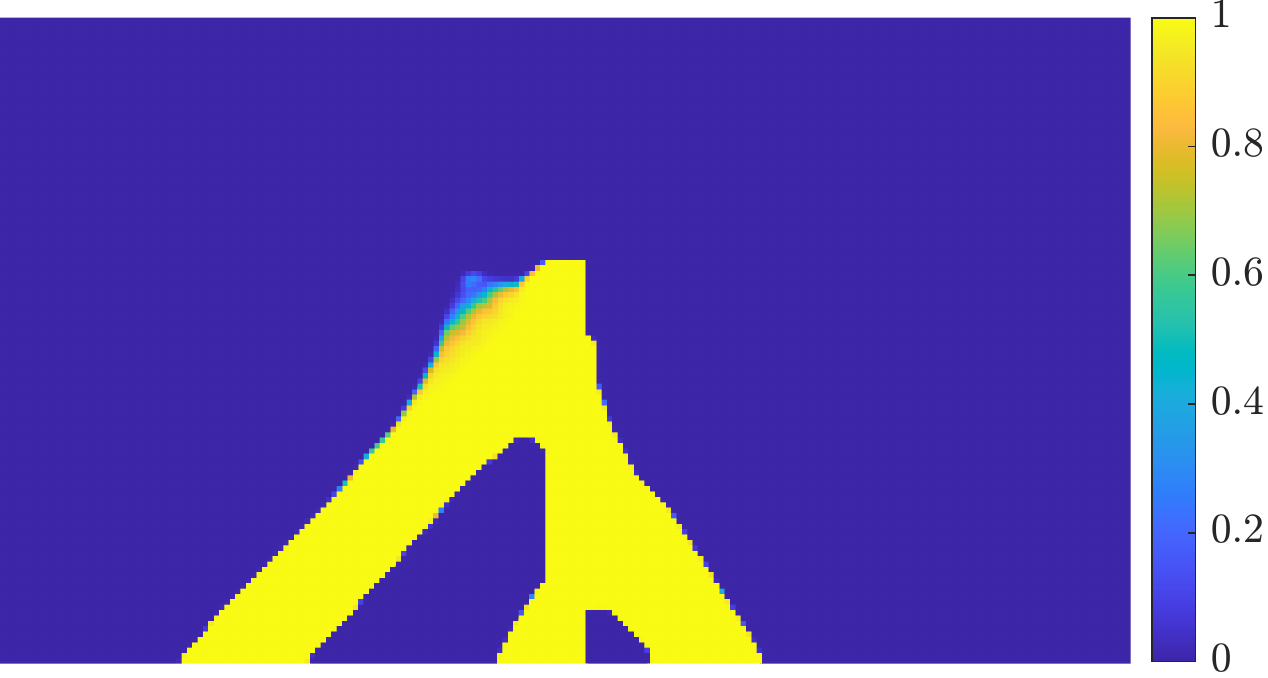}}
	\hspace{2pt}
	\subfloat[$p_\alpha = 9e-8, \ f = 0.00017952,$ \\ \vspace{1pt} $DM = 0.44 \%$.]{\includegraphics[width=\mywdth]{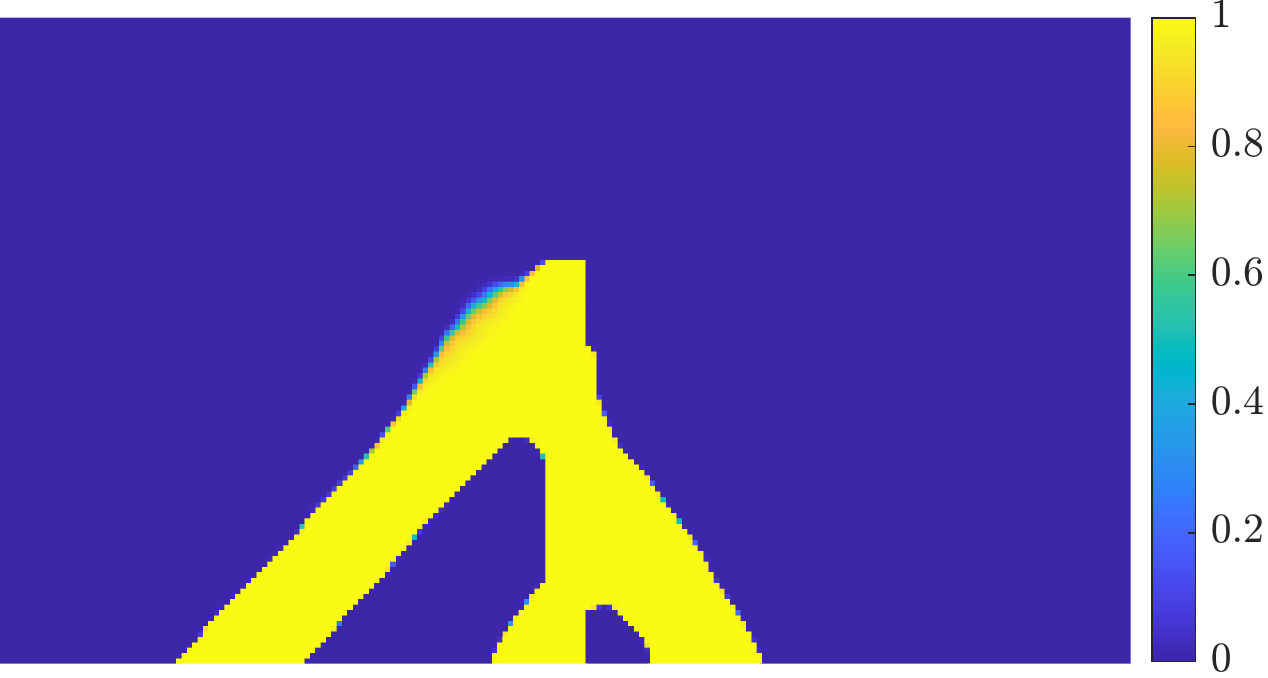}}
	\hspace{2pt}
	\subfloat[$p_\alpha = 8e-8, \ f = 0.00017985,$ \\ \vspace{1pt} $DM = 0.45 \%$.]{\includegraphics[width=\mywdth]{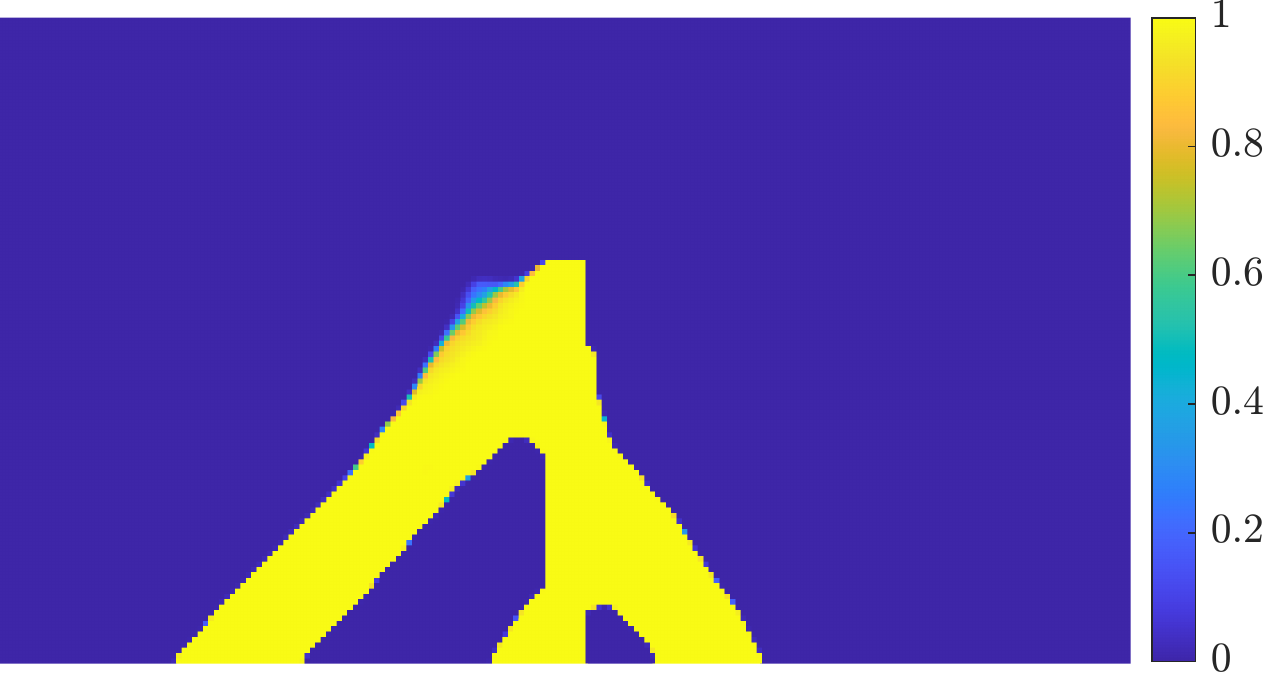}}\\[\myvspcng]
	\subfloat[$p_\alpha = 7e-8, \ f = 0.00017844,$ \\ \vspace{1pt} $DM = 0.40 \%$.]{\includegraphics[width=\mywdth]{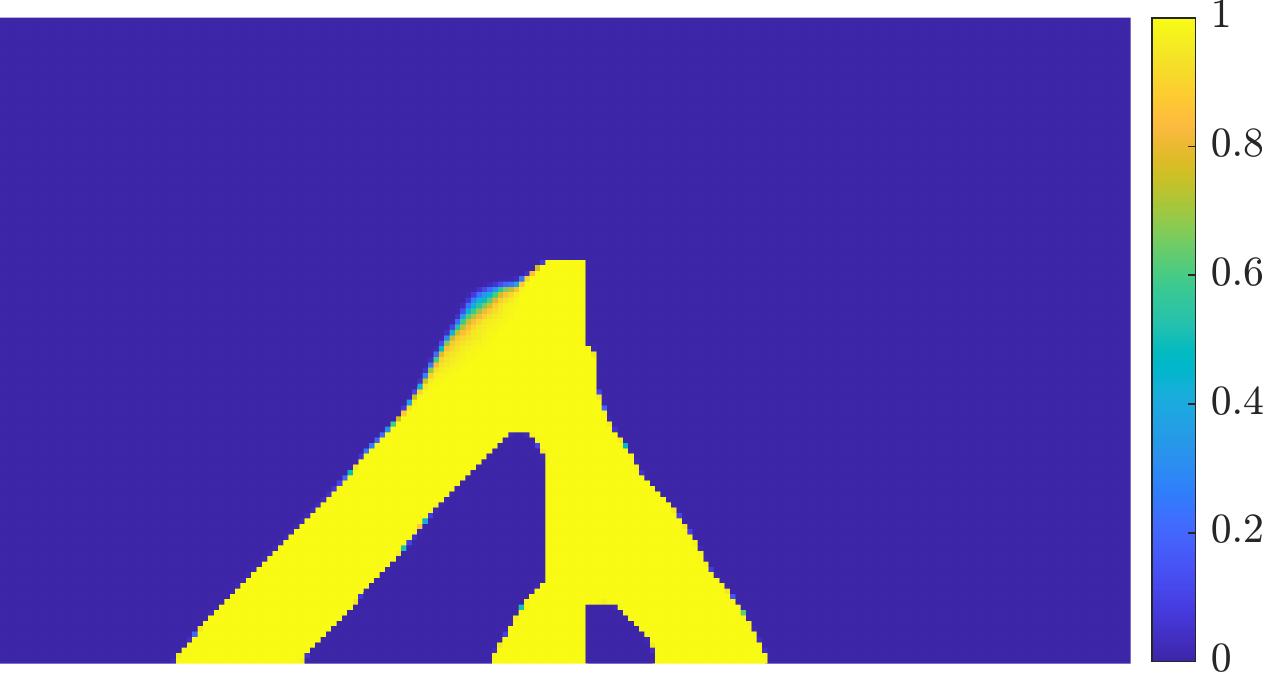}}
	\hspace{2pt}
	\subfloat[$p_\alpha = 6e-8, \ f = 0.00018026,$ \\ \vspace{1pt} $DM = 0.40 \%$.]{\includegraphics[width=\mywdth]{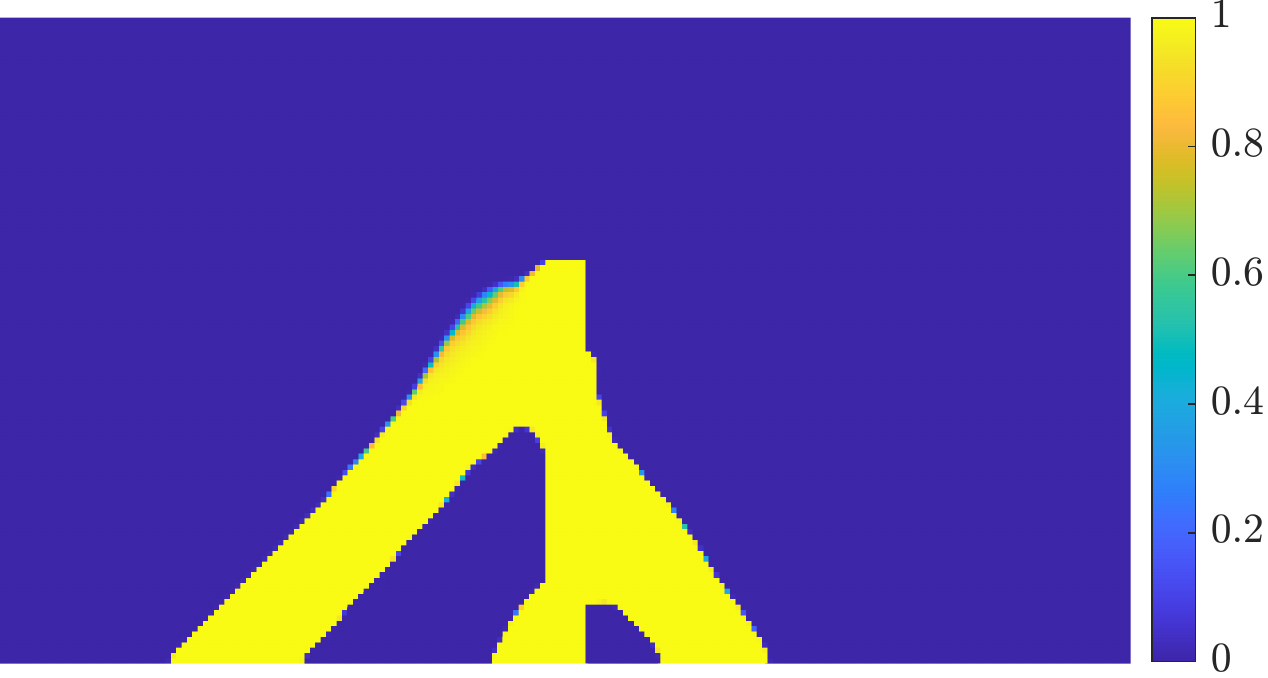}}
	\hspace{2pt}
	\subfloat[$p_\alpha = 5e-8, \ f = 0.00018077,$ \\ \vspace{1pt} $DM = 0.33 \%$. \label{ffig:9o}]{\includegraphics[width=\mywdth]{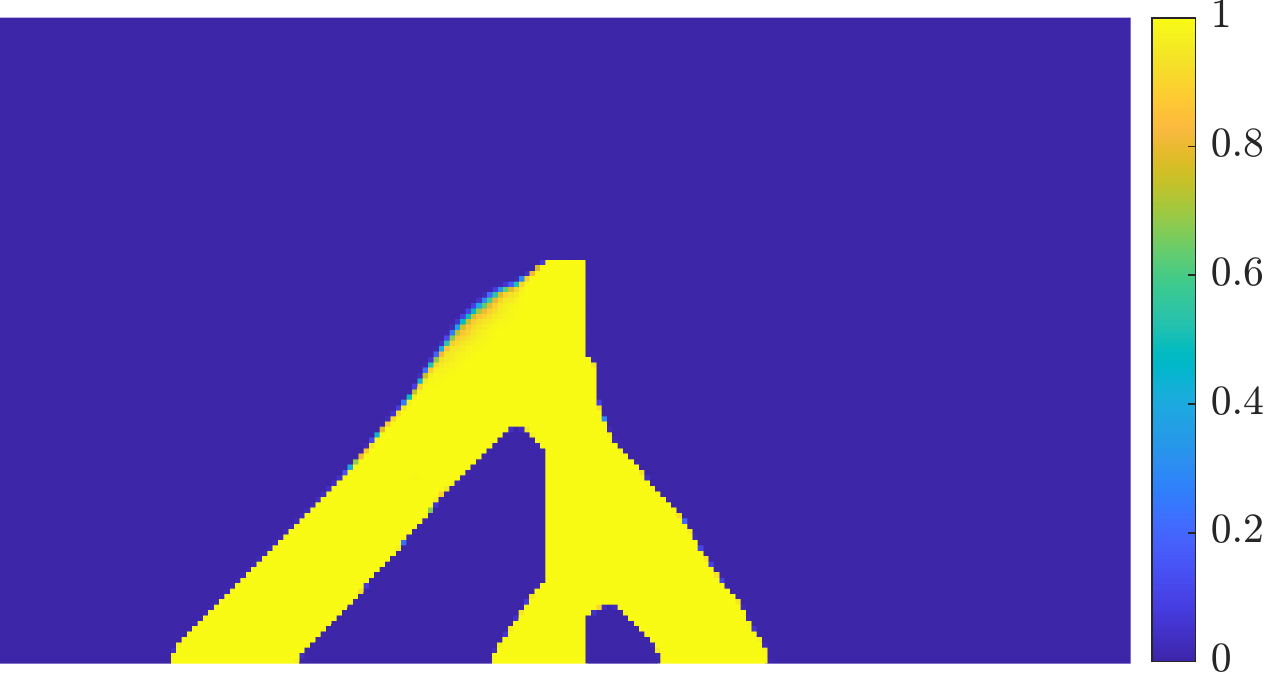}}
	\caption{Effect of $p_\alpha$ on optimized designs for the case $Re = 10$ with total stress coupling. At iterations \{21, 41, 61, 81\}, $p_E$ is set to \{1.5, 2.0, 3.0, 4.0\} while $p_\mathit{\Upsilon}$ is set to \{1.25, 1.50, 2.00, 2.50\}, i.e. $\delta_{E|\mathit{\Upsilon}} = 2$.}
	\label{fig:effect_of_p_alpha_re10_ttl}
\end{figure}

To further elaborate on the effect of $p_\alpha$ on the optimized designs, we solved the case of $Re = 10$ using both pressure and total stress coupling for a range of $p_\alpha$ values. A value of $\delta_{E|\mathit{\Upsilon}} = 2$ is used such that at iterations \{21, 41, 61, 81\}, $p_E$ is set to \{1.5, 2.0, 3.0, 4.0\} and $p_\mathit{\Upsilon}$ is set to \{1.25, 1.50, 2.00, 2.50\}. Some interesting conclusions can be drawn on the effect of $p_\alpha$:
\begin{enumerate}[a.]
	\item One important remark we noted in this work and in a previous work \cite{abdelhamid2023calculation}, is that the pressure is much more sensitive to the selection of $p_\alpha$ than the velocity. Meaning that a higher $p_\alpha$ is needed to get an accurate pressure field than is needed to get an accurate velocity field. By ``accurate" we mean in comparison to a segregated domain, body-fitted formulation approach.
	\item Large $p_\alpha$ means that the pressure is ``over-sensitive" to intermediate densities. Hence, the system tends to ``abuse" this sensitivity and build intermediate density ``walls" that redirect the flow and diffuse the pressure around the main structure. Thus reducing the objective function as the structure experiences a non-physically lower force than it should. This can be seen in Figs. \ref{fig:8a} to \ref{fig:8e} for the pressure coupling scenario and similarly in Figs. \ref{fig:9a} to \ref{fig:9i} for the total stress coupling scenario.
	\item Small $p_\alpha$ means the pressure is ``under-sensitive" to intermediate densities. The system then ``abuses" this insensitivity to pressure and extend thin structural members further to the left in a way that better supports the non-design column, which is still experiencing the correct fluidic forces. These thin structural members are ``immune" to pressure loads due to the small $p_\alpha$ selected thus rendering the design non-physical and wouldn't pass a thresholding test. This phenomenon can be observed in the pressure coupling scenario in Figs. \ref{fig:8l} to \ref{fig:8o}.
\end{enumerate}

At this point, it makes sense to formally define the appropriate $p_\alpha$ for a certain case as the value that doesn't result in non-physical, intermediate density walls and at the same time resolves the correct pressure field. Thus taking pressure loading into account during the design process and potentially passing a thresholding test. Building upon the above three remarks, it seems that the appropriate $p_\alpha$ for pressure coupling tends to be relatively higher than the appropriate $p_\alpha$ for total stress coupling. That is why the total stress coupling cases tend to display the same behavior as the pressure coupling cases towards changing $p_\alpha$, albeit at a lower $p_\alpha$. Examples of this observation are:
\begin{enumerate}[a.]
	\item The total stress coupling cases in Fig. \ref{fig:pe_vs_pu_2.5e-6_ttl} still show intermediate density walls at $\delta_{E|\mathit{\Upsilon}} = 1 \ \& \ 1.5$ in contrast to the pressure coupling cases (Fig. \ref{fig:pe_vs_pu_2.5e-6_prs}) not showing intermediate density walls even at $\delta_{E|\mathit{\Upsilon}} = 1$ for the same $p_\alpha$.
	\item The same observation is more pronounced in Fig. \ref{fig:effect_of_p_alpha_re10_prs} where the pressure coupling cases stop demonstrating intermediate density walls at $p_\alpha = 5e-7$ while the total stress coupling cases (Fig. \ref{fig:effect_of_p_alpha_re10_ttl}) still demonstrate intermediate density walls down till roughly $p_\alpha = 6e-8$.
\end{enumerate}

Based on these observations, we recommend - for calculating the appropriate $p_\alpha$ for a certain problem - to start from a high enough $p_\alpha$ that results in intermediate density walls, then reduce its value in steps until they disappear without any further reduction. The pressure field should be also monitored visually to ensure it can still ``see" the structure at the lower $p_\alpha$ values.

\subsection{Effect of Changing the Reynolds Number on the Optimized Designs}
Changing the Reynolds number seems to have a similar effect to what Lundgaard et al. \citep[p.~978]{Lundgaard2018} described. Increasing $Re$ tends to shift more material towards the upstream half to better strengthen the structural member on this side which faces increasing forces from the fluid flow. The same phenomenon is observed for both pressure and total stress coupling scenarios.

To check the suitability of the optimized designs for their respective $Re$, we performed cross-check calculations of each design for both $Re$ values. For the pressure coupling designs, both designs passed the cross-check as each appears to be suitable for their respective $Re$ as shown in table \ref{tab:prsr_crsschck}. However, for the total stress coupling designs, the design selected for $Re = 10$ failed its cross-check calculations. It seems that more trial-and-error solutions are needed to get the appropriate interpolation parameters for the increased $Re$. Probably a higher $\delta_{E|\mathit{\Upsilon}}$ would enable the use of a higher $p_\alpha$ without the appearance of intermediate density walls.

\begin{table}[t!]
	\centering
	\captionsetup{font+={color=black}}
	\caption{Cross-check for the optimized designs using pressure coupling against $Re$.}
	\label{tab:prsr_crsschck}
	\color{black}
	\begin{tabular}{clcc} \toprule
		& & \multicolumn{2}{c}{Obj. fun. analyzed for} \\
		Fig. & Designed for & \multicolumn{1}{c}{$Re = 1$} & \multicolumn{1}{c}{$Re = 10$} \\ \midrule
		\ref{ffig:5e} & $Re = 1$ & $\cellcolor{yellow} 1.16537e-2$ & $1.40134e-4$ \\ % TR 166.
		\ref{ffig:8k} & $Re = 10$ & $1.39845e-2$ & $\cellcolor{yellow} 1.05849e-4$ \\ % TR 212.
		\bottomrule
	\end{tabular}
\end{table}

\begin{table}[t!]
	\centering
	\captionsetup{font+={color=black}}
	\caption{Cross-check for the optimized designs using total stress coupling against $Re$.}
	\label{tab:ttlstrs_crsschck}
	\color{black}
	\begin{tabular}{clcc} \toprule
		& \multicolumn{2}{c}{Obj. fun. analyzed for} \\
		Fig. & Designed for & \multicolumn{1}{c}{$Re = 1$} & \multicolumn{1}{c}{$Re = 10$} \\ \midrule
		\ref{ffig:7c} & $Re = 1$ & $\cellcolor{yellow} 1.62095e-2$ & $\cellcolor{yellow} 1.76743e-4$ \\ % TR 168.
		\ref{ffig:9o} & $Re = 10$ & $1.76084e-2$ & $1.80766e-4$ \\ % TR 102.
		\bottomrule
	\end{tabular}
\end{table}

\begin{figure}[t!]
	\centering
	\captionsetup[subfigure]{justification=centering}
	\captionsetup{font+={color=black}}
	\captionsetup[subfigure]{font+={color=black}}
	\subfloat[Pressure coupling, $f = 0.00010301$, \\ \vspace{1pt} $DM = 0.39\%$. \label{ffig:prs_re10_vel}]{\includegraphics[width=0.328\textwidth]{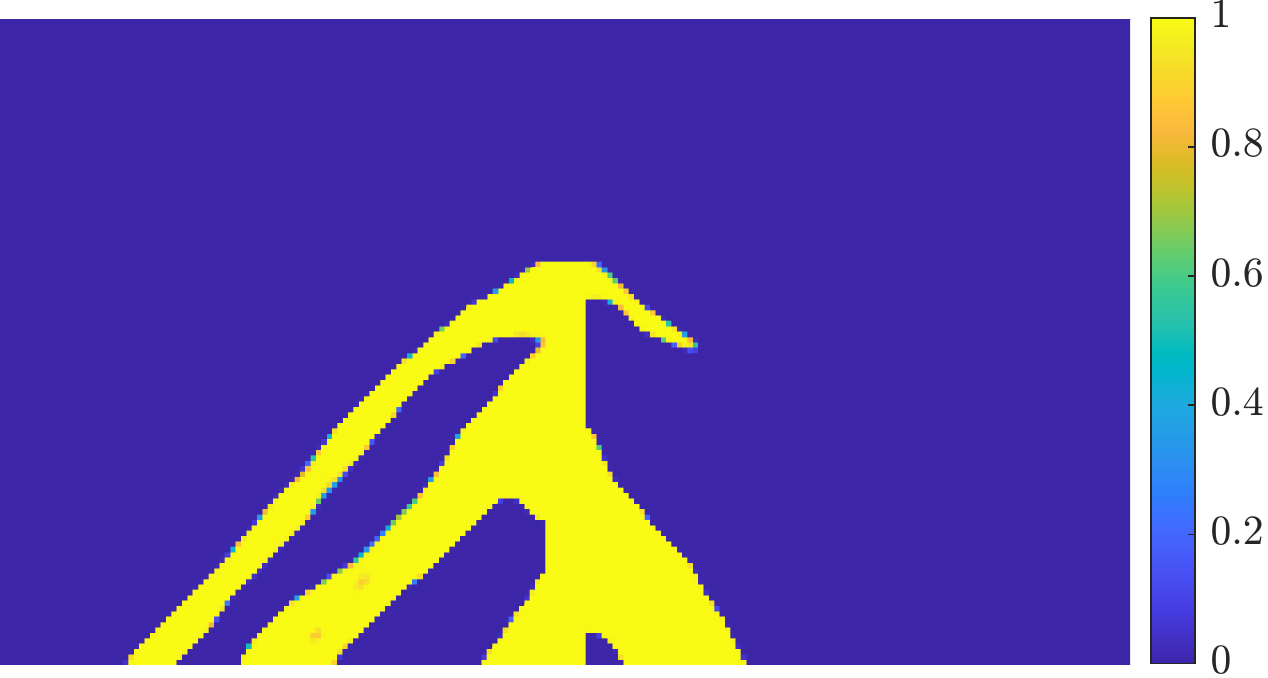}}
	\hspace{2pt}
	\subfloat[Total stress coupling, $f = 0.00018018$, \\ \vspace{1pt} $DM = 0.34\%$. \label{ffig:ttl_re10_vel}]{\includegraphics[width=0.328\textwidth]{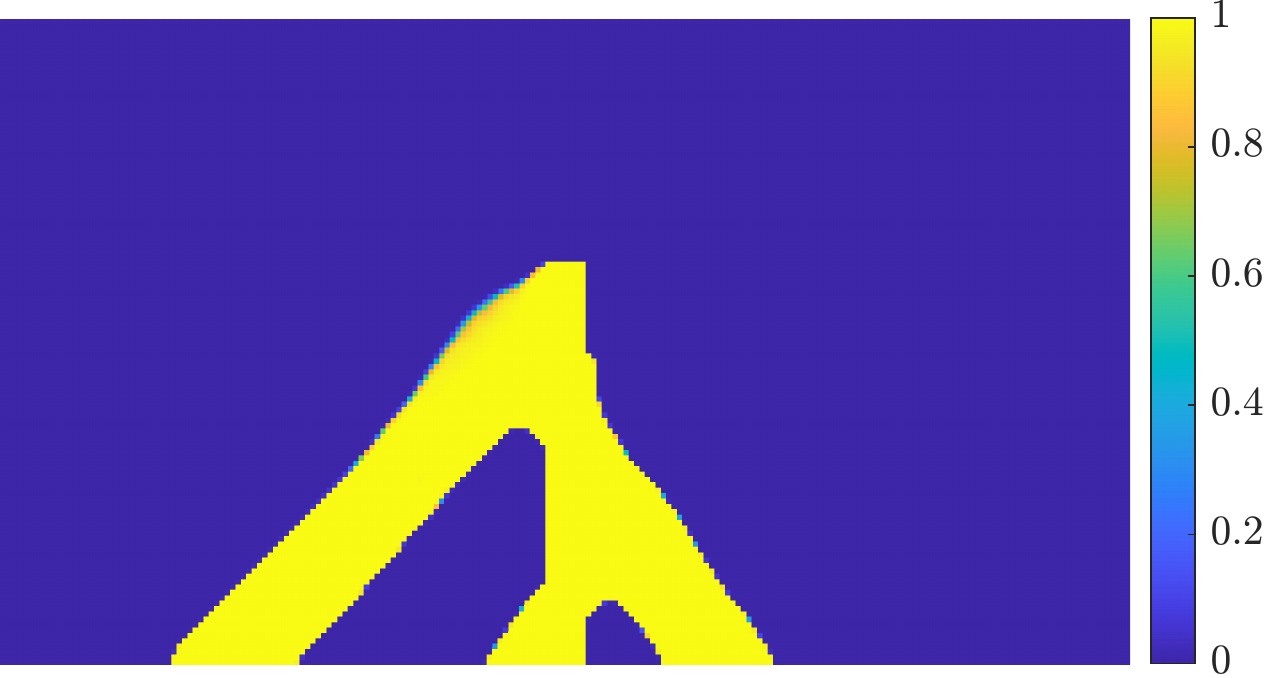}}
	\caption{Effect of using $v_\text{max}$ to change the Reynolds number instead of $\mu$.}
\end{figure}

In a previous work \cite{abdelhamid2023calculation}, we studied the effect of $Re$ on the appropriate Brinkman penalization maximum limit $\alpha_\text{max}$. This study was performed by looking at the maximum velocity in the solid regions, which made the investigation possible by looking at the individual terms of the Brinkman-penalized Navier-Stokes equations. We noted that the effect of $Re$ cannot be isolated from its constituents (i.e.  viscosity, velocity, density, and length); hence how $Re$ is changed was concluded to be more significant than the value of $Re$ itself. It is not straightforward to extend these conclusions to areas of intermediate and low density due to the strong non-linearity of the fluid flow equations, but we believe it's worth investigating in future work. To check this hypothesis, we solved two examples in which we used $v_\text{max} = 10$ m/s and $\mu$ = 1 Pa$\cdot$s to get $Re = 10$. This is in contrast to the previous $Re = 10$ designs (Figs. \ref{fig:effect_of_p_alpha_re10_prs} and \ref{fig:effect_of_p_alpha_re10_ttl}) in which we used $v_\text{max}$ = 1 m/s and $\mu = 0.1$ Pa$\cdot$s. We used $p_\alpha = 5e-7$, $\delta_{E|\mathit{\Upsilon}} = 2$, and 0.01 as a scaling factor for the objective function and its derivatives. Interestingly, the pressure coupling design (Fig. \ref{ffig:prs_re10_vel}) bears strong resemblance to the pressure coupling design obtained for $p_\alpha = 5e-8$ while using $\mu$ to change $Re$ (Fig. \ref{fig:8o}). Similarly, the total stress coupling design (Fig. \ref{ffig:ttl_re10_vel}) bears strong resemblance to the total stress coupling design obtained for $p_\alpha = 5e-8$ while using $\mu$ to change $Re$ (Fig. \ref{ffig:9o}). This confirms our hypothesis that $Re$ on its own is not enough to completely characterize fluid flow behavior in FSI problems.

\begin{figure}[b!]
	\centering
	\captionsetup{font+={color=black}}
	\captionsetup[subfigure]{font+={color=black}}
	\subfloat[Velocity (m/s).]{\includegraphics[width=0.49\textwidth]{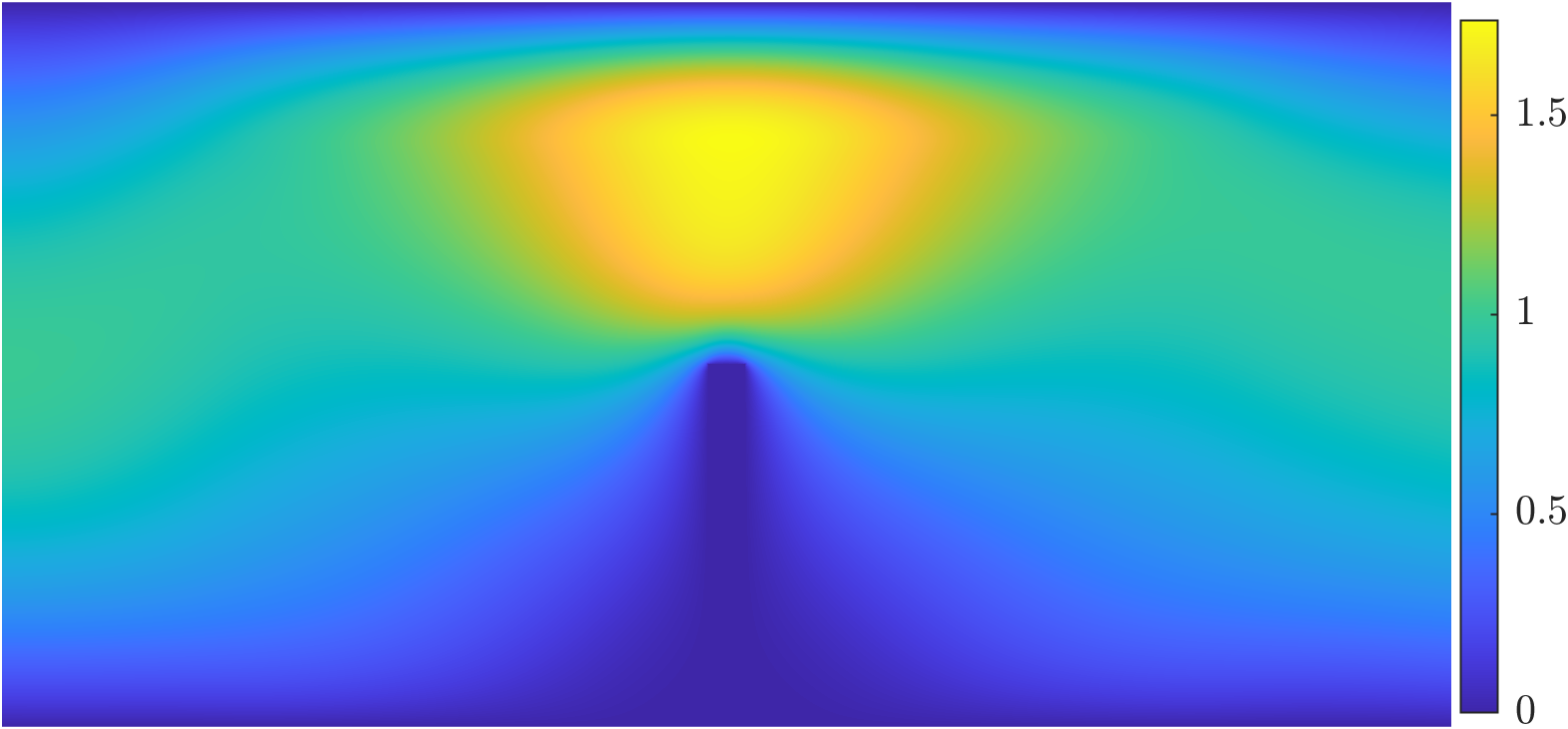}} 
	\hspace{2pt}
	\subfloat[Pressure (Pa).]{\includegraphics[width=0.49\textwidth]{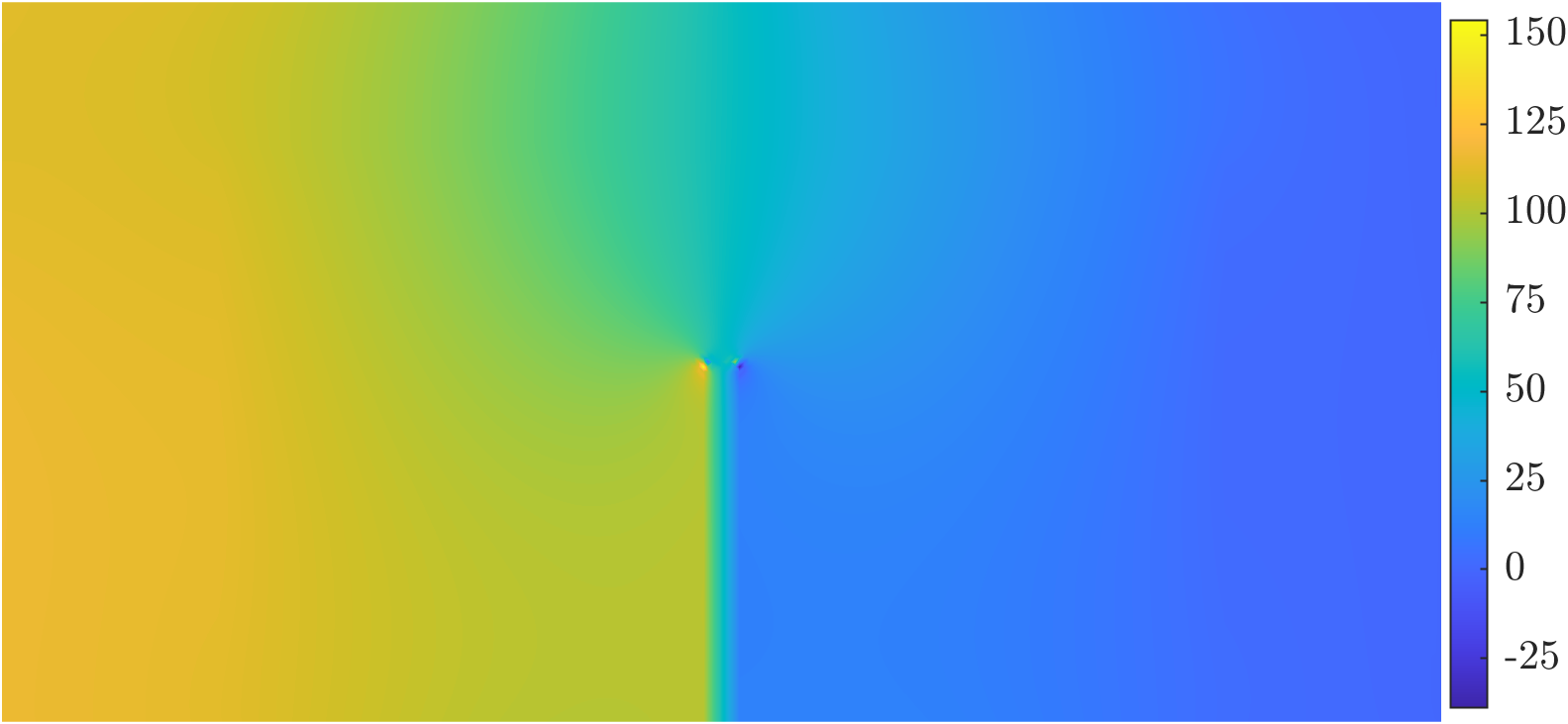}}	
	\caption{State fields at the first iteration for $Re = 1$ with $p_\alpha = 2.5e-6$.}
	\label{fig:state_fields_re0.04_initial}
\end{figure}

\begin{figure}[t!]
	\centering
	\captionsetup{font+={color=black}}
	\captionsetup[subfigure]{font+={color=black}}
	\subfloat[Pressure coupling.]{\includegraphics[width=0.49\textwidth]{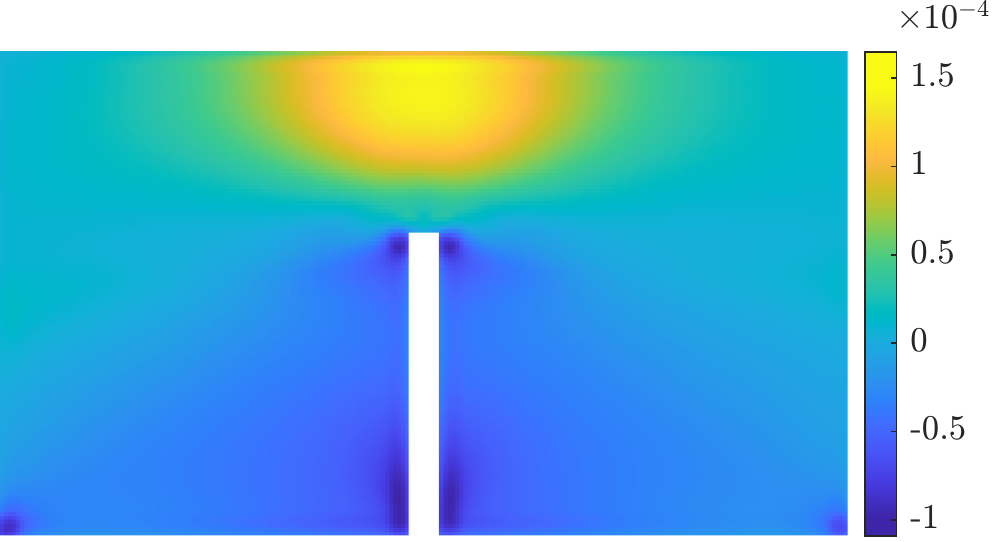}}
	\hspace{2pt}
	\subfloat[Total stress coupling.]{\includegraphics[width=0.49\textwidth]{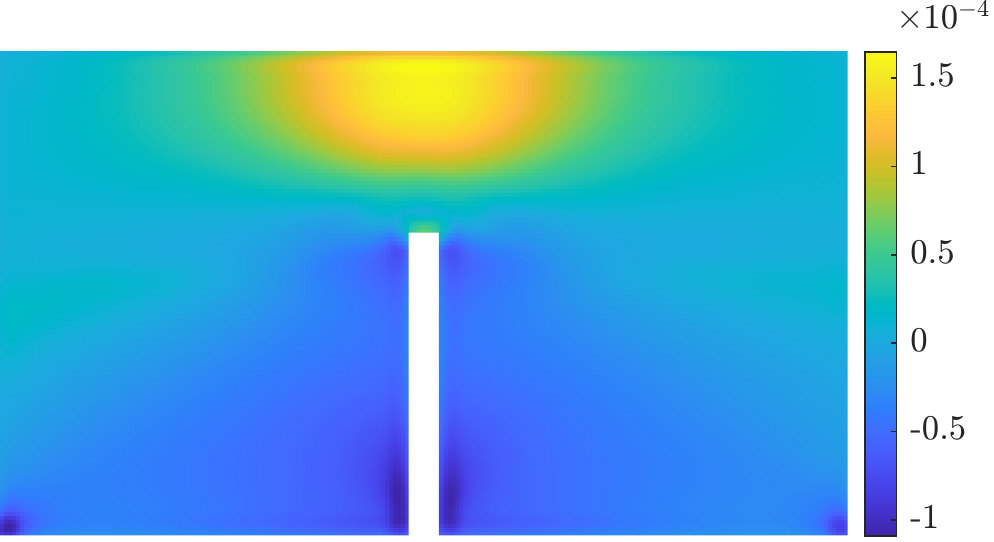}}
	\caption{Sensitivity of the column in a channel test problem at the second iteration for $Re = 1$ with $p_\alpha = 2.5e-6$.}
	\label{fig:snst_yoon_re0.04}
\end{figure}

\subsection{Effect of Pressure Spikes}

In the existence of external sharp corners, pressure spikes - in the form of extreme high and low localized values - are prone to appear. For instance, Fig. \ref{fig:state_fields_re0.04_initial} shows the velocity and pressure fields at the first iteration of the $Re = 1$ case, where these pressure spikes can be observed at the top sharp corners of the non-design column. Note that the type of force coupling used doesn't affect the resolved velocity and pressure fields, but rather how the fluidic forces are extracted after. Hence, the velocity and pressure fields of the initial design shown in Fig. \ref{fig:state_fields_re0.04_initial} are exactly similar for either type of force coupling. The sensitivities, on the other hand, are indeed affected by the type of force coupling.

In the pressure coupling cases, the sensitivities are greatly affected by these pressure spikes, which led the optimizer to place material in this vicinity to alleviate their effect. However, in the total stress coupling cases, the sensitivities were not affected as much by these pressure spikes, and the available material was used somewhere else to better support the design. This is evident in Fig. \ref{fig:snst_yoon_re0.04}, which shows the sensitivities for $Re = 1$ at the second iteration for the pressure vs. total stress coupling cases.

A pressure vs. total stress cross-check is performed in table \ref{tab:force_crsschck} by analyzing each design for pressure as well as total stress coupling (the yellow highlight indicates the best performance). It is clear that each design performs well for its intended purpose and that excluding the viscous stresses may alter the optimized result depending on the flow characteristics.

\begin{table}[t]
	\centering
	\captionsetup{font+={color=black}}
	\caption{Cross-check for the optimized designs against the type of force coupling.}
	\label{tab:force_crsschck}
	\color{black}
	\begin{tabular}{cclccrr} \toprule
		& & & & & \multicolumn{2}{c}{Obj. fun. analyzed for} \\
		Fig. & $Re$ & Designed for & $p_\alpha$ & Discreteness & \multicolumn{1}{c}{Pressure} & \multicolumn{1}{c}{Total Stress} \\ \midrule
		\ref{ffig:5e} & 1 & Pressure & $5.25e-6$ & 0.30\% & \cellcolor{yellow} $1.16537e-2$ & $2.03826e-2$  \\ % TR 166.
		\ref{ffig:7c} & 1 & Total Stress & $2.50e-6$ & 0.19\% & $1.24263e-2$ & \cellcolor{yellow} $1.62095e-2$ \\ % TR 168.
		& & & & \\[-5pt]
		\ref{ffig:8k} & 10 & Pressure & $9.00e-8$ & 0.37\% & \cellcolor{yellow} $1.05849e-4$ & $2.48147e-4$ \\ % TR 212.
		\ref{ffig:9o} & 10 & Total Stress & $5.00e-8$ & 0.33\% & $1.21663e-4$ & \cellcolor{yellow} $1.80766e-4$ \\ % TR214.
		% PRB.ofsc is e+0, e+2 respectively.
		\bottomrule
	\end{tabular}
\end{table}

\section{Conclusions}
\label{sec:cnclsns}
In this work, we extended the force coupling in the topology optimization of fluid-structure interaction problems from pressure stresses to total stresses through the inclusion of viscous stress components, thus increasing the fidelity and versatility of TOFSI as a design tool. To overcome the discontinuities in the derivatives of the velocity across Lagrangian elemental boundaries, we utilized the superconvergent patch recovery technique by Zienkiewicz and Zhu to smoothen the velocity derivatives at the elemental boundaries. Sensitivity analysis equations of the implemented technique are derived and a sensitivity verification check is performed using the complex-step derivative approximation method. The key conclusions from this work are:
\begin{itemize}
	\item TOFSI problems using density-based methods are highly sensitive to the selection of the interpolation parameters in general and the Brinkman penalization interpolation parameter $p_\alpha$ in specific.
	\item Depending on the flow characteristics, pressure vs. total stress force coupling may show a differentiation in the optimized designs, where each design works well for its intended purpose. Pressure spikes at sharp corners are a clear example.
	\item With efficient coding practices and in the absence of mesh deformations, the increase in computational costs when using total stress coupling as opposed to pressure coupling in TOFSI can be kept below 10\%, which is a negligible factor when deciding the use of either coupling methods. With mesh deformations taken into account, the increase in computational costs is expected to be higher.
	\item The effect of the Reynolds number $Re$ on its own cannot be isolated from its constituents. We have demonstrated with examples that controlling $Re$ using the dynamic viscosity $\mu$ is entirely different from controlling $Re$ using the characteristic velocity $v_\text{max}$. We conjecture that the main reason behind this issue is that the forces experienced by the structure are not only related to $Re$ but rather to its constituents (i.e. characteristic velocity, characteristic length, density, and dynamic viscosity).
\end{itemize}

\section*{Acknowledgment}
The first author would like to thank Dr. Joe Alexandersen, University of South Denmark, for a lot of helpful advice and many fruitful discussions. The authors are grateful to Prof. Krister Svanberg, KTH Royal Institute of Technology in Stockholm, for providing the MATLAB implementation of the Method of Moving Asymptotes. This research was enabled in part by support provided by BC DRI Group and the Digital Research Alliance of Canada (alliancecan.ca).

\bibliography{library.bib}
\bibliographystyle{elsarticle-num-names}

\pagebreak

\appendix

\section{Numerical Verification of the Brinkman-Penalized Fluid-Structure Interaction Formulation}

In this appendix, we verify numerically the Brinkman-penalized, unified domain FSI formulation through a comparison with a body-fitted, segregated domain FSI formulation. To provide an additional layer of verification of our proprietary code, the body-fitted results are extracted from COMSOL Multiphysics \cite{comsol6}.

\begin{figure}[t]
	\centering
	\begingroup
	\captionsetup[subfigure]{width=0.45\textwidth}
	\subfloat[Velocity field with streamlines extracted from COMSOL Multiphysics.]{\includegraphics[width=0.49\textwidth]{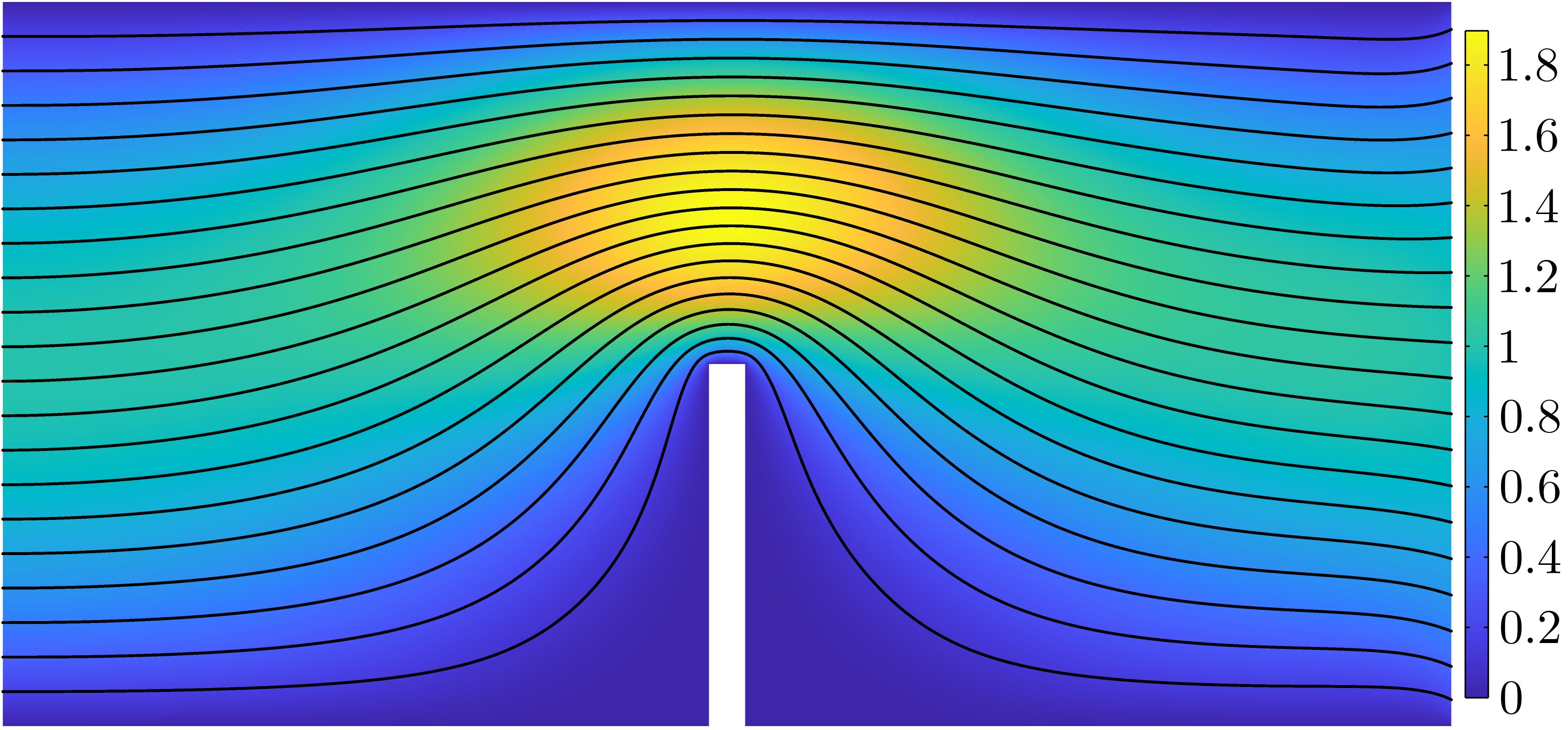}}
	\subfloat[Velocity field with streamlines extracted from our proprietary code.]{\includegraphics[width=0.49\textwidth]{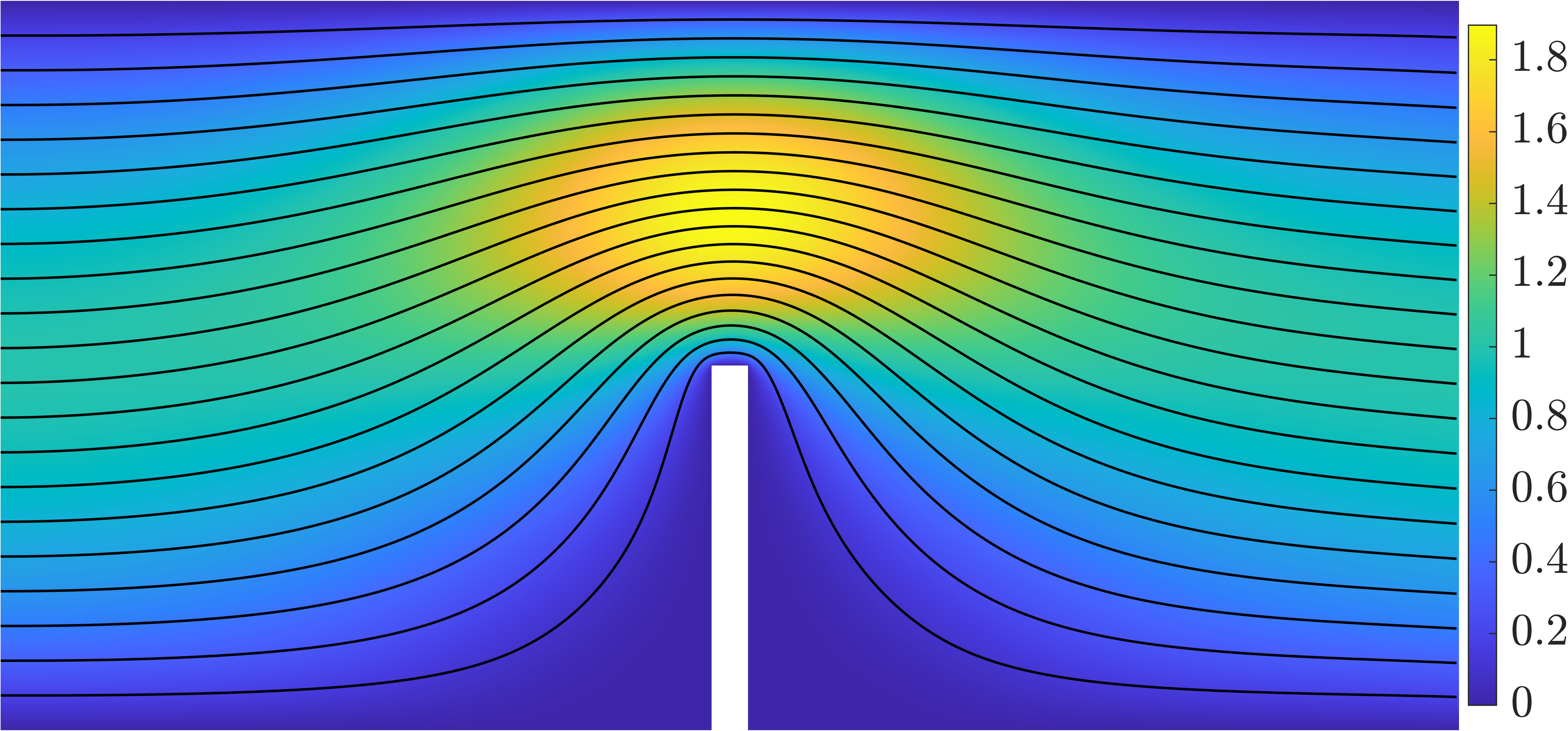}} \\
	\subfloat[Pressure field with contour lines extracted from COMSOL Multiphysics.]{\includegraphics[width=0.49\textwidth]{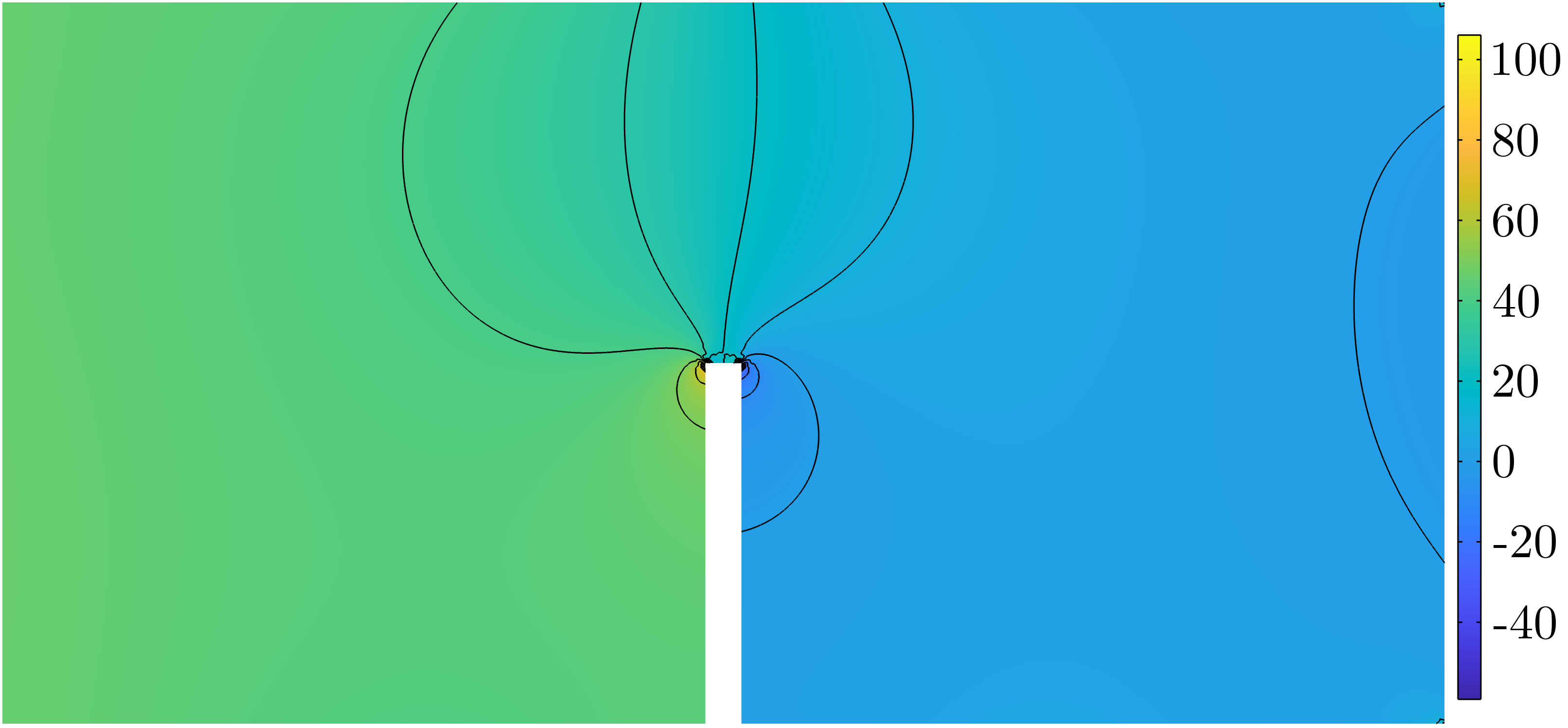}}
	\subfloat[Pressure field with contour lines extracted from our proprietary code.]{\includegraphics[width=0.49\textwidth]{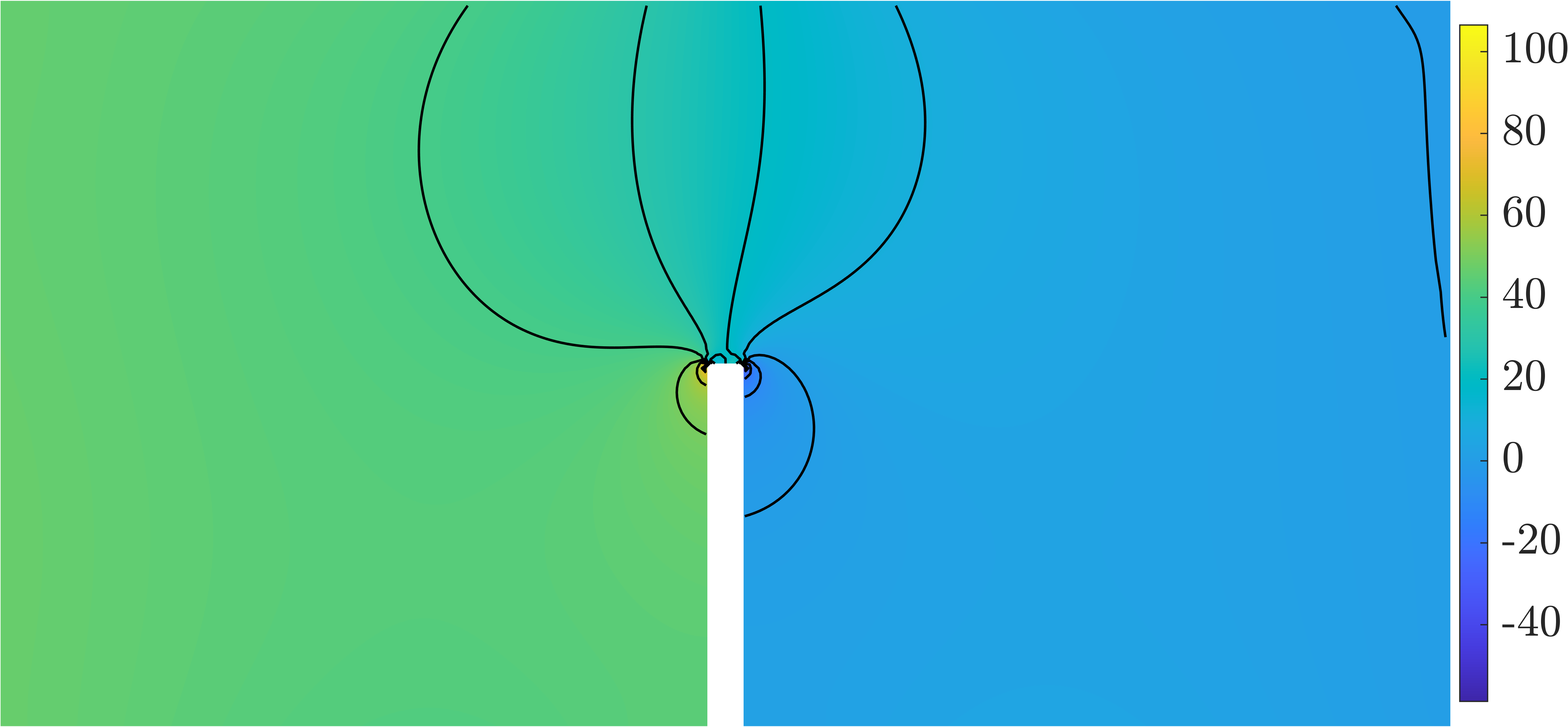}}
	\endgroup
	\caption{Numerical verification of the Brinkman-penalized FSI formulation for the initial design of the column in a channel problem. Velocity fields with streamlines and pressure fields and contour lines are shown for comparison.}
	\label{fig:crossvldt_initdsgn}
\end{figure}

The \textbf{first case} is an initial design of the column in a channel problem described in Sec. \ref{sec:tst_prblm} assuming pressure coupling. The design domain is set to zero (i.e. $\rho = 0$ in $\mathit{\Omega_d}$) while keeping the non-design column as solid (i.e. $\rho = 1$ in $\mathit{\Omega_{nd}}$). To maintain the infinitesimal deformation assumption at this reduced material state, the elastic modulus is raised to $1e+7$ Pa. In COMSOL Multiphysics, we implement the same mesh and finite element type as our proprietary code, and the fluid-structure coupling is set to pressure forces only through modifying the weak expression under the fluid-structure interaction node. The results to be compared are the maximum and minimum $x$ and $y$ velocities, the maximum and minimum pressures, the total compliance of the structure, and the $x$ and $y$ displacements of a point `A' located at the top left corner of the non-design column such that $\{x_A, y_A\} = \{0.975 \ \mathrm{m}, 0.5 \ \mathrm{m}\}$. We also compare the fluid velocity field with streamlines and the pressure field with contour lines. As can be seen in Table \ref{tab:crossvldt_initdsgn}, there is good agreement between the unified and segregated formulations for the parameters of interest. Similarly, the velocity and pressure fields show good agreement as shown in Fig. \ref{fig:crossvldt_initdsgn}.

As for the total stress coupling, the forces experienced by the structure on the fluid-structure interface in the segregated domain formulation is from one side only, hence the discontinuity issue is not as severe as in the unified domain formulation implemented in this work. To the best of our knowledge, commercial software packages typically don't treat these discontinuities by default. In addition, discontinuity treatment using different techniques are not directly comparable. Hence, a similar verification with the total stress coupling case doesn't offer any value.

\begin{table}[b!]
	\centering
	\captionsetup{font+={color=black}}
	\captionsetup{width=0.78\textwidth}
	\caption{Numerical verification of the Brinkman-penalized FSI formulation for the initial design of the column in a channel problem. Some parameters of interest are selected for comparison.}
	\label{tab:crossvldt_initdsgn}
	\begingroup
	\color{black}
	\begin{tabular}{wc{1in} wc{1.2in} wc{1.2in} wc{1in}} \toprule
		Parameter & Segregated Formulation & Unified Formulation & Normalized Error \% \\ \midrule
		max($\mathbf{v_1}$) (m/s) & $+1.89732955e+0$ & $+1.89729912e+0$ & $+1.60e-3\%$ \\
		& & & \\[-5pt]
		min($\mathbf{v_1}$) (m/s) & $-2.90659936e-3$ & $-2.87661716e-3$ & $+1.03e+0\%$ \\
		& & & \\[-5pt]
		max($\mathbf{v_2}$) (m/s) & $+7.10320533e-1$ & $+7.10314492e-1$ & $+8.50e-4\%$ \\
		& & & \\[-5pt]
		min($\mathbf{v_2}$) (m/s) & $-6.64991818e-1$ & $-6.65954950e-1$ & $-1.45e-1\%$ \\
		& & & \\[-5pt]
		max($\mathbf{p}$) (Pa) & $+106.106710e+0$ & $+106.487052e+0$ & $-3.58e-1\%$ \\
		& & & \\[-5pt]
		min($\mathbf{p}$) (Pa) & $-59.1725096e+0$ & $-58.7279842e+0$ & $+7.51e-1\%$ \\
		& & & \\[-5pt]
		$\mathbf{\hat{U}} ^T \ \mathbf{K}^s \ \mathbf{\hat{U}}$ (J) & $+4.88811847e-2$ & $+4.88899052e-2$ & $-1.78e-2\%$ \\
		& & & \\[-5pt]
		$u_{1|A}$ (m) & $+4.13584942e-3$ & $+4.13607818e-3$ & $-5.53e-3\%$ \\
		& & & \\[-5pt]
		$u_{2|A}$ (m) & $+2.79011498e-4$ & $+2.79012811e-4$ & $-4.71e-4\%$ \\
		\bottomrule
	\end{tabular}
	\endgroup
\end{table}

The \textbf{second case} is an optimized design of the column in a channel test problem using pressure coupling at $Re = 1$; namely Fig. \ref{ffig:6c}. To generate the geometry needed for COMSOL Multiphysics, the optimized design is traced manually in an external CAD software then imported into COMSOL Multiphysics. The fluid-structure coupling is set to pressure similar to the first case. An unavoidable byproduct of using the unified domain formulation, is the uncontrolled pressure inside the enclosed cavities. In the segregated domain formulation, this internal pressure doesn't typically occur. To generate the same boundary conditions in COMSOL Multiphysics, the fluid flow is first solved using the Brinkman-penalized Navier-Stokes to calculate these internal pressure values. This is performed by appending the fluid flow node with a volume force (i.e. $-\alpha_\text{max} \cdot v_1$ in $x$ and $-\alpha_\text{max} \cdot v_2$ in $y$). Once these internal pressures are calculated in COMSOL, the volume force node is disabled and the entire problem is solved as a fluid-structure interaction problem in the segregated domain formulation where the calculated pressures are applied as boundary loads at the enclosed cavities. The calculated compliance using the segregated formulation is 0.01251 J against 0.01137 J for the unified domain formulation resulting in an error of 9.08\%. There are mainly two reasons for this discrepancy; \textbf{(i)} the effect of the smooth interfaces implemented in COMSOL vs the original jagged interface in the unified domain solution, and \textbf{(ii)} the effect of some intermediate density elements on the outside interface that strongly affect the pressure field resolved. We believe both effects can be alleviated if some sort of adaptive mesh refinement is implemented during the optimization process, cf. \cite{Lambe2018a, SalazardeTroya2018}.

\end{document}